\documentclass[12pt]{article}
\usepackage{amssymb}
\usepackage{latexsym}
\usepackage{amsthm}

\newtheorem{definition}{Definition}[section]
\newtheorem{theorem}[definition]{Theorem}
\newtheorem{lemma}[definition]{Lemma}
\newtheorem{corollary}[definition]{Corollary}
\newtheorem{remark}[definition]{Remark}
\newtheorem{example}[definition]{Example}

\newtheorem{note}[definition]{Note}

\newtheorem{proposition}[definition]{Proposition}

\newtheorem{notation}[definition]{Notation}

\newcommand{\beast}{\begin{eqnarray*}}
\newcommand{\eeast}{\end{eqnarray*}}

\renewcommand{\b}[1]{\langle #1 \rangle}

\def\K{\mathbb K}

\def\K{\mathbb K}
\def\fld{\mathbb K}
\def\Mdf{{\rm Mat}_{d+1}(\K)}

\begin{document}

\title{\bf Thin Hessenberg Pairs and Double Vandermonde Matrices}
\author{Ali Godjali}


\date{July 26, 2011}

\maketitle

\begin{abstract}
\noindent 
A square matrix is called {\it Hessenberg} whenever each entry below the subdiagonal is zero and each entry on the subdiagonal is nonzero. Let $V$ denote a nonzero finite-dimensional vector space over a field $\fld$. We consider an ordered pair of linear transformations $A: V
\rightarrow V$ and $A^*: V \rightarrow V$ which satisfy both (i), (ii) below.
\begin{enumerate}
\item[\rm (i)] There exists a basis for $V$ with respect to which the matrix representing $A$ is Hessenberg and the matrix representing $A^*$ is diagonal. 

\item[\rm (ii)] There exists a basis for $V$ with respect to which the matrix representing $A$ is diagonal and the matrix representing $A^*$ is Hessenberg.  
\end{enumerate}

\noindent We call such a pair a {\it thin Hessenberg pair} (or {\it TH pair}). By the {\it diameter} of the pair we mean the dimension of $V$ minus one. There is an ``oriented'' version of a TH pair called a TH system. In this paper we investigate a connection between TH systems
and double Vandermonde matrices. We have two main results. For the first result we give a bijection between any two of the following three sets: 

\begin{itemize}
 \item
The set of isomorphism classes of TH systems over $\K$ of diameter $d$.
 \item
The set of normalized west-south Vandermonde systems in $\Mdf$. 
 \item
The set of parameter arrays over $\K$ of diameter $d$. 
\end{itemize}

\noindent For the second result we give a bijection between any two of the following five sets:

\begin{itemize}
 \item
The set of affine isomorphism classes of TH systems over $\K$ of diameter $d$.
 \item
The set of isomorphism classes of RTH systems over $\K$ of diameter $d$. 
 \item
The set of affine classes of normalized west-south Vandermonde systems in $\Mdf$.
 \item
The set of normalized west-south Vandermonde matrices in $\Mdf$. 
 \item
The set of reduced parameter arrays over $\K$ of diameter $d$.
\end{itemize}

\bigskip


\end{abstract}

\section{Introduction}
\label{sec:intro}
This paper is about a linear algebraic object called a thin Hessenberg pair \cite{thinhess}. 
To recall its definition, we will use the following term. A square matrix is called {\it Hessenberg} whenever each entry below the subdiagonal is zero and each entry on the subdiagonal is nonzero. Throughout the paper, $\fld$ will denote a field.

\begin{definition}
{\rm \cite[Definition 1.1]{thinhess}}
\label{def:thinhess} 
\rm 
Let $V$ denote a nonzero finite-dimensional vector space over $\fld$. By a \emph{thin Hessenberg pair} (or {\it TH pair}) on $V$, we mean an ordered pair of linear transformations $A:V\rightarrow V$ and $A^*:V\rightarrow V$ which satisfy both (i), (ii) below. 
\begin{enumerate}
\item[\rm (i)] 
There exists a basis for $V$ with respect to which the matrix representing $A$ is Hessenberg and the matrix representing $A^*$ is diagonal. 
\item[\rm (ii)]
There exists a basis for $V$ with respect to which the matrix representing $A$ is diagonal and the matrix representing $A^*$ is Hessenberg.
\end{enumerate}
We call $V$ the {\it underlying} vector space and say that $A, A^*$ is {\it over} $\fld$. By the {\it diameter} of $A, A^*$ we mean the dimension of $V$ minus one.
\end{definition}

\begin{note}  
\samepage
\rm
It is a common notational convention to use $A^*$ to represent the conjugate-transpose of $A$. We are not using this convention. In a TH pair $A, A^*$ the linear transformations $A$ and $A^*$ are arbitrary subject to (i), (ii) above.
\end{note}

A TH pair is a generalization of a Leonard pair \cite{LS99}. Roughly speaking, a Leonard pair is a pair of linear transformations as in Definition \ref{def:thinhess}, with the Hessenberg requirement replaced by an irreducible tridiagonal requirement. Leonard pairs have been extensively studied; for more information see \cite{madrid} and the references therein. 

In \cite{thinhess} we introduced the concept of a TH pair and began a systematic study of these objects. We now summarize the
content of \cite{thinhess}. In \cite[Definition 2.2]{thinhess} we introduced an ``oriented'' version of a TH pair called a TH system. A TH system is described as follows. Let $A, A^*$ denote a TH pair on $V$ of diameter $d$. By definition $A$ is diagonalizable. It turns out that each eigenspace of $A$ has dimension one \cite[Lemma 2.1]{thinhess}.
Therefore a basis from Definition \ref{def:thinhess}(ii) induces an ordering $\{ V_i \}_{i=0}^d$ of the eigenspaces of $A$. For $0 \leq i \leq d$, let $E_i$ denote the primitive idempotent of $A$ that corresponds to $V_i$. 
We call $\lbrace E_i\rbrace_{i=0}^d$ a standard ordering of the primitive idempotents of $A$. A standard ordering of the primitive idempotents of $A^*$ is defined similarly. A TH system is a TH pair $A, A^*$ together with a standard ordering of the primitive idempotents of $A$ and a standard ordering of the primitive idempotents of $A^*$. Let $(A;\lbrace E_i\rbrace_{i=0}^d;A^*;\lbrace E^*_i\rbrace_{i=0}^d)$ denote a TH system on $V$. In \cite{thinhess} we investigated six bases for $V$ with respect to which the matrices representing $A$ and $A^*$ are attractive. We displayed these matrices along with the transition matrices relating the bases. We classified the TH systems up to isomorphism.

In the present paper, we continue our study of TH pairs and TH systems. Our focus is on a connection between TH systems and double Vandermonde matrices. We establish two main results. These results have the following form. In the first result we display three sets and show any two are in bijection. In the second result we display five sets and show any two are in bijection. We now describe the first result. To do this we display the three sets and then discuss the meaning. The three sets are:

\begin{itemize}
 \item
The set of isomorphism classes of TH systems over $\K$ of diameter $d$.
 \item
The set of normalized west-south Vandermonde systems in $\Mdf$. 
 \item
The set of parameter arrays over $\K$ of diameter $d$. 
\end{itemize}

We now describe the above three sets in more detail. The first set is clear, so
consider the second set. For an indeterminate $\lambda$ let $\K[\lambda]$ denote the $\K$-algebra consisting of the polynomials in $\lambda$ that have all coefficients in $\K$. Let $\{ f_i \}_{i=0}^d$ denote a sequence of polynomials in $\K[\lambda]$. We say that $\{ f_i \}_{i=0}^d$ is graded whenever $f_0 = 1$ and $f_i$ has degree $i$ for $0 \leq i \leq d$.
By a normalized west-south Vandermonde system in $\Mdf$ we mean a sequence $(X, \lbrace \theta_i \rbrace_{i=0}^d, \lbrace \theta_i^* \rbrace_{i=0}^d)$ such that: (i) $X$ is a matrix in $\Mdf$; (ii) $\lbrace \theta_i \rbrace_{i=0}^d$ is a sequence of mutually distinct scalars in $\K$; (iii) $\lbrace \theta_i^* \rbrace_{i=0}^d$ is a sequence of mutually distinct scalars in $\K$; (iv)
there exists a graded sequence of polynomials $\{ f_i \}_{i=0}^d$ in $\K[\lambda]$ such that $X_{ij} = f_j(\theta_i)$ for $0 \leq i,j \leq d$; (v) there exists a graded sequence of polynomials $\{ f_i^* \}_{i=0}^d$ in $\K[\lambda]$ such that $X_{ij} = f_{d-i}^*(\theta_j^*)$ for $0 \leq i,j \leq d$. We now describe the third set. 
By a parameter array over $\fld$ of diameter $d$ we mean a sequence $(\lbrace \theta_i \rbrace_{i=0}^d, \lbrace \theta^*_i \rbrace_{i=0}^d, \lbrace \phi_i \rbrace_{i=1}^d)$ of scalars taken from $\fld$ such that: (i) $\lbrace \theta_i \rbrace_{i=0}^d$ are mutually distinct; (ii) $\lbrace \theta^*_i \rbrace_{i=0}^d$ are mutually distinct; (iii) $\lbrace \phi_i \rbrace_{i=1}^d$ are all nonzero.
We have now described the three sets. We now describe the bijections between these sets. We start by describing the bijection from the first set to the second set. Let $\Phi = (A;\lbrace E_i\rbrace_{i=0}^d;A^*;\lbrace E^*_i\rbrace_{i=0}^d)$ denote a TH system on $V$. Associated with $\Phi$ is a certain matrix $\mathcal P \in \Mdf$. This is the transition matrix from a basis in Definition \ref{def:thinhess}(ii) to a basis in Definition \ref{def:thinhess}(i), where the bases are normalized so that each entry in the leftmost column and the bottom row of $\mathcal P$ is $1$. For $0 \leq i \leq d$ let $\theta_i$ (resp. $\theta^*_i$) denote the eigenvalue of $A$ (resp. $A^*$) that corresponds to $E_i$ (resp. $E_i^*$). Our bijection sends the isomorphism class of $\Phi$ to $(\mathcal P, \lbrace \theta_i \rbrace_{i=0}^d, \lbrace \theta^*_i \rbrace_{i=0}^d)$. 
We now describe the bijection from the third set to the first set. Let $(\lbrace \theta_i \rbrace_{i=0}^d, \lbrace \theta^*_i \rbrace_{i=0}^d, \lbrace \phi_i \rbrace_{i=1}^d)$ denote a parameter array over $\K$ of diameter $d$. Let $A$ denote the lower bidiagonal matrix in $\Mdf$ with entries $A_{ii} = \theta_{d-i}$ for $0 \leq i \leq d$ and $A_{i,i-1} = \phi_i$ for $1 \leq i \leq d$.  
Let $A^*$ denote the upper bidiagonal matrix in $\Mdf$ with entries $A^*_{ii} = \theta_{i}^*$ for $0 \leq i \leq d$ and $A^*_{i-1,i} = 1$ for $1 \leq i \leq d$. Observe that $\lbrace \theta_i \rbrace_{i=0}^d$ (resp. $\lbrace \theta_i^* \rbrace_{i=0}^d$) is an ordering of the eigenvalues of $A$ (resp. $A^*$). For $0 \leq i \leq d$ let $E_i$ (resp. $E_i^*$) denote the primitive idempotent of $A$ (resp. $A^*$) that corresponds to $\theta_i$ (resp. $\theta_i^*$). We show that $\Phi = (A; \{E_i\}_{i=0}^d;A^*; \{E^*_i\}_{i=0}^d)$ is a TH system. 
Our bijection sends $(\lbrace \theta_i \rbrace_{i=0}^d, \lbrace \theta^*_i \rbrace_{i=0}^d, \lbrace \phi_i \rbrace_{i=1}^d)$ to the isomorphism class of $\Phi$. 

\smallskip

We now describe our second result, which is a variation on the first result. We mentioned above that the second result involves five sets. The five sets are:

\begin{itemize}
 \item
The set of affine isomorphism classes of TH systems over $\K$ of diameter $d$.
 \item
The set of isomorphism classes of RTH systems over $\K$ of diameter $d$. 
 \item
The set of affine classes of normalized west-south Vandermonde systems in $\Mdf$.
 \item
The set of normalized west-south Vandermonde matrices in $\Mdf$. 
 \item
The set of reduced parameter arrays over $\K$ of diameter $d$.
\end{itemize}

\noindent We now describe the above five sets in more detail. Throughout the description let $\alpha, \beta, \alpha^*, \beta^*$ denote scalars in $\K$ with $\alpha, \alpha^*$ nonzero. We now describe the first set. Let \\ $\Phi = (A;\lbrace E_i\rbrace_{i=0}^d;A^*;\lbrace E^*_i\rbrace_{i=0}^d)$ denote a TH system over $\K$. Observe that the sequence $(\alpha A + \beta I; \{E_i\}_{i=0}^d; \alpha^* A^* +\beta^* I; \{E^*_i\}_{i=0}^d)$ is a TH system over $\K$, said to be an affine transformation of $\Phi$. We now describe the second set. By an RTH system over $\K$ we mean the sequence $(\lbrace E_i\rbrace_{i=0}^d; \lbrace E^*_i\rbrace_{i=0}^d)$ induced by a TH system $(A;\lbrace E_i\rbrace_{i=0}^d;A^*;\lbrace E^*_i\rbrace_{i=0}^d)$ over $\K$. We now describe the third set. 
Let $(X, \lbrace \theta_i \rbrace_{i=0}^d, \lbrace \theta^*_i \rbrace_{i=0}^d)$ denote a normalized west-south Vandermonde system in $\Mdf$. One checks that $(X, \lbrace \alpha \theta_i + \beta \rbrace_{i=0}^d, \lbrace \alpha^* \theta^*_i + \beta^* \rbrace_{i=0}^d)$ is a normalized west-south Vandermonde system in $\Mdf$. These two systems are said to be affine related. 
We now describe the fourth set. By a normalized west-south Vandermonde matrix in $\Mdf$ we mean the matrix $X$ induced by a normalized west-south Vandermonde system $(X, \lbrace \theta_i \rbrace_{i=0}^d, \lbrace \theta^*_i \rbrace_{i=0}^d)$ in $\Mdf$. We now describe the fifth set. 
Let $(\lbrace \theta_i \rbrace_{i=0}^d, \lbrace \theta^*_i \rbrace_{i=0}^d, \lbrace \phi_i \rbrace_{i=1}^d)$ denote a parameter array over $\K$. Observe that $(\lbrace \alpha \theta_i + \beta \rbrace_{i=0}^d, \lbrace \alpha^* \theta^*_i + \beta^* \rbrace_{i=0}^d, \lbrace \alpha \alpha^* \phi_i \rbrace_{i=1}^d)$ is a parameter array over $\K$. These two parameter arrays are said to be affine related. 
This affine relation is an equivalence relation; the equivalence classes are called reduced parameter arrays. We have now described the five sets. We omit the description of the bijections between these sets as they are not hard to guess. 

This paper is organized as follows. In Sections 2, 3 we review some basic concepts regarding TH pairs and TH systems. In Section 4 we summarize the classification of TH systems given in \cite{thinhess}. In Section 5 we discuss affine transformations of a TH system. In Sections 6, 7 we discuss how a given TH system yields three more TH systems called the relatives. In Sections 8, 9 we discuss some scalars that are helpful in describing a given TH system. In Section 10 we use these scalars to describe the relatives of a given TH system. 
In Sections 11, 12 we discuss the transition matrix $\mathcal P$ and a related matrix $P$. In Section 13 we define the notion of a Vandermonde system. In Sections 14--16 we discuss the connection between Vandermonde systems, graded sequences of polynomials, and Hessenberg matrices. In Section 17 we discuss the double Vandermonde structure of the transition matrices $\mathcal P$ and $P$.
Sections 18, 19  contain the main results of the paper.

\section{Thin Hessenberg systems}
\label{sec:thsystem}
In our study of a TH pair, it is often helpful to consider a closely related object called a TH system. Before defining this notion, we make some definitions and observations. 
For the rest of the paper, fix an integer $d \geq 0$. Let $\Mdf$ denote the $\K$-algebra consisting of the $(d+1)\times(d+1)$ matrices that have all entries in $\K$. We index the rows and columns by $0,1, \ldots, d$. Let $\K^{d+1}$ denote the $\K$-vector space consisting of the $(d+1)\times1$ matrices that have all entries in $\K$. We index the columns by $0,1,\ldots, d$. Observe that $\Mdf$ acts on $\K^{d+1}$ by left multiplication.
For the rest of the paper, fix a vector space $V$ over $\fld$ with dimension $d+1$. 
Let ${\rm End}(V)$ denote the $\K$-algebra consisting of the linear transformations from $V$ to $V$. 
Suppose that $\lbrace v_i \rbrace_{i=0}^d$ is a basis for $V$. For $X \in \Mdf$ and $Y \in {\rm End}(V)$, we say that $X$ \emph{represents} $Y$ \emph{with respect to} $\lbrace v_i \rbrace_{i=0}^d$ whenever $Y v_j = \sum_{i=0}^d X_{ij} v_i$ for $0 \leq j \leq d$. 
For $A \in {\rm End}(V)$ and $W \subseteq V$,  we call $W$ an {\it eigenspace} of
$A$ whenever $W\not=0$ and there exists $\theta \in \K$ such that $W=\lbrace v \in V \;\vert \;Av = \theta v\rbrace$.
In this case $\theta$ is called the {\it eigenvalue} of $A$ corresponding to $W$. We say that $A$ is {\it diagonalizable} whenever $V$ is spanned by the eigenspaces of $A$.  We say that $A$ is {\it multiplicity-free} whenever $A$ is diagonalizable and each eigenspace of $A$ has dimension one. 

\begin{lemma}
\label{lem:multfree}
{\rm \cite[Lemma 2.1]{thinhess}}
Let $A, A^*$ denote a TH pair on $V$. Then each of $A, A^*$ is multiplicity-free. 
\end{lemma}

\noindent We recall a few more concepts from linear algebra. Let $A$ denote a multiplicity-free element of ${\rm End}(V)$. Let $\{V_i\}_{i=0}^d$ denote an ordering of the eigenspaces of $A$ and let $\{\theta_i\}_{i=0}^d$ denote the corresponding ordering of the eigenvalues of $A$. For $0 \leq i \leq d$, define $E_i \in {\rm End}(V)$ such that $(E_i-I)V_i=0$ and $E_iV_j=0$ for $j \neq i$ $(0 \leq j \leq d)$. Here $I$ denotes the identity of ${\rm End}(V)$. We call $E_i$ the {\it primitive idempotent} of $A$ corresponding to $V_i$ (or $\theta_i$).
Observe that
(i)  $I = \sum_{i=0}^d E_i$;
(ii) $E_iE_j=\delta_{i,j}E_i$ $(0 \leq i,j \leq d)$;
(iii) $V_i=E_iV$ $(0 \leq i \leq d)$;
(iv) $A=\sum_{i=0}^d \theta_i E_i$.
Moreover
\begin{equation}        
\label{eq:defEi}
E_i=\prod_{\stackrel{0 \leq j \leq d}{j \neq i}}
          \frac{A-\theta_jI}{\theta_i-\theta_j} \qquad \qquad \qquad (0 \leq i \leq d).
\end{equation}
Note that each of $\{A^i\}_{i=0}^d$, $\{E_i\}_{i=0}^d$ is a basis for the $\K$-subalgebra of $\mbox{\rm End}(V)$ generated by $A$. Moreover $\prod_{i=0}^d(A-\theta_iI)=0$.

\medskip
\noindent We now define a TH system. 

\begin{definition} 
\label{def:HS}
\rm
By a {\it thin Hessenberg system} (or {\it TH system}) on $V$ we mean a sequence
\[
\Phi=(A;\lbrace E_i\rbrace_{i=0}^d;A^*;\lbrace E^*_i\rbrace_{i=0}^d)
\]
which satisfies (i)--(v) below. 
\begin{enumerate}
\item[\rm (i)] 
Each of $A,A^*$ is a multiplicity-free element of ${\rm End}(V)$.
\item[\rm (ii)]
$\{E_i\}_{i=0}^d$ is an ordering of the primitive idempotents of $A$.
\item[\rm (iii)]
$\{E^*_i\}_{i=0}^d$ is an ordering of the primitive idempotents of $A^*$.
\item[\rm (iv)]
${\displaystyle{
E_iA^*E_j = \cases{0, &if $\; i-j > 1$\cr
\not=0, &if $\; i-j  = 1$\cr}
\qquad \qquad 
(0 \leq i,j\leq d)}}$.
\item[\rm (v)]
${\displaystyle{
 E^*_iAE^*_j = \cases{0, &if $\; i-j > 1$\cr
\not=0, &if $\; i-j  = 1$\cr}
\qquad \qquad 
(0 \leq i,j\leq d).}}$
\end{enumerate}
We refer to $d$ as the {\it diameter} of $\Phi$. We call $V$ the {\it underlying} vector space and say that $\Phi$ is {\it over} $\K$.  
\end{definition}

\noindent We comment on how TH pairs and TH systems are related.
Let $(A;\lbrace E_i\rbrace_{i=0}^d;A^*;\lbrace E^*_i\rbrace_{i=0}^d)$ denote a TH system on $V$. 
For $0 \leq i \leq d$, let $v_i$ (resp. $v^*_i$) denote a nonzero vector in $E_iV$ (resp. $E^*_iV$). Then the sequence $\lbrace v_i \rbrace_{i=0}^d$ (resp. $\lbrace v^*_i \rbrace_{i=0}^d$) is a basis for $V$ which satisfies Definition \ref{def:thinhess}(ii) (resp. Definition \ref{def:thinhess}(i)). Therefore the pair $A, A^*$ is a TH pair on $V$. Conversely, let $A, A^*$ denote a TH pair on $V$. Then each of $A, A^*$ is multiplicity-free by Lemma \ref{lem:multfree}. Let $\lbrace v_i \rbrace_{i=0}^d$ (resp. $\lbrace v^*_i \rbrace_{i=0}^d$) denote a basis for $V$ which satisfies Definition \ref{def:thinhess}(ii) (resp. Definition \ref{def:thinhess}(i)). For $0 \leq i \leq d$, the vector $v_i$ (resp. $v^*_i$) is an eigenvector for $A$ (resp. $A^*$); let $E_i$ (resp. $E^*_i$) denote the corresponding primitive idempotent. Then $(A;\lbrace E_i\rbrace_{i=0}^d;A^*;\lbrace E^*_i\rbrace_{i=0}^d)$ is a TH system on $V$. 

\begin{definition}
\label{def:thpairassthsytem}
\rm
Let $\Phi = (A;\lbrace E_i\rbrace_{i=0}^d;A^*;\lbrace E^*_i\rbrace_{i=0}^d)$ denote a TH system on $V$.
Observe that $A, A^*$ is a TH pair on $V$. We say that this pair is {\it associated} with $\Phi$. 
\end{definition}

\begin{remark}
\rm
With reference to Definition \ref{def:thpairassthsytem}, conceivably a given TH pair is associated with many TH systems. 
\end{remark}

\noindent We now recall several definitions and results on TH systems. 

\begin{definition}
\label{def:evseq}
\rm
Let $\Phi=(A;\lbrace E_i\rbrace_{i=0}^d;A^*;\lbrace E^*_i\rbrace_{i=0}^d)$ denote a TH system on $V$.
For $0 \leq i \leq d$, let $\theta_i $ (resp. $\theta^*_i$) denote the eigenvalue of $A$ (resp. $A^*$) corresponding to $E_i$ (resp. $E^*_i$). We refer to $\lbrace \theta_i \rbrace_{i=0}^d$ as the {\it eigenvalue sequence} of $\Phi$. We refer to  $\lbrace \theta^*_i \rbrace_{i=0}^d$ as the {\it dual eigenvalue sequence} of $\Phi$. We observe that $\lbrace \theta_i \rbrace_{i=0}^d$ are mutually distinct and contained in $\K$. Similarly $\lbrace \theta^*_i \rbrace_{i=0}^d$ are mutually distinct and contained in $\K$. 
\end{definition}

\begin{definition}
\label{def:evseqhp}
\rm 
Let $A, A^*$ denote a TH pair. By an {\it eigenvalue sequence} (resp. {\it dual eigenvalue sequence}) of $A, A^*$, we mean the eigenvalue sequence (resp. dual eigenvalue sequence) of an associated TH system. We emphasize that a given TH pair could have many eigenvalue and dual eigenvalue sequences.
\end{definition}

\smallskip

\noindent Let $\K[\lambda]$ denote the $\K$-algebra consisting of the polynomials in $\lambda$ that have all coefficients in $\K$.

\begin{notation}
\label{def:tau}
\rm
Let $\lbrace \theta_i \rbrace_{i=0}^d, \lbrace \theta^*_i \rbrace_{i=0}^d$ denote two sequences of scalars taken from $\K$. For $0 \leq i \leq d+1$, let $\tau_i$, $\tau^*_i$, $\eta_i$, $\eta^*_i$
denote the following polynomials in $\K[\lambda]$.
\begin{eqnarray*}
&&\tau_i = \prod_{h=0}^{i-1} (\lambda - \theta_h),
\qquad \qquad \;\;
\eta_i = \prod_{h=0}^{i-1}(\lambda - \theta_{d-h}),
\\
&&\tau^*_i = \prod_{h=0}^{i-1} (\lambda - \theta^*_h),
\qquad \qquad
\eta^*_i = \prod_{h=0}^{i-1}(\lambda - \theta^*_{d-h}).
\end{eqnarray*}
\noindent We observe that each of $\tau_i$, $\eta_i$, $\tau^*_i$, $\eta^*_i$ is monic with degree $i$. 
\end{notation} 

\noindent By (\ref{eq:defEi}), for $0 \leq i \leq d$ 
\begin{eqnarray}
\label{eq:primidtau}
E_i = \frac{\tau_i(A)\eta_{d-i}(A)}{\tau_i(\theta_i)\eta_{d-i}(\theta_i)}, \qquad \qquad  \ E^*_i =\frac{\tau^*_i(A^*)\eta^*_{d-i}(A^*)}{\tau^*_i(\theta^*_i)\eta^*_{d-i}(\theta^*_i)}.
\end{eqnarray}

By a {\it decomposition} of $V$ we mean a sequence $\lbrace U_i \rbrace_{i=0}^d$ consisting of one-dimensional subspaces of $V$ such that
\begin{eqnarray*}
\qquad \qquad V=U_0+U_1+\cdots + U_d \qquad
\qquad (\hbox{direct sum}).
\end{eqnarray*}
For notational convenience, set $U_{-1}=0$ and $U_{d+1}=0$. 

Let $\Phi=(A;\lbrace E_i\rbrace_{i=0}^d;A^*;\lbrace E^*_i\rbrace_{i=0}^d)$ denote a TH system on $V$. Then $\lbrace E^*_iV \rbrace_{i=0}^d$ is a decomposition of $V$, said to be {\it $\Phi$-standard}. 
Let $0 \neq \xi_0 \in E_0V$. The sequence $\lbrace E^*_i \xi_0 \rbrace_{i=0}^d$ is a basis for $V$ \cite[Lemma 8.1]{thinhess}, said to be {\it $\Phi$-standard}.  We recall another decomposition of $V$ associated with $\Phi$. For $0 \leq i \leq d$, let

\begin{equation}
\label{eq:defui}
U_i = (E^*_0V + E^*_1V + \cdots + E^*_iV) \cap (E_0V + E_{1}V + \cdots + E_{d-i}V).
\end{equation}

\noindent The sequence $\lbrace U_i \rbrace_{i=0}^d$ is a decomposition of $V$ \cite[Section 4]{thinhess}, said to be {\it $\Phi$-split}. Moreover for $0 \leq i \leq d$, both
\begin{eqnarray}
(A-\theta_{d-i} I)U_i &=& U_{i+1},
\label{eq:raise}
\\
(A^*-\theta^*_i I)U_i &=& U_{i-1}.
\label{eq:lower}
\end{eqnarray}

\noindent Setting $i=d$ in (\ref{eq:defui}) we find $U_d=E_0V$. Combining this with (\ref{eq:lower}) we find  
\begin{equation}
U_i = \eta^*_{d-i}(A^*) E_0V
\qquad \qquad (0 \leq i\leq d).
\label{eq:uialt}
\end{equation}
Recall $0 \neq \xi_0 \in E_0V$. From (\ref{eq:uialt}) we find that for $0 \leq i \leq d$, the vector $\eta^*_{d-i}(A^*) \xi_0$ is a basis for $U_i$. By this and since $\lbrace U_i \rbrace_{i=0}^d$ is a decomposition of $V$, the sequence
\begin{eqnarray}
\label{eq:splitbasis}
\eta^*_{d-i}(A^*) \xi_0  \qquad \qquad (0 \leq i\leq d)
\end{eqnarray}
is a basis for $V$, said to be {\it $\Phi$-split}. Let $1 \leq i \leq d$. By (\ref{eq:lower}) we have $(A^*-\theta^*_i I)U_i = U_{i-1}$, and by (\ref{eq:raise}) we have $(A-\theta_{d-i+1} I)U_{i-1} = U_i$. Therefore $U_i$ is invariant under $(A-\theta_{d-i+1}I)(A^*-\theta^*_i I)$ and the corresponding eigenvalue is a nonzero element of $\K$. We denote this eigenvalue by $\phi_i$. We call the sequence $\lbrace \phi_i \rbrace_{i=1}^d$ the {\it split sequence} of $\Phi$. For notational convenience, set $\phi_0=0$ and $\phi_{d+1} = 0$.

\begin{proposition}
{\rm \cite[Proposition 4.4]{thinhess}}
\label{prop:aastarsplitrep}
Let $\Phi=(A;\lbrace E_i\rbrace_{i=0}^d;A^*;\lbrace E^*_i\rbrace_{i=0}^d)$ denote a TH system on $V$ with eigenvalue sequence 
$\lbrace \theta_i \rbrace_{i=0}^d$, dual eigenvalue sequence $\lbrace \theta^*_i \rbrace_{i=0}^d$, and split sequence 
$\lbrace \phi_i \rbrace_{i=1}^d$. Then the matrices representing $A$ and $A^*$ with respect to a $\Phi$-split basis for $V$ are  
\begin{equation}
\label{eq:matrepaastar}
\left(
\begin{array}{c c c c c c}
\theta_d & & & & & {\bf 0} \\
\phi_1 & \theta_{d-1} &  & & & \\
& \phi_2 & \theta_{d-2} &  & & \\
& & \cdot & \cdot &  &  \\
& & & \cdot & \cdot &  \\
{\bf 0}& & & & \phi_d & \theta_0
\end{array}
\right),
\qquad \quad 
\left(
\begin{array}{c c c c c c}
\theta^*_0 & 1 & & & & {\bf 0} \\
& \theta^*_1 & 1 & & & \\
& & \theta^*_2 & \cdot & & \\
& &  & \cdot & \cdot &  \\
& & &  & \cdot & 1 \\
{\bf 0}& & & &  & \theta^*_d
\end{array}
\right)
\end{equation}
\noindent respectively. 
\end{proposition}

\noindent Next we describe the matrices representing the primitive idempotents of $A$, $A^*$ with respect to a $\Phi$-split basis for $V$. 

\begin{proposition}
\label{prop:formofprimidS99}
Let $\Phi=(A;\lbrace E_i\rbrace_{i=0}^d;A^*;\lbrace E^*_i\rbrace_{i=0}^d)$ denote a TH system on $V$ with eigenvalue sequence $\lbrace \theta_i \rbrace_{i=0}^d$, dual eigenvalue sequence $\lbrace \theta^*_i \rbrace_{i=0}^d$, and split sequence $\lbrace \phi_i \rbrace_{i=1}^d$. For $0 \leq r \leq d$, consider the matrices in $\Mdf$ that represent $E_r$ and $E^*_r$ with respect to a $\Phi$-split basis. For $0 \leq i,j\leq d$, their $(i,j)$-entry is described as follows.
For $E_r$ this entry is 
\begin{eqnarray}
\label{eq:entriesofefreefinS99}
{{\phi_1 \phi_2 \cdots \phi_i}
\over{\phi_1 \phi_2 \cdots \phi_j}} \,
{{\tau_{d-i}(\theta_r) \eta_j(\theta_r)}\over
{\tau_r(\theta_r)\eta_{d-r}(\theta_r)}},
\end{eqnarray}
and for $E^*_r$ this entry is
\begin{eqnarray}
{{\tau^*_i(\theta^*_r)\eta^*_{d-j}(\theta^*_r)}\over{
\tau^*_r(\theta^*_r)\eta^*_{d-r}(\theta^*_r)}}.
\label{eq:theestarentriesfinS99}
\end{eqnarray}
\end{proposition}

\noindent {\it Proof:} Fix a $\Phi$-split basis for $V$. For notational convenience, identify each element of ${\rm End}(V)$ with the matrix in $\Mdf$ that represents it with respect to this basis. We first show that the $(i,j)$-entry of $E_r^*$ is given by (\ref{eq:theestarentriesfinS99}). 
Computing the $(i,j)$-entry of $A^*E^*_r = \theta^*_r E^*_r$ using matrix multiplication, and taking into account the form of $A^*$ in (\ref{eq:matrepaastar}), we find
\begin{eqnarray*}
(E^*_r)_{i+1,j} = (\theta^*_r-\theta^*_i)(E^*_r)_{ij} 
\end{eqnarray*}
if $i \leq d-1$.
Replacing $i$ by $i-1$ in the above line, we find
\begin{eqnarray}
(E^*_r)_{ij} = (\theta^*_r-\theta^*_{i-1})(E^*_r)_{i-1,j} 
\label{eq:firstrecS99}
\end{eqnarray}
if $i \geq 1$.
Using the recursion (\ref{eq:firstrecS99}), we routinely find
\begin{eqnarray}
(E^*_r)_{ij} &=& (\theta^*_r-\theta^*_{i-1})(\theta^*_r-\theta^*_{i-2}) \cdots (\theta^*_r- \theta^*_0)(E^*_r)_{0j}
\nonumber
\\
&=& \tau^*_{i}(\theta^*_r)(E^*_r)_{0j}.
\label{eq:firstpartrecS99}
\end{eqnarray}
Computing the $(0,j)$-entry of $E^*_r A^* = \theta^*_r E^*_r$ using matrix multiplication, and taking into account the form of $A^*$, we find
\begin{eqnarray*}
(E^*_r)_{0,j-1} = (\theta^*_r - \theta^*_j)(E^*_r)_{0j} 
\end{eqnarray*}
if $j \geq 1$.  Replacing $j$ by $j+1$ in the above line we find 
\begin{eqnarray}
(E^*_r)_{0j} = (\theta^*_r - \theta^*_{j+1})(E^*_r)_{0,j+1} 
\label{eq:secondrecS99}
\end{eqnarray}
if $j \leq d-1$.
Using the recursion (\ref{eq:secondrecS99}), we routinely find
\begin{eqnarray}
(E^*_r)_{0j} &=& (\theta^*_r - \theta^*_{j+1})(\theta^*_r - \theta^*_{j+2}) \cdots (\theta^*_r - \theta^*_d)(E^*_r)_{0d}
\nonumber
\\
&=& \eta^*_{d-j}(\theta^*_r)(E^*_r)_{0d}.
\label{eq:secondpartrecS99}
\end{eqnarray}
Combining (\ref{eq:firstpartrecS99}), (\ref{eq:secondpartrecS99}), we find
\begin{eqnarray}
(E^*_r)_{ij} = \tau^*_{i}(\theta^*_r) \eta^*_{d-j}(\theta^*_r)c,
\label{eq:togetherS99}
\end{eqnarray}
where we abbreviate $c = (E^*_r)_{0d}$. We now find $c$. Since $A^*$ is upper triangular, and since $E^*_r$ is a polynomial in $A^*$, we see 
$E^*_r$ is upper triangular. Recall $E_r^{*2}=E^*_r$, so the diagonal entry of $(E^*_r)_{rr}$ equals 0 or 1. We show $(E^*_r)_{rr} = 1$.
Setting $i=r$, $j=r$ in (\ref{eq:togetherS99}),
\begin{eqnarray}
(E^*_r)_{rr} = \tau^*_{r}(\theta^*_r) \eta^*_{d-r}(\theta^*_r)c.
\label{eq:rrentryoferS99}
\end{eqnarray}
Observe $\tau^*_{r}(\theta^*_r) \neq 0$ and $\eta^*_{d-r}(\theta^*_r) \neq 0$ by Notation \ref{def:tau}, and since $\lbrace \theta^*_i \rbrace_{i=0}^d$ are distinct. Observe $c \neq 0$; otherwise $E^*_r = 0$ in view of (\ref{eq:togetherS99}). Apparently the right side of
(\ref{eq:rrentryoferS99}) is not 0, so $(E^*_r)_{rr} \neq 0$, and we conclude $(E^*_r)_{rr} = 1$. Setting $(E^*_r)_{rr} = 1$ in 
(\ref{eq:rrentryoferS99}), solving for $c$, and evaluating  (\ref{eq:togetherS99}) using the result, we find the $(i,j)$-entry 
of $E^*_r$ is given by (\ref{eq:theestarentriesfinS99}). \\
We now show that the $(i,j)$-entry of $E_r$ is given by (\ref{eq:entriesofefreefinS99}). Let $G \in \Mdf$ denote the diagonal matrix with $(i,i)$-entry $\phi_1 \phi_2 \cdots \phi_i$ for $0 \leq i \leq d$ and set $A':=GA^{t}G^{-1}$, where $A$ is the matrix on the left of (\ref{eq:matrepaastar}). The matrix $A'$ is equal to
\begin{eqnarray}
\label{eq:mataprime}
\left(
\begin{array}{c c c c c c}
\theta_d & 1 & & & & {\bf 0} \\
& \theta_{d-1} & 1 & & & \\
& & \theta_{d-2} & \cdot & & \\
& &  & \cdot & \cdot &  \\
& & &  & \cdot & 1 \\
{\bf 0}& & & &  & \theta_0
\end{array}
\right)
\end{eqnarray} 
Let $E'_r$ denote the primitive idempotent of $A'$ associated with the eigenvalue $\theta_r$. We find $E'_r$ in two ways. On one hand, applying 
(\ref{eq:theestarentriesfinS99}) to $A'$, we find $E'_r$ has $(i,j)$-entry 
\begin{eqnarray}
{{\tau_{d-j}(\theta_r)\eta_{i}(\theta_r)}\over
{\tau_r(\theta_r)\eta_{d-r}(\theta_r)}}
\label{eq:entriesofestarAfreeS99}
\end{eqnarray}
for $0 \leq i,j\leq d$. On the other hand, by elementary linear algebra
\begin{eqnarray*}
E'_r = GE^{t}_rG^{-1},
\end{eqnarray*}
so $E'_r$ has $(i,j)$-entry 
\begin{eqnarray}
G_{ii}(E_r)_{ji}G^{-1}_{jj} = {{\phi_1 \phi_2 \cdots \phi_i} \over {\phi_1 \phi_2 \cdots \phi_j}} \, (E_r)_{ji}
\label{eq:theestarentriesyzS99}
\end{eqnarray}
for $0 \leq i,j \leq d$. Equating (\ref{eq:entriesofestarAfreeS99}) and the right side of (\ref{eq:theestarentriesyzS99}),
and solving for $(E_r)_{ji}$, we routinely obtain the result. \hfill $\Box$

\begin{example}
\label{ex:primidsdeq2freeS99}
\rm
Referring to Proposition \ref{prop:formofprimidS99}, assume $d = 2$. With respect to a $\Phi$-split basis, the matrices representing $E_0, E_1, E_2$ are

\bigskip

$
\left(
\begin{array}{c c c }
0 & 0 & 0 \\
0 & 0 & 0 \\
{{\phi_1 \phi_2}\over{(\theta_0 - \theta_2)(\theta_0 - \theta_1)}} & {{\phi_2}\over {\theta_0-\theta_1}} & 1
\end{array}
\right), \ \ \left(
\begin{array}{c c c }
0 & 0  &  0
\\
{{\phi_1}\over {\theta_1-\theta_2}} & 1 & 0 
\\
{{\phi_1\phi_2}\over {(\theta_1 - \theta_0)(\theta_1 - \theta_2)}} & {{\phi_2}\over {\theta_1-\theta_0}} & 0
\end{array}
\right), \ \ \left(
\begin{array}{c c c }
1 & 0  &  0
\\
{{\phi_1}\over {\theta_2-\theta_1}} & 0 & 0
\\
{{\phi_1\phi_2}\over {(\theta_2 - \theta_0)(\theta_2 - \theta_1)}} & 0 & 0
\end{array}
\right),$ 

\bigskip

\noindent respectively. Moreover the matrices representing $E_0^*, E_1^*, E_2^*$ are

\bigskip

$
\left(
\begin{array}{c c c }
1 & {{1}\over{\theta^*_0-\theta^*_1}} & {{1}\over{(\theta^*_0-\theta^*_1)(\theta^*_0-\theta^*_2)}}
\\
0 & 0 & 0
\\
0 & 0 & 0
\end{array}
\right), \ \ 
\left(
\begin{array}{c c c }
0 & {{1}\over{\theta^*_1-\theta^*_0}} & {{1}\over{(\theta^*_1-\theta^*_0)(\theta^*_1-\theta^*_2)}}
\\
0 & 1 & {{1}\over{\theta^*_1-\theta^*_2}}
\\
0 & 0 & 0
\end{array}
\right), \ \ 
\left(
\begin{array}{c c c }
0 & 0 & {{1}\over{(\theta^*_2-\theta^*_1)(\theta^*_2-\theta^*_0)}}
\\
0 & 0 & {{1}\over{\theta^*_2-\theta^*_1}}
\\
0 & 0 & 1
\end{array}
\right),$ 
\smallskip
respectively. 
\end{example}

\noindent We now give some characterizations of the split sequence. 

\begin{lemma}
\label{lem:splitchar}
Let $(A;\lbrace E_i\rbrace_{i=0}^d;A^*;\lbrace E^*_i\rbrace_{i=0}^d)$ denote a TH system with eigenvalue sequence $\lbrace \theta_i \rbrace_{i=0}^d$, dual eigenvalue sequence $\lbrace \theta^*_i \rbrace_{i=0}^d$, and split sequence $\lbrace \phi_i \rbrace_{i=1}^d$. Then  
\begin{eqnarray}
\label{eq:label100}
E_0^* \eta_i(A) E_0^* = \frac{\phi_1 \phi_2 \cdots \phi_i}{(\theta^*_0 - \theta^*_1)(\theta^*_0 - \theta^*_2) \cdots (\theta^*_0 - \theta^*_i)} E^*_0 \qquad \ \ (0 \leq i \leq d).
\end{eqnarray}
\end{lemma}

\noindent {\it Proof:} Let $\Phi$ denote the TH system in question and assume $V$ is the underlying vector space. Let $\{U_i\}_{i=0}^d$ denote the $\Phi$-split decomposition of $V$. Setting $i=0$ in (\ref{eq:defui}) we find $U_0=E^*_0V$. By this and (\ref{eq:raise}), (\ref{eq:lower}) we obtain
\begin{equation}    
\label{eq:zetaiaux}
(A^*-\theta^*_1I)(A^*-\theta^*_2I)\cdots(A^*-\theta^*_iI) \eta_i(A) = \phi_1 \phi_2 \cdots \phi_i I
\end{equation}
on $E^*_0V$.
To obtain (\ref{eq:label100}), multiply both sides of (\ref{eq:zetaiaux}) on the left by $E^*_0$ and use $E^*_0A^*= \theta^*_0E^*_0$. \hfill $\Box$

\begin{corollary}
\label{cor:splitchar}
Let $(A;\lbrace E_i\rbrace_{i=0}^d;A^*;\lbrace E^*_i\rbrace_{i=0}^d)$ denote a TH system with eigenvalue sequence $\lbrace \theta_i \rbrace_{i=0}^d$, dual eigenvalue sequence $\lbrace \theta^*_i \rbrace_{i=0}^d$, and split sequence $\lbrace \phi_i \rbrace_{i=1}^d$. Then for $0 \leq i \leq d$, 
\begin{eqnarray}
\label{eq:lab1}
\phi_1 \phi_2 \cdots \phi_i = (\theta^*_0 - \theta^*_1)(\theta^*_0 - \theta^*_2) \cdots (\theta^*_0 - \theta^*_i) {\rm trace}(\eta_i(A) E_0^*).
\end{eqnarray}
Moreover $\eta_i(A) E_0^*$ has nonzero trace. 
\end{corollary}

\noindent {\it Proof:} To obtain (\ref{eq:lab1}), in (\ref{eq:label100}) take the trace of each side and simplify the result using the fact that ${\rm trace}(E_0^*) = 1$ and ${\rm trace}(E_0^* \eta_i(A) E_0^*) = {\rm trace}(\eta_i(A)E_0^*E_0^*) = {\rm trace}( \eta_i(A)E_0^*)$. This gives (\ref{eq:lab1}). The last assertion follows since $\phi_i \neq 0$ for $1 \leq i \leq d$. \hfill $\Box$

\begin{corollary}
\label{cor:splitcharref} 
Let $(A;\lbrace E_i\rbrace_{i=0}^d;A^*;\lbrace E^*_i\rbrace_{i=0}^d)$ denote a TH system with eigenvalue sequence $\lbrace \theta_i \rbrace_{i=0}^d$, dual eigenvalue sequence $\lbrace \theta^*_i \rbrace_{i=0}^d$, and split sequence $\lbrace \phi_i \rbrace_{i=1}^d$. Then 
\begin{eqnarray}
\label{eq:lab3}
\phi_i = (\theta^*_0 - \theta^*_i) {\rm trace}(\eta_i(A) E_0^*)/{\rm trace}(\eta_{i-1}(A) E_0^*)  \qquad \ \ (1 \leq i \leq d).
\end{eqnarray}
\end{corollary}

\noindent {\it Proof:} Routine by Corollary \ref{cor:splitchar}.  \hfill $\Box$

\medskip

\noindent In Section \ref{sec:z2action} we give some more characterizations of the split sequence.

\section{Isomorphisms for TH pairs and TH systems}
In this section we discuss the notion of isomorphism for TH pairs and TH systems. 

\begin{lemma}
\label{lem:multfreecomm}
For $X \in \Mdf$ the following {\rm (i)--(iii)} are equivalent. 
\begin{enumerate}
\item[\rm (i)]
$X$ is diagonal.  
\item[\rm (ii)]
$D X = X D$ for all diagonal $D \in \Mdf$.
\item[\rm (iii)]
There exists a diagonal $D \in \Mdf$ that has mutually distinct diagonal entries and $D X = X D$.
\end{enumerate}
\end{lemma}

\noindent {\it Proof:} (i) $\Rightarrow$ (ii) Clear. \\
(ii) $\Rightarrow$ (iii) Clear. \\
(iii) $\Rightarrow$ (i) For $0 \leq i, j \leq d$ with $i \neq j$, we show $X_{ij} = 0$. Comparing the $(i,j)$-entry of $DX$ and $XD$, we find $D_{ii}X_{ij} = X_{ij} D_{jj}$. By assumption $D_{ii} \neq D_{jj}$, so $X_{ij} = 0$. \hfill $\Box$

\medskip

Let $A, A^*$ denote a TH pair on $V$. In general, ${\rm End}(V)$ may not be generated by $A, A^*$. Moreover there may exist a subspace $W$ of $V$ such that $AW \subseteq W, A^*W \subseteq W, W \neq 0, W \neq V$. However we do have the following result.

\begin{lemma}
\label{lem:scalar}
Let $A, A^*$ denote a TH pair on $V$. Let $\Delta$ denote an element of ${\rm End}(V)$ such that $\Delta A = A \Delta$ and  $\Delta A^* = A^* \Delta$. Then $\Delta \in \K I$.   
\end{lemma}

\noindent {\it Proof:} Pick a basis for $V$ from Definition \ref{def:thinhess}(i). For notational convenience, identify each element of ${\rm End}(V)$ with the matrix that represents it with respect to this basis. Thus the matrix $A$ is Hessenberg and the matrix $A^*$ is diagonal. Moreover the diagonal entries of $A^*$ are mutually distinct by Lemma \ref{lem:multfree}. Applying Lemma \ref{lem:multfreecomm} with $D = A^*$ and $X = \Delta$, we find $\Delta$ is diagonal. For $1 \leq i \leq d$, comparing the $(i,i-1)$-entry of $\Delta A$ and $A \Delta$, we find $\Delta_{ii} A_{i,i-1} = A_{i,i-1} \Delta_{i-1,i-1}$. Observe that $A_{i,i-1} \neq 0$ since $A$ is Hessenberg, so $\Delta_{ii} = \Delta_{i-1,i-1}$. Therefore $\Delta_{ii}$ is independent of $i$ for $0 \leq i \leq d$. Consequently $\Delta \in \K I$. \hfill $\Box$

\medskip

\noindent For the rest of this section, let $W$ denote a vector space over $\fld$ with dimension $d+1$. Let $\Gamma:V \to W$ denote a $\mathbb{K}$-vector space isomorphism. Then there exists a unique $\fld$-algebra isomorphism $\gamma :{\rm End}(V) \to {\rm End}(W)$ such that $S^\gamma=\Gamma S \Gamma^{-1}$ for all $S \in {\rm End}(V)$. Conversely let $\gamma: {\rm End}(V) \to {\rm End}(W)$ denote a $\fld$-algebra isomorphism. By the Skolem-Noether theorem \cite[Corollary 9.122]{Rot} there exists a  $\fld$-vector space isomorphism $\Gamma: V \to W$ such that $S^\gamma = \Gamma S \Gamma^{-1}$ for all $S \in {\rm End}(V)$. Moreover $\Gamma$ is unique up to multiplication by a nonzero scalar in $\K$.

\begin{definition}
\label{def:isompair}
\rm
Let $A, A^*$ denote a TH pair on $V$ and let $B, B^*$ denote a TH pair on $W$. By an {\it isomorphism of TH pairs} from $A, A^*$ to $B, B^*$ we mean a $\fld$-algebra isomorphism $\gamma: {\rm End}(V) \rightarrow {\rm End}(W)$ such that $B = A^\gamma$ and $B^* = A^{*\gamma}$. We say that the TH pairs $A, A^*$ and $B, B^*$ are {\it isomorphic} whenever there exists an isomorphism of TH pairs from $A, A^*$ to $B, B^*$.
\end{definition}

\begin{lemma}
\label{lem:isouniquepair}
Let $A, A^*$ and $B, B^*$ denote isomorphic TH pairs over $\fld$. Then the isomorphism of TH pairs from $A, A^*$ to $B, B^*$ is unique. 
\end{lemma}

\noindent {\it Proof:} Let $\gamma$ and $\gamma'$ denote isomorphisms of TH pairs from $A, A^*$ to $B, B^*$. We show that $\gamma = \gamma'$. By the comments above Definition \ref{def:isompair}, there exists a  $\fld$-vector space isomorphism $\Gamma: V \to W$ (resp. $\Gamma': V \to W$) such that $S^\gamma = \Gamma S \Gamma^{-1}$ (resp. $S^{\gamma'} = \Gamma' S \Gamma^{\prime -1}$) for all $S \in {\rm End}(V)$.
Consider the composition $\Delta = \Gamma^{-1} \Gamma'$.
Observe that $\Delta$ is an invertible element of ${\rm End}(V)$. By construction, $\Delta A = A \Delta$ and $\Delta A^* = A^* \Delta$. Therefore $\Delta \in \K I$ by Lemma \ref{lem:scalar}. By these comments, there exists $0 \neq \alpha \in \K$ such that 
$\Delta = \alpha I$. Hence $\Gamma' = \alpha \Gamma$, so $\gamma = \gamma'$. \hfill $\Box$

\begin{definition}
\label{def:isomsystem}
\rm
Let $\Phi=(A;\lbrace E_i\rbrace_{i=0}^d;A^*;\lbrace E^*_i\rbrace_{i=0}^d)$ denote a TH system on $V$ and let $\Psi=(B;\lbrace F_i \rbrace_{i=0}^d;B^*;\lbrace F^*_i\rbrace_{i=0}^d)$ denote a TH system on $W$. By an {\it isomorphism of TH systems} from $\Phi$ to $\Psi$ we mean 
a $\fld$-algebra isomorphism $\gamma: {\rm End}(V) \rightarrow {\rm End}(W)$ such that 
\begin{eqnarray*}
B = A^\gamma, \qquad B^* = A^{*\gamma}, \qquad F_i = E_i^\gamma, \qquad F^*_i = E^{*\gamma}_i   \ \ \ (0 \leq i \leq d).
\end{eqnarray*}
We say that the TH systems $\Phi$ and $\Psi$ are {\it isomorphic} whenever there exists an isomorphism of TH systems from $\Phi$ to $\Psi$. 
\end{definition}

\begin{lemma}
\label{lem:isouniquesystem}
Let $\Phi$ and $\Psi$ denote isomorphic TH systems over $\fld$. Then the isomorphism of TH systems from $\Phi$ to $\Psi$ is unique. 
\end{lemma}
\noindent {\it Proof:} Similar to the proof of Lemma \ref{lem:isouniquepair}. \hfill $\Box$

\medskip

\noindent We give another interpretation of isomorphism for TH pairs and TH systems.  

\begin{lemma}          
\label{lem:isomaltpair}  
Let $A, A^*$ denote a TH pair on $V$ and let $B, B^*$ denote a TH pair on $W$. Then the following {\rm (i)}, {\rm (ii)} are equivalent.
\begin{itemize}
\item[\rm (i)]
The TH pairs $A, A^*$ and $B, B^*$ are isomorphic.
\item[\rm (ii)]
There exists a $\K$-vector space isomorphism $\Gamma: V \rightarrow W$ such that $B\Gamma = \Gamma A$ and $B^*\Gamma = \Gamma A^*$.
\end{itemize}
Moreover assume {\rm (i)}, {\rm (ii)} hold. Then $\Gamma$ is unique up to a multiplication by a nonzero scalar in $\K$. 
\end{lemma}

\begin{lemma}          
\label{lem:isomaltsystem}  
Let $\Phi=(A;\lbrace E_i\rbrace_{i=0}^d;A^*;\lbrace E^*_i\rbrace_{i=0}^d)$ denote a TH system on $V$ and let $\Psi=(B;\lbrace F_i \rbrace_{i=0}^d;B^*;\lbrace F^*_i\rbrace_{i=0}^d)$ denote a TH system on $W$. Then the following {\rm (i)}, {\rm (ii)} are equivalent.
\begin{itemize}
\item[\rm (i)]
The TH systems $\Phi$ and $\Psi$ are isomorphic.
\item[\rm (ii)]
There exists a $\K$-vector space isomorphism $\Gamma: V \rightarrow W$ such that 
\begin{eqnarray*}
\label{eq:isosys}
B \Gamma = \Gamma A, \qquad B^*\Gamma = \Gamma A^*, \qquad F_i \Gamma = \Gamma E_i, \qquad  F^*_i \Gamma = \Gamma E^*_i  \ \ \ (0 \leq i \leq d).
\end{eqnarray*}
\end{itemize}
Moreover assume {\rm (i)}, {\rm (ii)} hold. Then $\Gamma$ is unique up to a multiplication by a nonzero scalar in $\K$. 
\end{lemma}

\section{The classification of TH systems}
In \cite{thinhess} we classified the TH systems up to isomorphism. We recall the result in this section.  

\begin{definition}
\label{def:paofphi}
\rm
Let $\Phi$ denote a TH system. By the {\it parameter array} of $\Phi$ we mean the sequence $(\lbrace \theta_i \rbrace_{i=0}^d, \lbrace \theta^*_i \rbrace_{i=0}^d, \lbrace \phi_i \rbrace_{i=1}^d)$, where 
$\lbrace \theta_i \rbrace_{i=0}^d$ (resp. $\lbrace \theta^*_i \rbrace_{i=0}^d$) is the eigenvalue (resp. dual eigenvalue) sequence of $\Phi$ and $\lbrace \phi_i \rbrace^d_{i=1}$ is the split sequence of $\Phi$. 
\end{definition}

\begin{theorem}
{\rm \cite[Theorem 6.3]{thinhess}}
\label{thm:classificationpa}
Let  
\begin{eqnarray}
\label{eq:pa}
(\lbrace \theta_i \rbrace_{i=0}^d, \lbrace \theta^*_i \rbrace_{i=0}^d, \lbrace \phi_i \rbrace_{i=1}^d)
\end{eqnarray}
denote a sequence of scalars taken from $\K$. Then there exists a TH system $\Phi$ over $\K$ with parameter array (\ref{eq:pa}) if and only if {\rm (i)--(iii)} hold below. 
\begin{enumerate}
\item[\rm (i)] 
$ \theta_i\not=\theta_j \qquad $ \ \ \ if $\;\;i\not=j \qquad \qquad \qquad \ (0 \leq i,j\leq d)$.
\item[\rm (ii)]
$\theta^*_i\not=\theta^*_j \qquad $ \ \ if $\;\;i\not=j \qquad \qquad \qquad \ \ (0 \leq i,j\leq d)$.
\item[\rm (iii)] 
$\phi_i \not=0 \qquad \qquad \qquad \qquad \qquad \qquad \qquad (1 \leq i \leq d)$.
\end{enumerate}
\noindent Moreover assume {\rm (i)--(iii)} hold. Then $\Phi$ is unique up to isomorphism of TH systems.
\end{theorem}

\begin{definition}
\rm
By a {\it parameter array over $\fld$ of diameter $d$} we mean a sequence of scalars \\
$(\lbrace \theta_i \rbrace_{i=0}^d, \lbrace \theta^*_i \rbrace_{i=0}^d, \lbrace \phi_i \rbrace_{i=1}^d)$ taken from $\fld$ that satisfies conditions (i)--(iii) of Theorem \ref{thm:classificationpa}.  \end{definition}

\begin{corollary}
\label{cor:thsyspabij}
The map which sends a given TH system to its parameter array induces a bijection from the set of isomorphism classes of TH systems over $\K$ of diameter $d$, to the set of parameter arrays over $\fld$ of diameter $d$.  
\end{corollary}

\noindent {\it Proof:} Immediate from Theorem \ref{thm:classificationpa}. \hfill $\Box$

\medskip

To illuminate the bijection in Corollary \ref{cor:thsyspabij} we now describe its inverse in concrete terms. Let $\pi$ denote the bijection in Corollary \ref{cor:thsyspabij}.

\begin{proposition}
Let $(\lbrace \theta_i \rbrace_{i=0}^d, \lbrace \theta^*_i \rbrace_{i=0}^d, \lbrace \phi_i \rbrace_{i=1}^d)$ denote a parameter array over $\K$ of diameter $d$. Let $A$ (resp. $A^*$) denote the matrix on the left (resp. right) in (\ref{eq:matrepaastar}). Observe that $A$ (resp. $A^*$) is multiplicity-free with eigenvalues $\lbrace \theta_i \rbrace_{i=0}^d$ (resp. $\lbrace \theta_i^* \rbrace_{i=0}^d$). For $0 \leq i \leq d$ let $E_i$
(resp. $E_i^*$) denote the primitive idempotent of $A$ (resp. $A^*$) that corresponds to $\theta_i$ (resp. $\theta_i^*$). Then $\Phi=(A; \{E_i\}_{i=0}^d;A^*; \{E^*_i\}_{i=0}^d)$ is a TH system over $\K$. Moreover $\pi^{-1}$ sends $(\lbrace \theta_i \rbrace_{i=0}^d, \lbrace \theta^*_i \rbrace_{i=0}^d, \lbrace \phi_i \rbrace_{i=1}^d)$ to the isomorphism class of $\Phi$. 
\end{proposition}

\noindent {\it Proof:} This is proven as part of the proof of \cite[Theorem 6.3]{thinhess}. \hfill $\Box$

\section{The affine transformations of a TH system}
A given TH system can be modified in several ways to get a new TH system. In this section we describe one way. In the next section we describe another way.  

\begin{lemma}         
\label{lem:affine}  
Let $\Phi=(A; \{E_i\}_{i=0}^d;A^*; \{E^*_i\}_{i=0}^d)$ denote a TH system on $V$. Let $\alpha, \beta, \alpha^*, \beta^*$ denote scalars in $\K$ with $\alpha, \alpha^*$ nonzero. Then the sequence
\begin{equation}         
\label{eq:affPhi}
(\alpha A + \beta I; \{E_i\}_{i=0}^d; \alpha^* A^* +\beta^* I; \{E^*_i\}_{i=0}^d)
\end{equation}
is a TH system on $V$. 
\end{lemma}

\noindent {\it Proof:} Routine. \hfill $\Box$

\begin{definition}            
\label{def:aff}        
\rm
Referring to Lemma \ref{lem:affine}, we call the TH system (\ref{eq:affPhi}) the {\em affine transformation of $\Phi$ associated with $\alpha, \beta, \alpha^*, \beta^*$}.
\end{definition}

\begin{definition}          
\label{def:affiso}        
\rm 
Let $\Phi$ and $\Phi'$ denote TH systems over $\K$. We say that $\Phi$ and $\Phi'$ are {\em affine isomorphic} whenever $\Phi$ is isomorphic to an affine transformation of $\Phi'$. Observe that affine isomorphism is an equivalence relation.
\end{definition}

\begin{lemma}
\label{lem:affinepa1}
With reference to Lemma \ref{lem:affine}, let $(\lbrace \theta_i \rbrace_{i=0}^d, \lbrace \theta^*_i \rbrace_{i=0}^d, \lbrace \phi_i \rbrace_{i=1}^d)$ denote the parameter array of $\Phi$. Then the parameter array of the TH system (\ref{eq:affPhi}) is
$(\lbrace \alpha \theta_i + \beta \rbrace_{i=0}^d, \lbrace \alpha^* \theta^*_i + \beta^* \rbrace_{i=0}^d, \lbrace \alpha \alpha^* \phi_i \rbrace_{i=1}^d)$.
\end{lemma}

\noindent {\it Proof:} Let $\Phi'$ denote the TH system (\ref{eq:affPhi}). By Definition \ref{def:evseq}, for $0 \leq i \leq d$ the scalar $\theta_i$ is the eigenvalue of $A$ associated
with $E_i$, so $\alpha \theta_i + \beta$ is the eigenvalue of $\alpha A + \beta I$ associated with $E_i$. Thus $\lbrace \alpha \theta_i + \beta \rbrace_{i=0}^d$ is the eigenvalue sequence of $\Phi'$. Similarly $\lbrace \alpha^* \theta_i^* + \beta^* \rbrace_{i=0}^d$ is the dual eigenvalue sequence of $\Phi'$. 
In (\ref{eq:lab3}), if we replace $A$ by $\alpha A + \beta I$ and replace $\theta_j$ (resp. $\theta^*_j$) by $\alpha \theta_j + \beta$ (resp. $\alpha^* \theta_j^* + \beta^*$) for $0 \leq j \leq d$, then the left-hand side becomes $\alpha \alpha^* \phi_i$. Therefore $\lbrace \alpha \alpha^* \phi_i \rbrace_{i=1}^d$ is the split sequence of $\Phi'$. \hfill $\Box$

\section{The relatives of a TH system}
Let $\Phi$ denote a TH system. In the previous section we modified $\Phi$ in a certain way to get another TH system. In this section we modify $\Phi$ in a different way to obtain two more TH systems. These TH systems are called $\Phi^*$ and ${\tilde \Phi}$. We start with $\Phi^*$. 

\begin{definition}
\label{def:THdual}
\rm
Let $\Phi=(A;\lbrace E_i\rbrace_{i=0}^d;A^*;\lbrace E^*_i\rbrace_{i=0}^d)$ denote a TH system on $V$. Observe that $(A^*;\lbrace E^*_i\rbrace_{i=0}^d;A;\lbrace E_i\rbrace_{i=0}^d)$ is a TH system on $V$, which we denote by $\Phi^*$.
\end{definition}

\begin{lemma}
{\rm \cite[Lemma 6.4]{thinhess}}
\label{lem:THdualpa}
Let $\Phi$ denote a TH system with parameter array \\
$(\lbrace \theta_i \rbrace_{i=0}^d, \lbrace \theta^*_i \rbrace_{i=0}^d, \lbrace \phi_i \rbrace_{i=1}^d)$. Then the TH system $\Phi^*$ has parameter array $(\lbrace \theta^*_i \rbrace_{i=0}^d, \lbrace \theta_i \rbrace_{i=0}^d, \lbrace \phi_{d-i+1} \rbrace_{i=1}^d)$.
\end{lemma}

We now consider ${\tilde \Phi}$. For the rest of this section, let $W$ denote a vector space over $\fld$ with dimension $d+1$. For $\K$-algebras $\mathcal A$ and $\mathcal A'$, by a {\it $\fld$-algebra anti-isomorphism} from $\mathcal A$ to $\mathcal A'$ we mean a $\K$-vector space isomorphism $\dagger:\mathcal A \rightarrow \mathcal A'$ such that $(RS)^{\dagger} = S^{\dagger}R^{\dagger}$ for all $R, S \in \mathcal A$. 
By a {\it $\fld$-algebra anti-automorphism} of $\mathcal A$ we mean a $\mathbb{K}$-algebra anti-isomorphism from $\mathcal A$ to $\mathcal A$. The anti-automorphisms of $\Mdf$ are described as follows. Let $R$ denote an invertible element of $\Mdf$. Then there exists a unique $\fld$-algebra anti-automorphism $\dagger$ of $\Mdf$ such that $S^{\dagger} = RS^tR^{-1}$ for all $S \in \Mdf$. Conversely, let $\dagger$ denote a $\fld$-algebra anti-automorphism of $\Mdf$. By the Skolem-Noether theorem \cite[Corollary 9.122]{Rot}, there exists an invertible $R \in \Mdf$ such that $S^{\dagger} = RS^tR^{-1}$ for all $S \in \Mdf$. Moreover $R$ is unique up to a multiplication by a nonzero scalar in $\K$. 

Define $Z \in \Mdf$ such that $Z_{ij} = \delta_{i+j,d}$ for $0 \leq i,j \leq d$. Observe that $Z^{-1} = Z$. Define $\varsigma$ to be the $\fld$-algebra anti-automorphism of $\Mdf$ such that $S^{\varsigma} = ZS^tZ$ for all $S \in \Mdf$. For $S \in \Mdf$, $S^\varsigma$ is obtained from $S$ by reflecting about the diagonal connecting the top right corner of $S$ and the bottom left corner of $S$. In other words, $(S^{\varsigma})_{ij} = S_{d-j,d-i}$ for $0 \leq i, j \leq d$. For example, 
\begin{eqnarray*}
S =
\left(
\begin{array}{c c c}
1 & 2 & 3  \\
4 & 5 & 6  \\
7 & 8 & 9  \\
\end{array}
\right), \qquad 
S^{\varsigma} = 
\left(
\begin{array}{c c c}
9 & 6 & 3  \\
8 & 5 & 2  \\
7 & 4 & 1  \\
\end{array}
\right).
\end{eqnarray*}
\noindent Observe that $(S^{\varsigma})^{\varsigma} = S$ for all $S \in \Mdf$. Note that if $H \in \Mdf$ is Hessenberg then $H^{\varsigma}$ is Hessenberg.

\begin{lemma}
\label{lem:dualprim}
Let $A$ denote a multiplicity-free element of ${\rm End}(V)$ with eigenvalues $\{\theta_i\}_{i=0}^d$. For $0 \leq i \leq d$, let $E_i \in {\rm End}(V)$ denote the primitive idempotent of $A$ corresponding to $\theta_i$. For any anti-isomorphism $\dagger: {\rm End}(V) \rightarrow {\rm End}(W)$, the following {\rm (i)}, {\rm (ii)} hold. 
\begin{enumerate}
\item[\rm (i)]
$A^{\dagger}$ is a multiplicity-free element of ${\rm End}(W)$ with eigenvalues $\{\theta_i\}_{i=0}^d$.
\item[\rm (ii)]
For $0 \leq i \leq d$, $E_i^{\dagger}$ is the primitive idempotent of $A^{\dagger}$ corresponding to $\theta_i$.
\end{enumerate}
\end{lemma}

\noindent{\it Proof:} (i) For $f \in \mathbb{K}[\lambda]$ we have $f(A)=0$ if and only if $f(A^\dagger)=0$. Therefore $A$ and $A^\dagger$ have the same minimal polynomial. The minimal polynomial of $A$ is $\prod_{i=0}^d (\lambda - \theta_i)$ so the minimal polynomial of $A^\dagger$ is $\prod_{i=0}^d (\lambda - \theta_i)$. By this and since $\{\theta_i\}_{i=0}^d$ are mutually distinct, $A^\dagger$ is diagonalizable with eigenvalues $\{\theta_i\}_{i=0}^d$. Recall $\dim W = d+1$ so $A^\dagger$ is multiplicity-free. \\
(ii) Apply $\dagger$ to (\ref{eq:defEi}). \hfill $\Box$

\begin{proposition}  
\label{prop:antihs}
Let $(A;\lbrace E_i\rbrace_{i=0}^d;A^*;\lbrace E^*_i\rbrace_{i=0}^d)$ denote a TH system on $V$. Let $\dagger$ denote an anti-isomorphism from ${\rm End}(V)$ to ${\rm End}(W)$. Then $(A^{\dagger};\lbrace E^{\dagger}_{d-i} \rbrace_{i=0}^d; A^{*\dagger};\lbrace E^{*\dagger}_{d-i} \rbrace_{i=0}^d)$ is a TH system on $W$.
\end{proposition}

\noindent{\it Proof:} Define $\Psi = (A^{\dagger};\lbrace E^{\dagger}_{d-i} \rbrace_{i=0}^d; A^{*\dagger};\lbrace E^{*\dagger}_{d-i} \rbrace_{i=0}^d)$. In order to show that $\Psi$ is a TH system on $W$, we show that $\Psi$ satisfies conditions (i)--(v) of Definition \ref{def:HS}. By Lemma \ref{lem:dualprim}, $\Psi$ satisfies conditions (i)--(iii).  We now show that $\Psi$ satisfies condition (iv). Since $(A;\lbrace E_i\rbrace_{i=0}^d;A^*;\lbrace E^*_i\rbrace_{i=0}^d)$ is a TH system, we have
\begin{eqnarray}
\label{eq:trip1}
E_iA^*E_j = \cases{0, \ {\rm if} \; i-j > 1 \cr \not=0, \ {\rm if} \; i-j  = 1 \cr} \qquad \qquad (0 \leq i,j\leq d).
\end{eqnarray}
Applying $\dagger$ to (\ref{eq:trip1}), we find
\begin{eqnarray}
\label{eq:trip2}
E_j^\dagger A^{*\dagger} E_i^\dagger = \cases{0, \ {\rm if} \; i-j > 1 \cr \not=0, \ {\rm if} \; i-j  = 1 \cr} \qquad \qquad (0 \leq i,j\leq d).
\end{eqnarray}
Relabelling the indices in (\ref{eq:trip2}), we find 
\begin{eqnarray*}
E_{d-i}^\dagger A^{*\dagger} E_{d-j}^\dagger = \cases{0, \ {\rm if} \; i-j > 1 \cr \not=0, \ {\rm if} \; i-j  = 1 \cr}
\qquad \qquad (0 \leq i,j \leq d).
\end{eqnarray*}
Therefore $\Psi$ satisfies condition (iv). Similarly $\Psi$ satisfies condition (v). Therefore $\Psi$ is a TH system on $W$. \hfill $\Box$

\begin{definition}
\label{def:anti-isomorphicsystem}
\rm
Let $\Phi=(A;\lbrace E_i\rbrace_{i=0}^d;A^*;\lbrace E^*_i\rbrace_{i=0}^d)$ denote a TH system on $V$ and let $\Psi=(B;\lbrace F_i \rbrace_{i=0}^d;B^*;\lbrace F^*_i\rbrace_{i=0}^d)$ denote a TH system on $W$. By an {\it anti-isomorphism of TH systems} from $\Phi$ to $\Psi$ we mean a $\fld$-algebra anti-isomorphism $\dagger: {\rm End}(V) \rightarrow {\rm End}(W)$ such that 
\begin{eqnarray*}
\label{eq:antiisosys}
B = A^\dagger, \qquad  B^* = A^{*\dagger}, \qquad F_{i} = E_{d-i}^\dagger, \qquad F^*_{i} = E^{*\dagger}_{d-i} \ \ \ (0 \leq i \leq d).
\end{eqnarray*}
Observe that if $\dagger$ is an anti-isomorphism of TH systems from $\Phi$ to $\Psi$, then $\dagger^{-1}$ is an anti-isomorphism of TH systems from $\Psi$ to $\Phi$. 
We say that the TH systems $\Phi$ and $\Psi$ are {\it anti-isomorphic} whenever there exists an anti-isomorphism of TH systems from $\Phi$ to $\Psi$. 
\end{definition}

\begin{lemma}
\label{lem:antiisounique}
Let $\Phi$ denote a TH system over $\fld$. Then there exists a TH system $\Psi$ over $\fld$ such that $\Phi$ and $\Psi$ are anti-isomorphic. Moreover $\Psi$ is unique up to isomorphism of TH systems.
\end{lemma}

\noindent {\it Proof:} We first show that $\Psi$ exists. Write $\Phi=(A;\lbrace E_i\rbrace_{i=0}^d;A^*;\lbrace E^*_i\rbrace_{i=0}^d)$ and assume $V$ is the vector space underlying $\Phi$. By elementary linear algebra, there exists a $\fld$-algebra anti-automorphism $\dagger$ of ${\rm End}(V)$. Define $\Psi = (A^{\dagger};\lbrace E^{\dagger}_{d-i} \rbrace_{i=0}^d; A^{*\dagger};\lbrace E^{*\dagger}_{d-i} \rbrace_{i=0}^d)$. By Proposition \ref{prop:antihs}, $\Psi$ is a TH system on $V$. By Definition \ref{def:anti-isomorphicsystem}, $\Phi$ and $\Psi$ are anti-isomorphic. We have shown that $\Psi$ exists. Next we show that $\Psi$ is unique. Suppose that $\Psi'$ is a TH system on $W$ such that $\Phi$ and $\Psi'$ are anti-isomorphic. We show that $\Psi$ and $\Psi'$ are isomorphic. Let $\dagger'$ denote an anti-isomorphism of TH systems from $\Phi$ to $\Psi'$. Then the composition $\dagger' \dagger^{-1}$ is an isomorphism of TH systems from $\Psi$ to $\Psi'$. Therefore $\Psi$ and $\Psi'$ are isomorphic. \hfill $\Box$

\begin{lemma}
\label{lem:antiisouniquesystem}
Let $\Phi$ and $\Psi$ denote anti-isomorphic TH systems over $\K$. Then the anti-isomorphism of TH systems from $\Phi$ to $\Psi$ is unique.
\end{lemma}

\noindent {\it Proof:} Let $\dagger$ and $\dagger'$ denote anti-isomorphisms of TH systems from $\Phi$ to $\Psi$. We show that $\dagger = \dagger'$. Observe that the composition $\dagger^{-1} \dagger'$ is an isomorphism of TH systems from $\Phi$ to $\Phi$. The map $\dagger^{-1} \dagger'$ is the identity by Lemma \ref{lem:isouniquesystem}, so $\dagger = \dagger'$. \hfill $\Box$

\begin{proposition}
\label{prop:antihess}
Let $\Phi$ denote a TH system over $\fld$ with parameter array \\ $(\lbrace \theta_i \rbrace_{i=0}^d, \lbrace \theta^*_i \rbrace_{i=0}^d, \lbrace \phi_i \rbrace_{i=1}^d)$. Let $\Psi$ denote a TH system over $\fld$. Then the following {\rm (i)}, {\rm (ii)} are equivalent. 
\begin{enumerate}
 \item[\rm (i)] 
$\Phi$ and $\Psi$ are anti-isomorphic.
 \item[\rm (ii)]
The parameter array of $\Psi$ is $(\lbrace \theta_{d-i} \rbrace_{i=0}^d, \lbrace \theta^*_{d-i} \rbrace_{i=0}^d, \lbrace \phi_{d-i+1} \rbrace_{i=1}^d)$.
\end{enumerate}
\end{proposition}

\noindent {\it Proof:} (ii)$\Rightarrow$(i) Write $\Phi=(A;\lbrace E_i\rbrace_{i=0}^d;A^*;\lbrace E^*_i\rbrace_{i=0}^d)$ and $\Psi=(B;\lbrace F_i \rbrace_{i=0}^d;B^*;\lbrace F^*_i\rbrace_{i=0}^d)$. Assume $V$ (resp. $W$) is the vector space underlying $\Phi$ (resp. $\Psi$). For notational convenience, fix a $\Phi$-split basis for $V$ (resp. $\Psi$-split basis for $W$) and identify each element of ${\rm End}(V)$ (resp. ${\rm End}(W)$) with the matrix in $\Mdf$ that represents it with respect to this basis. By Proposition \ref{prop:aastarsplitrep}, 
\begin{eqnarray*}
A = \left(
\begin{array}{c c c c c c}
\theta_d & & & & & {\bf 0} \\
\phi_1 & \theta_{d-1} &  & & & \\
& \phi_2 & \theta_{d-2} &  & & \\
& & \cdot & \cdot &  &  \\
& & & \cdot & \cdot &  \\
{\bf 0}& & & & \phi_d & \theta_0
\end{array}
\right),
\qquad \quad 
A^* = \left(
\begin{array}{c c c c c c}
\theta^*_0 & 1 & & & & {\bf 0} \\
& \theta^*_1 & 1 & & & \\
& & \theta^*_2 & \cdot & & \\
& &  & \cdot & \cdot &  \\
& & &  & \cdot & 1 \\
{\bf 0}& & & &  & \theta^*_d
\end{array}
\right).
\end{eqnarray*}
Moreover 
\begin{eqnarray*}
B = \left(
\begin{array}{c c c c c c}
\theta_0 & & & & & {\bf 0} \\
\phi_d & \theta_{1} &  & & & \\
& \phi_{d-1} & \theta_{2} &  & & \\
& & \cdot & \cdot &  &  \\
& & & \cdot & \cdot &  \\
{\bf 0}& & & & \phi_1 & \theta_d
\end{array}
\right),
\qquad \quad 
B^* = \left(
\begin{array}{c c c c c c}
\theta^*_d & 1 & & & & {\bf 0} \\
& \theta^*_{d-1} & 1 & & & \\
& & \theta^*_{d-2} & \cdot & & \\
& &  & \cdot & \cdot &  \\
& & &  & \cdot & 1 \\
{\bf 0}& & & &  & \theta^*_0
\end{array}
\right).
\end{eqnarray*}
Recall the $\fld$-algebra anti-automorphism $\varsigma$ of $\Mdf$ from above Lemma \ref{lem:dualprim}. Observe that $B = A^\varsigma$ and $B^* = A^{*\varsigma}$. By this and (\ref{eq:defEi}) we find $F_{i} = E_{d-i}^\varsigma$ and $F^*_{i} = E^{*\varsigma}_{d-i}$ for $0 \leq i \leq d$. Therefore $\varsigma$ is an anti-isomorphism of TH systems from $\Phi$ to $\Psi$, so $\Phi$ and $\Psi$ are anti-isomorphic. \\
(i)$\Rightarrow$(ii) Routine by Theorem \ref{thm:classificationpa}, Lemma \ref{lem:antiisounique}, and the previous part. \hfill $\Box$

\medskip

\noindent We now discuss the notion of anti-isomorphism for TH pairs. 

\begin{definition}
\label{def:anti-isomorphicpair}
\rm
Let $A, A^*$ denote a TH pair on $V$ and let $B, B^*$ denote a TH pair on $W$. By an {\it anti-isomorphism of TH pairs} from $A, A^*$ to $B, B^*$ we mean a $\fld$-algebra anti-isomorphism $\dagger: {\rm End}(V) \rightarrow {\rm End}(W)$ such that $B = A^\dagger$ and $B^* = A^{*\dagger}$. Observe that if $\dagger$ is an anti-isomorphism of TH pairs from $A, A^*$ to $B, B^*$, then $\dagger^{-1}$ is an anti-isomorphism of TH pairs from $B, B^*$ to $A, A^*$. We say that the TH pairs $A, A^*$ and $B, B^*$ are {\it anti-isomorphic} whenever there exists an anti-isomorphism of TH pairs from $A, A^*$ to $B, B^*$. 
\end{definition}

\begin{lemma}
\label{cor:antipair}
Let $A, A^*$ denote a TH pair over $\fld$. Then there exists a TH pair $B, B^*$ over $\fld$ such that $A, A^*$ and $B, B^*$ are anti-isomorphic. Moreover $B, B^*$ is unique up to isomorphism of TH pairs.  
\end{lemma}

\noindent {\it Proof:} Similar to the proof of Lemma \ref{lem:antiisounique}. \hfill $\Box$

\begin{lemma}
\label{lem:antiisouniquepair}
Let $A, A^*$ and $B, B^*$ denote anti-isomorphic TH pairs over $\fld$. Then the anti-isomorphism of TH pairs from $A, A^*$ to $B, B^*$ is unique. 
\end{lemma}

\noindent {\it Proof:} Similar to the proof of Lemma \ref{lem:antiisouniquesystem}. \hfill $\Box$

\medskip

We recall some more terms and facts from elementary linear algebra. A map $\b{\;,\,} : V \times W \to \mathbb{K}$ is called a {\em bilinear form} whenever the following conditions hold
for all $v, v' \in V$, $w, w' \in W$, and $\alpha \in \mathbb{K}$:
(i) $\b{v + v', w}=\b{v,w}+\b{v',w}$;
(ii) $\b{\alpha v, w}=\alpha \b{v,w}$;
(iii) $\b{v, w+w'}=\b{v,w}+\b{v,w'}$;
(iv) $\b{v,\alpha w}=\alpha \b{v,w}$.
We observe that a scalar multiple of a bilinear form is a bilinear form. 

Let $\b{\;,\,} : V \times W \to \mathbb{K}$ denote a bilinear form. Then the following are equivalent:
(i) there exists a nonzero $v \in V$ such that $\b{v,w}=0$ for all $w\in W$;
(ii) there exists a nonzero $w \in W$ such that $\b{v,w}=0$ for all $v\in V$.
The form $\b{\;,\,}$ is said to be {\em degenerate} whenever (i), (ii) hold and {\em nondegenerate} otherwise.

Bilinear forms are related to anti-isomorphisms as follows. Let $\b{\;,\,} : V\times W \to \mathbb{K}$ denote a nondegenerate
bilinear form. Then there exists a unique anti-isomorphism $\dagger :{\rm End}(V) \to {\rm End}(W)$ such that $\b{Sv,w}=\b{v,S^\dagger w}$ for all  $v \in V$, $w \in W$, and $S \in {\rm End}(V)$. 
Conversely, given an anti-isomorphism $\dagger :{\rm End}(V) \to {\rm End}(W)$ there exists a nonzero bilinear form $\b{\;,\,} : V \times W \to \mathbb{K}$ such that $\b{Sv,w}=\b{v,S^\dagger w}$ for all $v \in V$, $w \in W$, and $S \in {\rm End}(V)$. This bilinear form is nondegenerate, and uniquely determined by $\dagger$ up to multiplication by a nonzero scalar in $\mathbb{K}$.
We say that the form $\b{\;,\,}$ is {\em associated} with $\dagger$.

Define ${\tilde V}$ to be the dual space of $V$, consisting of all $\mathbb{K}$-linear transformations from $V$ to $\mathbb{K}$. By elementary linear algebra, ${\tilde V}$ is a vector space over $\K$ and $\dim {\tilde V} = \dim V$. Define a bilinear form $\b{\;,\,}: V \times \tilde{V} \to \mathbb{K}$ such that $\b{v,f}=f(v)$ for all $v \in V$ and $f \in \tilde{V}$. The form $\b{\;,\,}$ is nondegenerate. We call $\b{\;,\,}$ the {\em canonical bilinear form} between $V$ and $\tilde{V}$. Let $\sigma : {\rm End}(V) \to {\rm End}(\tilde{V})$ denote the anti-isomorphism associated with $\b{\;,\,}$. Thus $\b{Sv, f}=\b{v,S^\sigma f}$ for all $v \in V$, $f \in \tilde{V}$, and $S \in {\rm End}(V)$. We call $\sigma$ the {\em canonical anti-isomorphism} from ${\rm End}(V)$ to ${\rm End}(\tilde{V})$.

\begin{definition}
\label{def:dualphisetup}
\rm
Let $\Phi=(A;\lbrace E_i\rbrace_{i=0}^d;A^*;\lbrace E^*_i\rbrace_{i=0}^d)$ denote a TH system on $V$ with parameter array $(\lbrace \theta_i \rbrace_{i=0}^d, \lbrace \theta^*_i \rbrace_{i=0}^d, \lbrace \phi_i \rbrace_{i=1}^d)$. Define ${\tilde \Phi} =  (A^\sigma;\lbrace E^{\sigma}_{d-i} \rbrace_{i=0}^d; A^{*\sigma};\lbrace E^{*\sigma}_{d-i} \rbrace_{i=0}^d)$, where $\sigma: {\rm End}(V) \to {\rm End}(\tilde{V})$ is the canonical anti-isomorphism. By Proposition \ref{prop:antihs}, ${\tilde \Phi}$ is a TH system on ${\tilde V}$. By Definition \ref{def:anti-isomorphicsystem}, $\Phi$ and ${\tilde \Phi}$ are anti-isomorphic. By Proposition \ref{prop:antihess}, ${\tilde \Phi}$ has parameter array $(\lbrace \theta_{d-i} \rbrace_{i=0}^d, \lbrace \theta^*_{d-i} \rbrace_{i=0}^d, \lbrace \phi_{d-i+1} \rbrace_{i=1}^d)$. 
\end{definition}

\section{The ${\mathbb Z_2} \times {\mathbb Z_2}$ action}
\label{sec:z2action}
Let $\Phi=(A;\lbrace E_i\rbrace_{i=0}^d;A^*;\lbrace E^*_i\rbrace_{i=0}^d)$ denote a TH system. We saw in the previous section that each of the following is a TH system:
\begin{eqnarray*}
\Phi^{*}  &=& (A^*; \{E^*_i\}_{i=0}^d; A; \{E_i\}_{i=0}^d), \\
{\tilde \Phi} &=& (A^\sigma; \{E^\sigma_{d-i}\}_{i=0}^d; A^{*\sigma}; \{E^{*\sigma}_{d-i}\}_{i=0}^d).
\end{eqnarray*}
Viewing $*$, ${\sim}$ as permutations on the set of all TH systems,
\begin{eqnarray}
\label{eq:relation1}    
*^2\;=\;\sim^2 \;=\;1, \qquad \qquad \ast \sim \;=\; \sim *.
\end{eqnarray}
The group generated by symbols $*$, $\sim$ subject to the relations (\ref{eq:relation1}) is the group ${\mathbb Z_2} \times {\mathbb Z_2}$. Thus $*$, $\sim$ induce an action of ${\mathbb Z_2} \times {\mathbb Z_2}$ on the set of all TH systems.
Two TH systems will be called {\em relatives} whenever they are in the same orbit of this ${\mathbb Z_2} \times {\mathbb Z_2}$ action. 
The relatives of $\Phi$ are as follows:

\medskip
\noindent
\begin{center}
\begin{tabular}{c|c}
name & relative \\
\hline
$\Phi$ & $(A; \{E_i\}_{i=0}^d; A^*;  \{E^*_i\}_{i=0}^d)$  \\ 

$\Phi^{*}$  & $(A^*; \{E^*_i\}_{i=0}^d; A;  \{E_i\}_{i=0}^d)$  \\
 
${\tilde \Phi}$ & $(A^\sigma; \{E^\sigma_{d-i}\}_{i=0}^d; A^{*\sigma}; \{E^{*\sigma}_{d-i}\}_{i=0}^d)$ \\ 

${\tilde \Phi^{*}}$ & $(A^{*\sigma}; \{E^{*\sigma}_{d-i}\}_{i=0}^d; A^\sigma; \{E^\sigma_{d-i}\}_{i=0}^d)$  
\end{tabular}
\end{center}

\begin{corollary}
\label{cor:parelatives}
Let $\Phi$ denote a TH system with parameter array $(\lbrace \theta_i \rbrace_{i=0}^d, \lbrace \theta^*_i \rbrace_{i=0}^d, \lbrace \phi_i \rbrace_{i=1}^d)$. Then the parameter arrays of its relatives are as follows:

\medskip
\noindent
\begin{center}
\begin{tabular}{c|c}
name & parameter array \\
\hline
$\Phi$ &  $(\lbrace \theta_i \rbrace_{i=0}^d, \lbrace \theta^*_i \rbrace_{i=0}^d, \lbrace \phi_i \rbrace_{i=1}^d)$ \\ 

$\Phi^{*}$  & $(\lbrace \theta^*_i \rbrace_{i=0}^d, \lbrace \theta_i \rbrace_{i=0}^d, \lbrace \phi_{d-i+1} \rbrace_{i=1}^d)$ \\
 
${\tilde \Phi}$ & $(\lbrace \theta_{d-i} \rbrace_{i=0}^d, \lbrace \theta^*_{d-i} \rbrace_{i=0}^d, \lbrace \phi_{d-i+1} \rbrace_{i=1}^d)$ \\ 

${\tilde \Phi^{*}}$ & $(\lbrace \theta^*_{d-i} \rbrace_{i=0}^d, \lbrace \theta_{d-i} \rbrace_{i=0}^d, \lbrace \phi_i \rbrace_{i=1}^d)$ 
\end{tabular}
\end{center}
\end{corollary}

\noindent {\it Proof:} Immediate from Lemma \ref{lem:THdualpa} and Proposition \ref{prop:antihess}. \hfill $\Box$

\medskip

\noindent We will use the following notational convention.

\smallskip

\begin{definition}
\label{def:conv2}
\rm
Let $\Phi$ denote a TH system. For $g \in {\mathbb Z_2} \times {\mathbb Z_2}$ and for an object $f$ associated with $\Phi$, let $f^g$ denote the corresponding object associated with $\Phi^{g}$.
\end{definition}

\noindent We end this section by giving some more characterizations of the split sequence, as promised at the end of Section \ref{sec:thsystem}.

\begin{lemma}
\label{lem:splitchar2}
Let $(A;\lbrace E_i\rbrace_{i=0}^d;A^*;\lbrace E^*_i\rbrace_{i=0}^d)$ denote a TH system with parameter array \\ $(\lbrace \theta_i \rbrace_{i=0}^d, \lbrace \theta^*_i \rbrace_{i=0}^d, \lbrace \phi_i \rbrace_{i=1}^d)$. Then the following {\rm (i)--(iii)} hold for $0 \leq i \leq d$. 
\begin{enumerate}
 \item[\rm (i)]
${\displaystyle  E_0 \eta^*_i(A^*) E_0 = \frac{\phi_d \phi_{d-1} \cdots \phi_{d-i+1}}{(\theta_0 - \theta_1)(\theta_0 - \theta_2) \cdots (\theta_0 - \theta_i)} E_0.}$
 \item[\rm (ii)]
${\displaystyle  E_d^* \tau_i(A) E_d^* = \frac{\phi_d \phi_{d-1} \cdots \phi_{d-i+1}}{(\theta^*_d - \theta^*_{d-1})(\theta^*_d - \theta^*_{d-2}) \cdots (\theta^*_d - \theta^*_{d-i})} E^*_d.}$
 \item[\rm (iii)]
${\displaystyle  E_d \tau^*_i(A^*) E_d = \frac{\phi_1 \phi_2 \cdots \phi_{i}}{(\theta_d - \theta_{d-1})(\theta_d - \theta_{d-2}) \cdots (\theta_d - \theta_{d-i})} E_d.}$
\end{enumerate}
\end{lemma}

\noindent {\it Proof:} Let $\Phi$ denote the TH system in question.  \\
(i) Apply Lemma \ref{lem:splitchar} to $\Phi^*$. \\
(ii) Apply Lemma \ref{lem:splitchar} to ${\tilde \Phi}$ and then apply $\sigma^{-1}$ to each side of the resulting equations. \\
(iii) Apply (ii) to $\Phi^*$. \hfill $\Box$

\begin{corollary}
\label{cor:splitchar2}
Let $(A;\lbrace E_i\rbrace_{i=0}^d;A^*;\lbrace E^*_i\rbrace_{i=0}^d)$ denote a TH system with parameter array \\ $(\lbrace \theta_i \rbrace_{i=0}^d, \lbrace \theta^*_i \rbrace_{i=0}^d, \lbrace \phi_i \rbrace_{i=1}^d)$. Then the following {\rm (i)--(iii)} hold for $0 \leq i \leq d$. 
\begin{enumerate}
 \item[\rm (i)]
${\displaystyle \phi_d \phi_{d-1} \cdots \phi_{d-i+1} = (\theta_0 - \theta_1)(\theta_0 - \theta_2) \cdots (\theta_0 - \theta_i) {\rm trace}(\eta^*_i(A^*) E_0).}$
 \item[\rm (ii)]
${\displaystyle \phi_d \phi_{d-1} \cdots \phi_{d-i+1} = (\theta^*_d - \theta^*_{d-1})(\theta^*_d - \theta^*_{d-2}) \cdots (\theta^*_d - \theta^*_{d-i}) {\rm trace}(\tau_i(A) E^*_d).}$
 \item[\rm (iii)]
${\displaystyle \phi_1 \phi_2 \cdots \phi_i = (\theta_d - \theta_{d-1})(\theta_d - \theta_{d-2}) \cdots (\theta_d - \theta_{d-i}) {\rm trace}(\tau^*_i(A^*) E_d).}$
\end{enumerate}
Moreover each of $\eta^*_i(A^*) E_0, \tau_i(A) E^*_d, \tau^*_i(A^*) E_d$ has nonzero trace. 
\end{corollary}

\noindent {\it Proof:} Let $\Phi$ denote the TH system in question. \\
(i) Apply Corollary \ref{cor:splitchar} to $\Phi^*$. \\
(ii) In the equation of Lemma \ref{lem:splitchar2}(ii), take the trace of each side and simplify the result using the fact that ${\rm trace}(E_d^*) = 1$ and ${\rm trace}(E_d^* \tau_i(A) E_d^*) = {\rm trace}(\tau_i(A)E_d^*E_d^*) = {\rm trace}( \tau_i(A)E_d^*)$.  \\
(iii) Apply (ii) to $\Phi^*$. \\ 
The last assertion follows since $\phi_i \neq 0$ for $1 \leq i \leq d$. \hfill $\Box$

\begin{corollary}
Let $(A;\lbrace E_i\rbrace_{i=0}^d;A^*;\lbrace E^*_i\rbrace_{i=0}^d)$ denote a TH system with parameter array \\ $(\lbrace \theta_i \rbrace_{i=0}^d, \lbrace \theta^*_i \rbrace_{i=0}^d, \lbrace \phi_i \rbrace_{i=1}^d)$. Then the following {\rm (i)--(iii)} hold for $1 \leq i \leq d$. 
\begin{enumerate}
 \item[\rm (i)]
${\displaystyle \phi_{i} = (\theta_0 - \theta_{d-i+1}) {\rm trace}(\eta^*_{d-i+1}(A^*) E_0)/{\rm trace}(\eta^*_{d-i}(A^*) E_0).}$
 \item[\rm (ii)]
${\displaystyle \phi_{i} = (\theta^*_d - \theta^*_{i-1}) {\rm trace}(\tau_{d-i+1}(A) E^*_d)/{\rm trace}(\tau_{d-i}(A) E^*_d).}$
 \item[\rm (iii)]
${\displaystyle \phi_i = (\theta_d - \theta_{d-i}) {\rm trace}(\tau^*_i(A^*) E_d)/{\rm trace}(\tau^*_{i-1}(A^*) E_d).}$
\end{enumerate}
\end{corollary}

\noindent {\it Proof:}  Routine by Corollary \ref{cor:splitchar2}.  \hfill $\Box$

\section{The scalars $\lbrace \ell_i \rbrace_{i=0}^d$}
Let $\Phi$ denote a TH system. In this section we associate with $\Phi$ a sequence of scalars $\lbrace \ell_i \rbrace_{i=0}^d$ that will help us describe $\Phi$. 
 
\begin{definition}
\label{def:kappas}
\rm
Let $\Phi$ denote a TH system with dual eigenvalue sequence $\lbrace \theta^*_i \rbrace_{i=0}^d$. For $0 \leq i \leq d$, define   
\begin{eqnarray*}
 \ell_i &=& \frac{\eta^*_d(\theta^*_0)}{\tau^*_i(\theta^*_i)\eta^*_{d-i}(\theta^*_i)} \\
     &=& \frac{(\theta^*_0-\theta^*_1)(\theta^*_0-\theta^*_2)\cdots (\theta^*_0-\theta^*_d)}{(\theta^*_i-\theta^*_0)(\theta^*_i-\theta^*_1)\cdots (\theta^*_i-\theta^*_{i-1}) (\theta^*_i-\theta^*_{i+1}) \cdots (\theta^*_i - \theta^*_{d-1})(\theta^*_i-\theta^*_d)}. \\
\end{eqnarray*}
Observe that $\ell_0 = 1$. 
\end{definition}

\begin{lemma}
\label{lem:ellistar}
Let $\Phi$ denote a TH system with eigenvalue sequence $\lbrace \theta_i \rbrace_{i=0}^d$ and dual eigenvalue sequence $\lbrace \theta^*_i \rbrace_{i=0}^d$. Then for $0 \leq i \leq d$, 
\begin{eqnarray*}
\ell_i^*&=& \frac{\eta_d(\theta_0)}{\tau_i(\theta_i)\eta_{d-i}(\theta_i)} \\
&=& \frac{(\theta_0-\theta_1)(\theta_0-\theta_2) \cdots (\theta_0-\theta_d)}{(\theta_i-\theta_0)(\theta_i-\theta_1)\cdots (\theta_i-\theta_{i-1})(\theta_i-\theta_{i+1}) \cdots (\theta_i- \theta_{d-1})(\theta_i-\theta_d)}, 
\\
\\
{\tilde \ell}_i &=& \frac{\tau^*_d(\theta^*_d)}{\eta^*_i(\theta^*_{d-i})\tau^*_{d-i}(\theta^*_{d-i})} \\
  &=& \frac{(\theta^*_d-\theta^*_{d-1})(\theta^*_d-\theta^*_{d-2}) \cdots (\theta^*_d-\theta^*_{0})}
 {(\theta^*_{d-i}-\theta^*_d)(\theta^*_{d-i}-\theta^*_{d-1})\cdots (\theta^*_{d-i}-\theta^*_{d-i+1})
 (\theta^*_{d-i}-\theta^*_{d-i-1}) \cdots (\theta^*_{d-i} - \theta^*_{1})(\theta^*_{d-i}-\theta^*_0)},
\\
\\
{\tilde \ell_i^*} &=& \frac{\tau_d(\theta_d)}{\eta_i(\theta_{d-i})\tau_{d-i}(\theta_{d-i})} \\ 
 &=& \frac{(\theta_d-\theta_{d-1})(\theta_d-\theta_{d-2})\cdots (\theta_d-\theta_0)}
{(\theta_{d-i}-\theta_d)(\theta_{d-i}-\theta_{d-1})\cdots (\theta_{d-i}-\theta_{d-i+1})
(\theta_{d-i}-\theta_{d-i-1}) \cdots (\theta_{d-i} - \theta_{1})(\theta_{d-i}-\theta_0)}.
 \end{eqnarray*}
Moreover $\displaystyle{{\tilde \ell}_i = \frac{\tau^*_d(\theta^*_d)}{\eta^*_d(\theta^*_0)} \ell_{d-i}}$ and 
$\displaystyle{{\tilde \ell^*}_i = \frac{\tau_d(\theta_d)}{\eta_d(\theta_0)} \ell^*_{d-i}}$.
\end{lemma}

\noindent {\it Proof:} Combine Corollary \ref{cor:parelatives} and Definition \ref{def:kappas}. \hfill $\Box$

\medskip

\noindent We give one significance of the sequence $\lbrace \ell_i \rbrace_{i=0}^d$.

\begin{lemma}
\label{lem:bilformstdual}
Let $(A;\lbrace E_i\rbrace_{i=0}^d;A^*;\lbrace E^*_i\rbrace_{i=0}^d)$ denote a TH system. Then $E_dE_i^*E_0 = \ell_i E_dE_0^*E_0$ for $0 \leq i \leq d$.
\end{lemma}

\noindent {\it Proof:} Let $\Phi$ denote the TH system in question and assume $V$ is the underlying vector space. For notational convenience, fix a $\Phi$-split basis for $V$ and identify each element of ${\rm End}(V)$ with the matrix in $\Mdf$ that represents it with respect to this basis. We show that $E_d(E_i^* - \ell_i E_0^*)E_0=0$. By (\ref{eq:entriesofefreefinS99}) the entries of all but the first column of $E_d$ are zero and the entries of all but the last row of $E_0$ are zero. Therefore for $0 \leq m, n \leq d$,
the $(m,n)$-entry of $E_d(E_i^* - \ell_i E_0^*)E_0$ is 
\begin{eqnarray}
\label{eq:lab33}
(E_d)_{m0}(E_i^* - \ell_i E_0^*)_{0d}(E_0)_{dn}.  
\end{eqnarray}
By (\ref{eq:theestarentriesfinS99}) the middle factor in (\ref{eq:lab33}) is $0$, so (\ref{eq:lab33}) is $0$. Therefore $E_d(E_i^* - \ell_i E_0^*)E_0=0$ and the result follows. \hfill $\Box$


\begin{corollary}
\label{cor:bilformststdual2}
Let $(A;\lbrace E_i\rbrace_{i=0}^d;A^*;\lbrace E^*_i\rbrace_{i=0}^d)$ denote a TH system. Then the following {\rm (i)--(iii)} hold for $0 \leq i \leq d$. 
\begin{enumerate}
 \item[\rm (i)]
$E_d^*E_iE_0^* = \ell_i^* E_d^*E_0E_0^*$. 
 \item[\rm (ii)] 
$E_dE_i^*E_0 = {\tilde \ell_{d-i}} E_dE_d^*E_0$.
 \item[\rm (iii)]
$E^*_dE_iE^*_0 = {\tilde \ell^*_{d-i}} E^*_dE_dE^*_0$.
\end{enumerate}
\end{corollary}

\noindent {\it Proof:} Let $\Phi$ denote the TH system in question. \\
(i) Apply Lemma \ref{lem:bilformstdual} to $\Phi^*$. \\
(ii) Apply Lemma \ref{lem:bilformstdual} to ${\tilde \Phi}$ and then apply $\sigma^{-1}$ to each side of the resulting equations.  \\
(iii) Apply (ii) to $\Phi^*$. \hfill $\Box$

\medskip

\begin{definition}
\label{def:diagli}
\rm
Let $\Phi$ denote a TH system of diameter $d$. We associate with $\Phi$ a diagonal matrix $L \in \Mdf$ with $(i,i)$-entry $\ell_i$ for $0 \leq i \leq d$.  
\end{definition}

\section{The scalar $\nu$}
Let $\Phi$ denote a TH system. In this section we associate with $\Phi$ a scalar $\nu$ that will help us describe $\Phi$. 

\begin{definition}
\label{def:nu}
\rm
Let $(A;\lbrace E_i\rbrace_{i=0}^d;A^*;\lbrace E^*_i\rbrace_{i=0}^d)$ denote a TH system. By \cite[Lemma 7.5]{thinhess}, ${\rm trace}(E_0E^*_0)$ is nonzero. Let $\nu$ denote the reciprocal of ${\rm trace}(E_0E^*_0)$.
\end{definition}

\noindent We give one significance of the scalar $\nu$.

\begin{lemma}
\label{lem:nusig}
{\rm \cite[Lemma 7.4]{thinhess}}
Let $(A;\lbrace E_i\rbrace_{i=0}^d;A^*;\lbrace E^*_i\rbrace_{i=0}^d)$ denote a TH system. Then both
\begin{eqnarray*}
\nu E_0E^*_0E_0 = E_0, \qquad \qquad \nu E^*_0E_0E^*_0 = E^*_0. 
\end{eqnarray*}
\end{lemma}

\begin{lemma}
\label{lem:nupa}
Let $\Phi$ denote a TH system with parameter array $(\lbrace \theta_i \rbrace_{i=0}^d, \lbrace \theta^*_i \rbrace_{i=0}^d, \lbrace \phi_i \rbrace_{i=1}^d)$ and let $\nu$ denote the scalar from Definition \ref{def:nu}. Then 
\begin{eqnarray}
\label{eq:lab38}
\nu &=& \frac{(\theta_0-\theta_1)(\theta_0-\theta_2)\cdots(\theta_0-\theta_d)(\theta^*_0-\theta^*_1)(\theta^*_0-\theta^*_2) \cdots (\theta^*_0-\theta^*_d)}{\phi_1 \phi_2 \cdots \phi_d}, \\
{\tilde \nu} &=& \frac{(\theta_d-\theta_{d-1})(\theta_d-\theta_{d-2})\cdots(\theta_d-\theta_0)(\theta^*_d-\theta^*_{d-1})(\theta^*_d-\theta^*_{d-2}) \cdots (\theta^*_d-\theta^*_0)}{\phi_1 \phi_2 \cdots \phi_d}.
\nonumber
\end{eqnarray}
Moreover $\nu^* = \nu$ and ${\tilde \nu^*} = {\tilde \nu}$. 
\end{lemma}

\noindent {\it Proof:} Line (\ref{eq:lab38}) holds by \cite[Lemma 7.6]{thinhess}. The remaining assertions follow from Corollary \ref{cor:parelatives}.  \hfill $\Box$

\begin{lemma}
Let $(A;\lbrace E_i\rbrace_{i=0}^d;A^*;\lbrace E^*_i\rbrace_{i=0}^d)$ denote a TH system. Then both
\begin{eqnarray*}
{\tilde \nu} E_d  E^*_d E_d = E_d, \qquad \qquad {\tilde \nu} E_d^* E_d E_d^* = E_d^*.
\end{eqnarray*}
\end{lemma}

\noindent {\it Proof:} Let $\Phi$ denote the TH system in question. Apply Lemma \ref{lem:nusig} to ${\tilde \Phi}$ and then apply $\sigma^{-1}$ to each side of the resulting equations. \hfill $\Box$

\section{Anti-isomorphic TH systems and the bilinear form}
\label{sec:bilform}
Let $\Phi$ denote a TH system on $V$ and let ${\tilde \Phi}$ denote the relative of $\Phi$  from Definition \ref{def:dualphisetup}. Recall that $\Phi$ and ${\tilde \Phi}$ are anti-isomorphic.  In this section, we discuss further the relationship between $\Phi$ and ${\tilde \Phi}$. 
Recall the canonical bilinear form $\b{\;,\,}: V \times {\tilde V} \to \K$ from above Definition \ref{def:dualphisetup}.
Let $U$ (resp. ${\tilde U}$) denote a subspace of $V$ (resp. ${\tilde V}$). We say that $U$ and ${\tilde U}$ are {\it orthogonal} whenever $\b{x, y} = 0$ for all $x \in U$ and $y \in {\tilde U}$.
Let $\{ V_i \}_{i=0}^d$ (resp. $\{ {\tilde V_i} \}_{i=0}^d$) denote a decomposition of $V$ (resp. ${\tilde V}$). We say that $\{ V_i \}_{i=0}^d$ and $\{ {\tilde V_i} \}_{i=0}^d$ are {\it dual} whenever $V_i$ and ${\tilde V_j}$ are orthogonal for $0 \leq i,j \leq d$, $i \neq j$.  
By elementary linear algebra, for any decomposition $\{ V_i \}_{i=0}^d$ (resp. $\{ {\tilde V_i} \}_{i=0}^d$) of $V$ (resp. ${\tilde V}$) there exists a unique decomposition $\{{\tilde V_i} \}_{i=0}^d$ (resp. $\{ V_i \}_{i=0}^d$) of ${\tilde V}$ (resp. $V$) such that $\{ V_i \}_{i=0}^d$ and $\{ {\tilde V_i} \}_{i=0}^d$ are dual. 
Let $\{ v_i \}_{i=0}^d$ (resp. $\{{\tilde v_i} \}_{i=0}^d$) denote a basis for $V$ (resp. ${\tilde V}$). We say that $\{ v_i \}_{i=0}^d$ and $\{ {\tilde v_i} \}_{i=0}^d$ are {\it dual}  whenever $\b{v_i, {\tilde v_j}} = \delta_{ij}$ for $0 \leq i,j \leq d$. 
By elementary linear algebra, for any basis $\{ v_i \}_{i=0}^d$ (resp. $\{{\tilde v_i} \}_{i=0}^d$) for $V$ (resp. ${\tilde V}$) there exists a unique basis $\{{\tilde v_i} \}_{i=0}^d$ (resp. $\{ v_i \}_{i=0}^d$) for ${\tilde V}$ (resp. $V$) such that $\{ v_i \}_{i=0}^d$ and $\{ {\tilde v_i} \}_{i=0}^d$ are dual.
Given any sequence $\lbrace \alpha_i \rbrace_{i=0}^d$, by the {\it inversion} of $\lbrace \alpha_i \rbrace_{i=0}^d$ we mean the sequence $\lbrace \alpha_{d-i} \rbrace_{i=0}^d$. 

Recall the $\Phi$-standard decomposition $\lbrace E^*_iV \rbrace_{i=0}^d$ from above (\ref{eq:defui}). Observe by Definition \ref{def:dualphisetup} that $\lbrace E^{*\sigma}_{d-i} {\tilde V} \rbrace_{i=0}^d$ is the ${\tilde \Phi}$-standard decomposition. We now compare these two decompositions. 

\begin{lemma}
\label{lem:dualdecst2}
With reference to Definition \ref{def:dualphisetup}, the following {\rm (i)}, {\rm (ii)} are inverted dual.
\begin{enumerate}
 \item[\rm (i)]
The $\Phi$-standard decomposition of $V$.
 \item[\rm (ii)]
The ${\tilde \Phi}$-standard decomposition of ${\tilde V}$.
\end{enumerate}
\end{lemma}

\noindent {\it Proof:} For distinct $i, j \ (0 \leq i, j \leq d)$ we show that $E_i^*V$ and $E_j^{*\sigma}{\tilde V}$ are orthogonal. 
Let $u \in E_i^*V$ and $v \in E_j^{*\sigma} {\tilde V}$. Simplify the equation $\langle A^*u, v \rangle = \langle u, A^{*\sigma} v \rangle$ using  $A^*u = \theta_i^* u$ and $A^{*\sigma} v = \theta_j^* v$ to obtain $(\theta_i^* - \theta_j^*)\langle u, v \rangle = 0$. Now $\langle u, v \rangle = 0$ since $\theta_i^* \neq \theta_j^*$. Therefore $E_i^*V$ and $E_j^{*\sigma} {\tilde V}$ are orthogonal and the result follows. \hfill $\Box$

\medskip

Let $0 \neq \xi_0 \in E_0V$ and recall the $\Phi$-standard basis $\lbrace E^*_i \xi_0 \rbrace_{i=0}^d$ for $V$ from above (\ref{eq:defui}).  
Let $0 \neq {\tilde \xi_d} \in E^{\sigma}_d{\tilde V}$ and observe by Definition \ref{def:dualphisetup} that $\lbrace E^{*\sigma}_{d-i} {\tilde \xi_d} \rbrace_{i=0}^d$ is a ${\tilde \Phi}$-standard basis for ${\tilde V}$. These two bases are related as follows. 

\begin{proposition}
\label{prop:orthog}
With reference to Definition \ref{def:dualphisetup}, let $0 \neq \xi_0 \in E_0V$ and $0 \neq {\tilde \xi_d} \in E^{\sigma}_d{\tilde V}$. Then for $0 \leq i,j \leq d$, 
\begin{eqnarray*}
\langle E^*_i\xi_0, E^{*\sigma}_j {\tilde \xi_d} \rangle = \delta_{ij} \ell_i  \b{E_0^*\xi_0, {\tilde \xi_d}}, 
\end{eqnarray*}
where $\ell_i$ is from Definition \ref{def:kappas}.
\end{proposition}

\noindent {\it Proof:} Using the definition of $\sigma$ from above Definition \ref{def:dualphisetup} along with Lemma \ref{lem:bilformstdual}, we find
\begin{eqnarray*}
\langle E^*_i\xi_0, E^{*\sigma}_j {\tilde \xi_d}                       \rangle
&=&
\langle E^*_iE_0\xi_0, E^{*\sigma}_jE^{\sigma}_d {\tilde \xi_d}        \rangle \\
&=&
\langle E^*_jE^*_iE_0\xi_0,E^{\sigma}_d {\tilde \xi_d}                \rangle \\
&=&
\delta_{ij} \langle E^*_iE_0\xi_0,E^{\sigma}_d {\tilde \xi_d}         \rangle \\
&=&
\delta_{ij} \langle E_dE^*_iE_0\xi_0, {\tilde \xi_d}                   \rangle \\
&=&
\delta_{ij} \ell_i \langle E_dE^*_0E_0 \xi_0, {\tilde \xi_d}          \rangle  \\
&=&
\delta_{ij} \ell_i \langle E^*_0E_0 \xi_0, E^{\sigma}_d {\tilde \xi_d}  \rangle  \\
&=&
\delta_{ij} \ell_i \langle E^*_0 \xi_0, {\tilde \xi_d} \rangle.
\end{eqnarray*}
\hfill $\Box $ 

\begin{corollary}
\label{cor:dualvstilde}
With reference to Definition \ref{def:dualphisetup}, let $\lbrace v_i \rbrace_{i=0}^d$ (resp. $\lbrace w_i \rbrace_{i=0}^d$) denote a basis for $V$ (resp. $\tilde V$). Suppose that $\lbrace v_i \rbrace_{i=0}^d$ and $\lbrace w_i \rbrace_{i=0}^d$ are inverted dual. Then the following {\rm (i)}, {\rm (ii)} are equivalent.
\begin{enumerate}
 \item[\rm (i)] 
$\lbrace \ell_i v_i \rbrace_{i=0}^d$ is a $\Phi$-standard basis for $V$.
 \item[\rm (ii)]
$\lbrace w_{i} \rbrace_{i=0}^d$ is a ${\tilde \Phi}$-standard basis for ${\tilde V}$.
\end{enumerate}
\end{corollary}

\noindent {\it Proof:} Use Proposition \ref{prop:orthog}. \hfill $\Box$

\medskip

By Definition \ref{def:THdual}, the sequence $\lbrace E_iV \rbrace_{i=0}^d$ is the $\Phi^*$-standard decomposition. By Definition \ref{def:dualphisetup}, the sequence $\lbrace E^{\sigma}_{d-i} {\tilde V} \rbrace_{i=0}^d$ is the ${\tilde \Phi^*}$-standard decomposition. We now compare these two decompositions. 

\begin{lemma}
\label{lem:dualdecst}
With reference to Definition \ref{def:dualphisetup}, the following {\rm (i)}, {\rm (ii)} are inverted dual.
\begin{enumerate}
 \item[\rm (i)]
The $\Phi^*$-standard decomposition of $V$.
 \item[\rm (ii)]
The ${\tilde \Phi^*}$-standard decomposition of ${\tilde V}$.
\end{enumerate}
\end{lemma}

\noindent {\it Proof:} Apply Lemma \ref{lem:dualdecst2} to $\Phi^*$. \hfill $\Box$ 

\medskip

Let $0 \neq \xi^*_0 \in E^*_0V$ and observe by Definition \ref{def:THdual} that $\lbrace E_i \xi^*_0 \rbrace_{i=0}^d$ is a $\Phi^*$-standard basis for $V$.
Let $0 \neq {\tilde \xi^*_d} \in E^{*\sigma}_d{\tilde V}$ and observe by Definition \ref{def:dualphisetup} that $\lbrace E^{\sigma}_{d-i} {\tilde \xi^*_d} \rbrace_{i=0}^d$ is a ${\tilde \Phi^*}$-standard basis for ${\tilde V}$. These two bases are related as follows. 

\begin{proposition}
\label{prop:orthogstar}
With reference to Definition \ref{def:dualphisetup}, let $0 \neq \xi^*_0 \in E^*_0V$ and $0 \neq {\tilde \xi^*_d} \in E^{*\sigma}_d{\tilde V}$. Then for $0 \leq i,j \leq d$, 
\begin{eqnarray*}
\langle E_i\xi^*_0,E^{\sigma}_j {\tilde \xi^*_d} \rangle = \delta_{ij} \ell^*_i  \b{E_0 \xi^*_0, {\tilde \xi^*_d}}, 
\end{eqnarray*}
where $\ell^*_i$ is from Lemma \ref{lem:ellistar}.
\end{proposition}

\noindent {\it Proof:} Apply Proposition \ref{prop:orthog} to $\Phi^*$. \hfill $\Box$ 

\begin{corollary}
\label{cor:dualvstildestar}
With reference to Definition \ref{def:dualphisetup}, let $\lbrace v_i \rbrace_{i=0}^d$ (resp. $\lbrace w_i \rbrace_{i=0}^d$) denote a basis for $V$ (resp. $\tilde V$). Suppose that $\lbrace v_i \rbrace_{i=0}^d$ and $\lbrace w_i \rbrace_{i=0}^d$ are inverted dual. Then the following {\rm (i)}, {\rm (ii)} are equivalent.
\begin{enumerate}
 \item[\rm (i)]
$\lbrace \ell^*_{i} v_i \rbrace_{i=0}^d$ is a $\Phi^*$-standard basis for $V$.
 \item[\rm (ii)]
$\lbrace w_i \rbrace_{i=0}^d$ is a ${\tilde \Phi^*}$-standard basis for ${\tilde V}$.
\end{enumerate}
\end{corollary}

\noindent {\it Proof:} Apply Corollary \ref{cor:dualvstilde} to $\Phi^*$. \hfill $\Box$

\medskip

\noindent Let $\{U_i \}_{i=0}^d$ denote the $\Phi$-split decomposition of $V$. Recall from (\ref{eq:defui}) that for $0 \leq i \leq d$,
\begin{equation}
\label{eq:defui2}
U_i = (E^*_0V + E^*_1V + \cdots + E^*_iV) \cap (E_0V + E_{1}V + \cdots + E_{d-i}V).
\end{equation}
Let $\{\tilde U_i \}_{i=0}^d$ denote the ${\tilde \Phi}$-split decomposition of ${\tilde V}$. Combining Definition \ref{def:dualphisetup} and (\ref{eq:defui2}), we find that for $0 \leq i \leq d$,
\begin{equation}
\label{eq:defuidual}
{\tilde U_i} = (E^{*\sigma}_d{\tilde V} + E^{*\sigma}_{d-1}{\tilde V} + \cdots + E^{*\sigma}_{d-i}{\tilde V}) \cap (E^\sigma_d{\tilde V} + E^\sigma_{d-1}{\tilde V} + \cdots + E^\sigma_{i}{\tilde V}).
\end{equation}
\noindent We now compare these two decompositions. 

\begin{lemma}
\label{lem:dualdecsp}
With reference to Definition \ref{def:dualphisetup}, the following {\rm (i)}, {\rm (ii)} are inverted dual.
\begin{enumerate}
 \item[\rm (i)]
The $\Phi$-split decomposition of $V$.
 \item[\rm (ii)]
The ${\tilde \Phi}$-split decomposition of ${\tilde V}$.
\end{enumerate}
\end{lemma}

\noindent {\it Proof:} We use the notation from above this lemma. For $0 \leq i, j \leq d$ with $i + j \neq d$, we show that $U_i$ and ${\tilde U_j}$ are orthogonal. We consider two cases: $i + j < d$ and $i + j > d$. First suppose that $i + j < d$. Abbreviate $M = E^*_0V + E^*_1V + \cdots + E^*_iV$ and $N = E^{*\sigma}_d{\tilde V} + E^{*\sigma}_{d-1}{\tilde V} + \cdots + E^{*\sigma}_{d-j}{\tilde V}$.
Observe that $U_i \subseteq M$ by (\ref{eq:defui2}) and ${\tilde U_j} \subseteq N$ by (\ref{eq:defuidual}). Moreover $M$ and $N$ are orthogonal by our assumption and Lemma \ref{lem:dualdecst2}. Therefore $U_i$ and ${\tilde U_j}$ are orthogonal. Next suppose that $i + j > d$. Abbreviate $S = E_0V + E_{1}V + \cdots + E_{d-i}V$ and $T = E^\sigma_d{\tilde V} + E^\sigma_{d-1}{\tilde V} + \cdots + E^\sigma_{j}{\tilde V}$. Observe $U_i \subseteq S$ by (\ref{eq:defui2}) and ${\tilde U_j} \subseteq T$ by (\ref{eq:defuidual}). 
Moreover $S$ and $T$ are orthogonal by our assumption and Lemma \ref{lem:dualdecst}. Therefore $U_i$ and ${\tilde U_j}$ are orthogonal and the result follows. 
\hfill $\Box$

\medskip

\noindent Let $0 \neq \xi_0 \in E_0V$ and recall from (\ref{eq:splitbasis}) that the sequence
\begin{eqnarray}
\label{eq:splitbasis2}
\eta^*_{d-i}(A^*) \xi_0 \qquad \qquad (0 \leq i \leq d) 
\end{eqnarray}
is a $\Phi$-split basis for $V$. Let $0 \neq {\tilde \xi_d} \in E^\sigma_d{\tilde V}$ and observe by Definition \ref{def:dualphisetup} that the sequence
\begin{eqnarray}
\label{eq:splitbasisdual}
\tau^*_{d-i}(A^{*\sigma}) {\tilde \xi_d} \qquad \qquad (0 \leq i \leq d) 
\end{eqnarray}
is a ${\tilde \Phi}$-split basis for ${\tilde V}$.  The bases (\ref{eq:splitbasis2}), (\ref{eq:splitbasisdual}) are related as follows. 

\begin{proposition}
\label{prop:splitvssplitdual}
With reference to Definition \ref{def:dualphisetup}, let $0 \neq \xi_0 \in E_0V$ and $0 \neq {\tilde \xi_d} \in E_d^\sigma {\tilde V}$.  Then for $0 \leq i,j \leq d$, 
\begin{eqnarray*}
\langle \eta^*_{i}(A^*) \xi_0, \tau^*_{j}(A^{*\sigma}) {\tilde \xi_d} \rangle = \delta_{i+j,d} \, \eta_d^*(\theta^*_0) \b{E_0^*\xi_0, {\tilde \xi_d}}.
\end{eqnarray*}
\end{proposition}

\noindent {\it Proof:} First suppose $i + j \neq d$. Then the result holds by Lemma \ref{lem:dualdecsp} and the comments above this proposition. Next suppose that $i + j = d$. Using (\ref{eq:primidtau}), Proposition \ref{prop:orthog}, and the definition of $\sigma$ from above Definition \ref{def:dualphisetup}, we find
\begin{eqnarray*}
\langle \eta^*_{i}(A^*) \xi_0, \tau^*_{j}(A^{*\sigma}) {\tilde \xi_d} \rangle 
&=&
\langle \eta_{d-j}^*(A^*) \xi_0, \tau_{j}^*(A^{*\sigma}) {\tilde \xi_d} \rangle \\
&=&
\langle \tau_{j}^*(A^{*}) \eta_{d-j}^*(A^*) \xi_0,  {\tilde \xi_d}   \rangle \\
&=&
\tau_{j}^*(\theta_{j}^*) \eta_{d-j}^*(\theta_{j}^*)  \langle E_{j}^* \xi_0, {\tilde \xi_d}  \rangle \\
&=&
\tau_{j}^*(\theta_{j}^*) \eta_{d-j}^*(\theta_{j}^*)  \langle E_{j}^* E_{j}^* \xi_0,  {\tilde \xi_d}  \rangle \\
&=&
\tau_{j}^*(\theta_{j}^*) \eta_{d-j}^*(\theta_{j}^*)  \langle E_{j}^* \xi_0, E_{j}^{*\sigma} {\tilde \xi_d}  \rangle \\
&=&
\tau_{j}^*(\theta_{j}^*) \eta_{d-j}^*(\theta_{j}^*)  \ell_{j} \langle E_0^* \xi_0, {\tilde \xi_d}      \rangle  \\
&=&
\eta_d^*(\theta^*_0) \langle {E_0^*\xi_0, {\tilde \xi_d}}  \rangle.
\end{eqnarray*}
\hfill $\Box$ 

\begin{corollary}
\label{cor:splitvssplitdual}
With reference to Definition \ref{def:dualphisetup}, let $\lbrace v_i \rbrace_{i=0}^d$ (resp. $\lbrace w_i \rbrace_{i=0}^d$) denote a basis for $V$ (resp. $\tilde V$). Suppose that $\lbrace v_i \rbrace_{i=0}^d$ and $\lbrace w_i \rbrace_{i=0}^d$ are inverted dual. Then the following {\rm (i)}, {\rm (ii)} are equivalent.
\begin{enumerate}
 \item[\rm (i)]
$\lbrace v_i \rbrace_{i=0}^d$ is a $\Phi$-split basis for $V$.
 \item[\rm (ii)]
$\lbrace w_i \rbrace_{i=0}^d$ is a ${\tilde \Phi}$-split basis for ${\tilde V}$.
\end{enumerate}
\end{corollary}

\noindent {\it Proof:} Use Proposition \ref{prop:splitvssplitdual}. \hfill $\Box$

\medskip
At the end of of Section \ref{sec:transvand}, we give the relationship between a $\Phi$-standard (resp. $\Phi^*$-standard) basis for $V$ and a ${\tilde \Phi^*}$-standard (resp. ${\tilde \Phi}$-standard) basis for ${\tilde V}$. This relationship is of a different type than the ones in the present section. 

\section{The transition matrices for a TH system}
Let $\Phi$ denote a TH system. In this section we consider several transition matrices associated with $\Phi$. First we clarify our terms. Let $\lbrace u_i \rbrace_{i=0}^d$ and $\lbrace v_i \rbrace_{i=0}^d$ denote bases for $V$. By the {\it transition matrix} from $\lbrace u_i \rbrace_{i=0}^d$ to $\lbrace v_i \rbrace_{i=0}^d$, we mean the matrix $T \in \hbox{Mat}_{d+1}(\K)$ such that $v_j = \sum_{i=0}^d T_{ij}u_i$ for $0 \leq j \leq d$.

\begin{definition}
{\rm \cite[Definition 10.6]{thinhess}}
\label{def:transmat}
\rm
Let $\Phi=(A;\lbrace E_i\rbrace_{i=0}^d;A^*;\lbrace E^*_i\rbrace_{i=0}^d)$ denote a TH system on $V$. Let $0 \neq \xi_0 \in E_0V$ and $0 \neq \xi^*_0 \in E^*_0V$.  Recall the $\Phi$-standard basis $\lbrace E^*_i \xi_0 \rbrace_{i=0}^d$ for $V$ and the $\Phi^*$-standard basis $\lbrace E_i \xi^*_0 \rbrace_{i=0}^d$ for $V$. Let $P \in \Mdf$ denote the transition matrix from $\lbrace E_i \xi^*_0 \rbrace_{i=0}^d$ to $\lbrace E^*_i \xi_0 \rbrace_{i=0}^d$, with $\xi_0, \xi_0^*$ chosen so that $\xi^*_0=E^*_0\xi_0$.
\end{definition}

\begin{theorem}
{\rm \cite[Theorem 10.8]{thinhess}}
\label{thm:pentries}
Let $\Phi$ denote a TH system with parameter array \\ $(\lbrace \theta_i \rbrace_{i=0}^d, \lbrace \theta^*_i \rbrace_{i=0}^d, \lbrace \phi_i \rbrace_{i=1}^d)$ and let $P$ denote the matrix from Definition \ref{def:transmat}. For $0 \leq i,j \leq d$, the $(i, j)$-entry of $P$ is equal to $\ell_j$ times
\begin{eqnarray}
\label{eq:psum}
\sum_{h=0}^d
\frac{(\theta_i-\theta_d)(\theta_i-\theta_{d-1})\cdots (\theta_i-\theta_{d-h+1})
(\theta^*_j-\theta^*_0)(\theta^*_j-\theta^*_1)
\cdots (\theta^*_j-\theta^*_{h-1})
}{\phi_1 \phi_2 \cdots \phi_h},
\end{eqnarray}
where $\ell_j$ is from Definition \ref{def:kappas}.
\end{theorem}

\begin{corollary}
\label{cor:pentries}
With reference to Definition \ref{def:transmat}, for $0 \leq i \leq d$ both
$$P_{i0} = 1, \qquad \qquad  P_{di} = \ell_i,$$ where $\ell_i$ is from Definition \ref{def:kappas}.
\end{corollary}

\noindent {\it Proof:} Use Theorem \ref{thm:pentries}. \hfill $\Box$

\medskip
\noindent The following definition is motivated by Definition \ref{def:transmat} and Corollary \ref{cor:pentries}. 

\begin{definition}
\rm
We call the matrix $P$ from Definition \ref{def:transmat} the {\it west normalized transition matrix} of $\Phi$. 
\end{definition}

\noindent Motivated by the sum (\ref{eq:psum}), we make a definition. Let $\lambda, \mu$ denote commuting indeterminates. Let $\K[\lambda, \mu]$ denote the $\K$-algebra consisting of the polynomials in $\lambda$, $\mu$ that have all coefficients in $\K$.

\begin{definition}
\label{def:scriptp}
\rm
Let $\Phi$ denote a TH system with parameter array $(\lbrace \theta_i \rbrace_{i=0}^d, \lbrace \theta^*_i \rbrace_{i=0}^d, \lbrace \phi_i \rbrace_{i=1}^d)$. Define $p \in \K[\lambda, \mu]$ by 
\begin{eqnarray}
\label{eq:p}
p = \sum_{h=0}^d \frac{\eta_h(\lambda)\tau^*_h(\mu)}{\phi_1\phi_2 \cdots \phi_h},
\end{eqnarray}
where $\lbrace \tau^*_i \rbrace_{i=0}^d$ and $\lbrace \eta_i \rbrace_{i=0}^d$ are from Notation \ref{def:tau}. We call $p$ the {\it two-variable polynomial} of $\Phi$. 
\end{definition}

\begin{example}
\rm
With reference to Definition \ref{def:scriptp}, assume $d = 2$. Then 
$$p = 1 +  \frac{(\lambda - \theta_2)(\mu - \theta_0^*)}{\phi_1} +  \frac{(\lambda - \theta_2)(\lambda - \theta_1)(\mu - \theta_0^*)(\mu - \theta_1^*)}{\phi_1 \phi_2}.$$
\end{example}

\begin{remark}
\label{rem:1}
\rm
With reference to Definition \ref{def:scriptp}, for $0 \leq i,j \leq d$ the scalar $p(\theta_i, \theta_j^*)$ is the sum (\ref{eq:psum}). 
\end{remark}

\begin{definition}
\label{def:transmatnorm}
\rm
Let $\Phi$ denote a TH system with parameter array $(\lbrace \theta_i \rbrace_{i=0}^d, \lbrace \theta^*_i \rbrace_{i=0}^d, \lbrace \phi_i \rbrace_{i=1}^d)$. Let $\mathcal P \in \Mdf$ denote the matrix with $(i,j)$-entry $p(\theta_i, \theta_j^*)$ for $0 \leq i, j \leq d$. Observe that by Theorem \ref{thm:pentries} and Remark \ref{rem:1}, the matrix $P$ from Definition \ref{def:transmat} is equal to $\mathcal P L$, where $L$ is the matrix from Definition \ref{def:diagli}.
\end{definition}

\begin{example}
\rm
With reference to Definition \ref{def:transmatnorm}, assume $d = 2$. Then 
\begin{eqnarray*}
\mathcal P = \left(
\begin{array}{c c c}
1 & 1 + \frac{(\theta_0 - \theta_2)(\theta^*_1-\theta^*_0)}{\phi_1} & 1 + \frac{(\theta_0 - \theta_2)(\theta^*_2-\theta^*_0)}{\phi_1} +   \frac{(\theta_0 - \theta_2)(\theta_0 - \theta_1)(\theta^*_2-\theta^*_0)(\theta^*_2-\theta^*_1)}{\phi_1\phi_2}  \\
1 & 1 + \frac{(\theta_1 - \theta_2)(\theta^*_1-\theta^*_0)}{\phi_1}  &  1 + \frac{(\theta_1 - \theta_2)(\theta^*_2-\theta^*_0)}{\phi_1} \\
1 & 1 & 1 
\end{array}
\right).
\end{eqnarray*}
\end{example}

\begin{corollary}
\label{cor:scriptpendpoint}
With reference to Definition \ref{def:transmatnorm}, for $0 \leq i \leq d$ both
$$\mathcal P_{i0} = 1, \qquad \qquad \ \mathcal P_{di} = 1.$$
\end{corollary}

\noindent {\it Proof:} Routine. \hfill $\Box$

\medskip

\noindent We now interpret the matrix $\mathcal P$ as a transition matrix. Let $\lbrace u_i \rbrace_{i=0}^d$ denote the inverted dual of a ${\tilde \Phi}$-standard basis for ${\tilde V}$. By Corollary \ref{cor:dualvstilde}, $\lbrace \ell_i u_i \rbrace_{i=0}^d$ is a $\Phi$-standard basis for $V$, where $\ell_i$ is from Definition \ref{def:kappas}. Recall the canonical bilinear form 
$\b{\;,\,}: V \times \tilde{V} \to \mathbb{K}$ from above Definition \ref{def:dualphisetup}. 

\begin{corollary}
\label{cor:transmatnorm2}
Let $\Phi=(A;\lbrace E_i\rbrace_{i=0}^d;A^*;\lbrace E^*_i\rbrace_{i=0}^d)$ denote a TH system on $V$. Let $0 \neq \xi^*_0 \in E^*_0V$ and $0 \neq {\tilde \xi_d} \in E_d^{\sigma}{\tilde V}$. Note by Lemma \ref{lem:dualdecst2} that $\langle \xi^*_0, {\tilde \xi_d} \rangle \neq 0$. Recall the $\Phi^*$-standard basis $\lbrace E_i \xi^*_0 \rbrace_{i=0}^d$ for $V$ and the ${\tilde \Phi}$-standard basis $\lbrace E^{*\sigma}_{d-i} {\tilde \xi_d} \rbrace_{i=0}^d$ for ${\tilde V}$. Let $\mathcal P$ denote the matrix from Definition \ref{def:transmatnorm}. Then $\alpha \mathcal P$ is the transition matrix from $\lbrace E_i \xi^*_0 \rbrace_{i=0}^d$ to the inverted dual of $\lbrace E^{*\sigma}_{d-i} {\tilde \xi_d} \rbrace_{i=0}^d$, where $\alpha$ is the reciprocal of $\langle \xi^*_0, {\tilde \xi_d} \rangle$. In particular if we choose $\xi_0^*, {\tilde \xi_d}$ so that $\langle \xi^*_0, {\tilde \xi_d} \rangle = 1$, then $\mathcal P$ is the transition matrix from $\lbrace E_i \xi^*_0 \rbrace_{i=0}^d$ to the inverted dual of $\lbrace E^{*\sigma}_{d-i} {\tilde \xi_d} \rbrace_{i=0}^d$.
\end{corollary}

\noindent {\it Proof:} Choose $\xi_0 \in E_0V$ so that $\xi^*_0 = E^*_0 \xi_0$, and recall the $\Phi$-standard basis $\lbrace E_i^* \xi_0 \rbrace_{i=0}^d$. The matrix $P$ from Definition \ref{def:transmat} is the transition matrix from $\lbrace E_i \xi^*_0 \rbrace_{i=0}^d$ to $\lbrace E_i^* \xi_0 \rbrace_{i=0}^d$. Now by Definition \ref{def:transmatnorm}, $\alpha \mathcal P$ is the transition matrix from $\lbrace E_i \xi^*_0 \rbrace_{i=0}^d$ to $\lbrace \alpha \ell_i^{-1} E_i^* \xi_0 \rbrace_{i=0}^d$. Moreover by Proposition \ref{prop:orthog}, $\{ \alpha \ell_i^{-1} E_i^* \xi_0 \}_{i=0}^d$ is the inverted dual of $\lbrace E^{*\sigma}_{d-i} {\tilde \xi_d} \rbrace_{i=0}^d$. Therefore $\alpha \mathcal P$ is the transition matrix from $\lbrace E_i \xi^*_0 \rbrace_{i=0}^d$ to the inverted dual of $\lbrace E^{*\sigma}_{d-i} {\tilde \xi_d} \rbrace_{i=0}^d$. \hfill $\Box$

\medskip

\noindent The following definition is motivated by Corollary \ref{cor:scriptpendpoint} and Corollary \ref{cor:transmatnorm2}.

\begin{definition}
\rm
We call the matrix $\mathcal P$ from Definition \ref{def:transmatnorm} the {\it west-south normalized transition matrix} of $\Phi$.
\end{definition}

\section{The transition matrices $P, \mathcal P$ and their relatives}
Let $\Phi$ denote a TH system. In the previous section we discussed two closely related transition matrices $P, \mathcal P$ associated with $\Phi$. In this section we find the relationship between $P, \mathcal P$ and their relatives. There are two types of relations; one type is best expressed in terms of $P$ and its relatives, while the other is best expressed in terms of $\mathcal P$ and its relatives. 

\begin{proposition}
\label{prop:ppstar}
{\rm \cite[Proposition 10.9]{thinhess}}
With reference to Definition \ref{def:transmat}, both
$$PP^* = P^*P = \nu I, \qquad \qquad  {\tilde P}{\tilde P^*} = {\tilde P^*}{\tilde P} = {\tilde \nu} I,$$
where $\nu, {\tilde \nu}$ are from Definition \ref{def:nu}.
\end{proposition}

Let $\{ u_i \}_{i=0}^d$ and $\{ v_i \}_{i=0}^d$ denote bases for $V$. Let $T \in \Mdf$ denote the transition matrix from $\{ u_i \}_{i=0}^d$ to $\{ v_i \}_{i=0}^d$. By elementary linear algebra, $T^t$ is the transition matrix from the dual of $\{ v_i \}_{i=0}^d$ to the dual of $\{ u_i \}_{i=0}^d$. Therefore $T^\varsigma$ is the transition matrix from the inverted dual of $\{ v_i \}_{i=0}^d$ to the inverted dual of $\{ u_i \}_{i=0}^d$.

\begin{proposition}
\label{prop:pgamma}
With reference to Definition \ref{def:transmatnorm}, both
\begin{eqnarray}
\label{eq:lab56}
\mathcal P^{\varsigma} = {\tilde \mathcal P^*}, \qquad \qquad (\mathcal P^*)^{\varsigma} = {\tilde \mathcal P}.
\end{eqnarray}

\end{proposition}

\noindent {\it Proof:} Choose $0 \neq \xi^*_0 \in E^*_0V$ and $0 \neq {\tilde \xi_d} \in E_d^{\sigma}{\tilde V}$ so that $\langle \xi^*_0, {\tilde \xi_d} \rangle = 1$. Recall the  $\Phi^*$-standard basis $\lbrace E_i \xi^*_0 \rbrace_{i=0}^d$ for $V$ and the ${\tilde \Phi}$-standard basis $\lbrace E^{*\sigma}_{d-i} {\tilde \xi_d} \rbrace_{i=0}^d$ for ${\tilde V}$.  By Corollary \ref{cor:transmatnorm2}, $\mathcal P$ is the transition matrix from $\lbrace E_i \xi^*_0 \rbrace_{i=0}^d$ to the inverted dual of $\lbrace E^{*\sigma}_{d-i} {\tilde \xi_d} \rbrace_{i=0}^d$. Moreover by Corollary \ref{cor:transmatnorm2} applied to ${\tilde \Phi^*}$, ${\tilde \mathcal P^*}$ is the transition matrix from $\lbrace E^{*\sigma}_{d-i} {\tilde \xi_d} \rbrace_{i=0}^d$ to the inverted dual of $\lbrace E_i \xi^*_0 \rbrace_{i=0}^d$. By these comments and the ones above this proposition we obtain the equation on the left in (\ref{eq:lab56}). The equation on the right
in (\ref{eq:lab56}) is similarly obtained. \hfill $\Box$

\medskip

\noindent The following lemma could be used to give another proof of Proposition \ref{prop:pgamma}.
 
\begin{lemma}
\label{lem:twovarpolydualrel}
With reference to Definition \ref{def:scriptp}, both $$p(\mu, \lambda) = {\tilde p^*}(\lambda, \mu), \qquad \qquad p^*(\mu, \lambda) = {\tilde p}(\lambda, \mu).$$
\end{lemma}

\noindent {\it Proof:} Combining Corollary \ref{cor:parelatives} and Definition \ref{def:scriptp}, we find 
$${\tilde p^*}(\lambda, \mu) = \sum_{h=0}^d \frac{\tau_h^*(\lambda)\eta_h(\mu)}{\phi_1 \phi_{2} \cdots \phi_{h}}. $$
Therefore $p(\mu, \lambda) = {\tilde p^*}(\lambda, \mu)$. The proof for the other claim is similar. \hfill $\Box$

\medskip

\noindent We will continue our discussion of the transition matrices $P, \mathcal P$ in Section \ref{sec:transvand}. In Sections \ref{sec:vand}--\ref{sec:invvand} we recall some linear algebra that will be needed in the discussion.  

\section{Vandermonde matrices and systems}
\label{sec:vand}
Let $\Phi$ denote a TH system. In Section \ref{sec:transvand} we will show that each of the transition matrices $P, \mathcal P$ of $\Phi$ has a certain structure said to be double Vandermonde. To prepare for that, over the next few sections we discuss some linear algebra related to Vandermonde matrices.

\begin{definition}
\label{def:gradedpoly}
\rm 
Let $n$ denote a nonnegative integer. Let $\lbrace f_i \rbrace_{i=0}^{n}$ denote a sequence of polynomials in $\K[\lambda]$. We say that $\lbrace f_i \rbrace_{i=0}^{n}$ is \emph{graded} whenever 
\begin{enumerate}
\item[\rm (i)]
$f_0 = 1$;
\item[\rm (ii)]
the degree of $f_i$ is equal to $i$ for $0 \leq i \leq n$. 
\end{enumerate} 
\end{definition}

\begin{definition}
\label{def:GVandermonde}
\rm
A matrix $X \in \Mdf$ is called \emph{west} \emph{Vandermonde} whenever the following {\rm (i)}, {\rm (ii)} hold.
\begin{enumerate}
 \item[\rm (i)]
There exists a sequence of mutually distinct scalars $\lbrace \theta_i \rbrace_{i=0}^d$ taken from $\K$ and a graded sequence of polynomials $\lbrace f_i \rbrace_{i=0}^d$ in $\K[\lambda]$ such that
\begin{eqnarray} 
\label{eq:vand}
X_{ij}=X_{i0}f_j(\theta_i) \qquad  (0 \leq i,j\leq d). 
\end{eqnarray}
 \item[\rm (ii)]
$X_{i0} \neq 0$ for $0 \leq i \leq d$. 
\end{enumerate}
\end{definition}

With reference to Definition \ref{def:GVandermonde}, assume $X$ is west Vandermonde. As we will see, the polynomials $\lbrace f_i \rbrace_{i=0}^d$ are uniquely determined by the sequence of scalars $\lbrace \theta_i \rbrace_{i=0}^d$ but the sequence $\lbrace \theta_i \rbrace_{i=0}^d$ is not unique. To facilitate our discussion of this issue, we introduce the following term. 

\begin{definition}
\label{def:parseq}
\rm
Let $X \in \Mdf$ denote a west Vandermonde matrix. Let $\lbrace \theta_i \rbrace_{i=0}^d$ denote a sequence of scalars taken from $\K$. We say that $X$ and $\lbrace \theta_i \rbrace_{i=0}^d$ are {\it compatible} whenever 
\begin{enumerate}
 \item[\rm (i)] 
$\theta_i \neq \theta_j$ if $i \neq j \qquad  (0 \leq i,j\leq d)$;
 \item[\rm (ii)]
there exists a graded sequence of polynomials $\lbrace f_i \rbrace_{i=0}^d$ in $\K[\lambda]$ that satisfies (\ref{eq:vand}).
\end{enumerate}
Observe that if $X$ and $\lbrace \theta_i \rbrace_{i=0}^d$ are compatible, then $X$ and $\lbrace \alpha \theta_i + \beta \rbrace_{i=0}^d$ are compatible for any $\alpha, \beta \in \K$ with $\alpha \neq 0$. 
\end{definition}

\begin{lemma}
\label{lem:parseq}
Let $X \in \Mdf$ denote a west Vandermonde matrix. Let $\lbrace \theta_i \rbrace_{i=0}^d$ denote a sequence of scalars taken from $\K$. Then the following {\rm (i)}, {\rm (ii)} are equivalent provided $d \geq 1$.
\begin{enumerate}
 \item[\rm (i)]
$X$ and $\lbrace \theta_i \rbrace_{i=0}^d$ are compatible.
 \item[\rm (ii)]
There exists $a, b \in \K$ with $a \neq 0$ such that $\theta_i = a X_{i1}/X_{i0} + b$ for $0 \leq i \leq d$.
\end{enumerate}
\end{lemma}

\noindent {\it Proof:} (i) $\Rightarrow$ (ii) By Definition \ref{def:parseq}, there exists a polynomial $f_1 \in \K[\lambda]$ of degree 1 such that $X_{i1} = X_{i0}f_1(\theta_i)$ for $0 \leq i \leq d$. Write $f_1 = \alpha \lambda + \beta$ for some $\alpha, \beta \in \K$ with $\alpha \neq 0$. Thus  
$X_{i1} = X_{i0} (\alpha \theta_i + \beta)$ for $0 \leq i \leq d$. Rearranging terms, we find that there exists $a, b \in \K$ with $a \neq 0$ such that $\theta_i = a X_{i1}/X_{i0} + b$ for $0 \leq i \leq d$. \\
(ii) $\Rightarrow$ (i) By Definition \ref{def:GVandermonde}, there exists a sequence of scalars $\lbrace \theta_i' \rbrace_{i=0}^d$ taken from $\K$ that is compatible with $X$. By the previous part, there exists $a', b' \in \K$ with $a' \neq 0$ such that $\theta_i' = a' X_{i1}/X_{i0} + b'$ for $0 \leq i \leq d$. Thus there exists $\alpha, \beta \in \K$ with $\alpha \neq 0$ such that $\theta_i = \alpha \theta_i' + \beta$ for $0 \leq i \leq d$. Now $X$ and $\lbrace \theta_i \rbrace_{i=0}^d$ are compatible by the observation at the end of Definition \ref{def:parseq}. \hfill $\Box$

\begin{lemma}
\label{lem:thetarelated}
Let $X \in \Mdf$ denote a west Vandermonde matrix. Let $\lbrace \theta_i \rbrace_{i=0}^d$ denote a sequence of scalars taken from $\K$ that is compatible with $X$. 
Let $\lbrace \theta_i' \rbrace_{i=0}^d$ denote a sequence of scalars taken from $\K$. Then the following {\rm (i)}, {\rm (ii)} are equivalent. 
\begin{enumerate}
 \item[\rm (i)]
$X$ and $\lbrace \theta_i' \rbrace_{i=0}^d$ are compatible.
 \item[\rm (ii)]
There exists $\alpha, \beta \in \K$ with $\alpha \neq 0$ such that $\theta_i' = \alpha \theta_i + \beta$ for $0 \leq i \leq d$.
\end{enumerate}
\end{lemma}

\noindent {\it Proof:} (i) $\Rightarrow$ (ii) Routine by Lemma \ref{lem:parseq}. \\
(ii) $\Rightarrow$ (i) This is the observation at the end of Definition \ref{def:parseq}. \hfill $\Box$

\begin{lemma}
\label{lem:polyvand}
Let $\lbrace \theta_i \rbrace_{i=0}^d$ denote a sequence of mutually distinct scalars taken from $\fld$. Let $X \in \Mdf$ denote a west Vandermonde matrix that is compatible with $\lbrace \theta_i \rbrace_{i=0}^d$. Then there exists a unique graded sequence of polynomials $\lbrace f_i \rbrace_{i=0}^d$ in $\K[\lambda]$ that satisfies (\ref{eq:vand}). 
\end{lemma}

\noindent {\it Proof:} By Definition \ref{def:parseq} there exists a graded sequence of polynomials $\lbrace f_i \rbrace_{i=0}^d$ in $\K[\lambda]$ that satisfies (\ref{eq:vand}). We show that this sequence is unique. Suppose that $\lbrace f_i' \rbrace_{i=0}^d$ is a graded sequences of polynomials in $\K[\lambda]$ that satisfies (\ref{eq:vand}). We show that $f_i' = f_i$ for $0 \leq i \leq d$. Let $i$ be given and define $g_i = f_i' - f_i$. Using (\ref{eq:vand}), we find $g_i(\theta_j) = 0$ for $0 \leq j \leq d$. Since $\lbrace \theta_i \rbrace_{i=0}^d$ are mutually distinct and $g_i$ has degree at most $i$, it follows that $g_i = 0$. Therefore $f_i' = f_i$. We have shown that the sequence $\lbrace f_i \rbrace_{i=0}^d$ is unique. \hfill $\Box$ 

\begin{lemma}
\label{lem:vandfrompolysca}
Let $\lbrace \theta_i \rbrace_{i=0}^d$ denote a sequence of mutually distinct scalars taken from $\K$. Let $\lbrace f_i \rbrace_{i=0}^{d}$ denote a graded sequence of polynomials in $\K[\lambda]$. Let $\lbrace c_i \rbrace_{i=0}^d$ denote a sequence of nonzero scalars taken from $\K$. Define $X \in \Mdf$ such that $X_{ij}=c_if_j(\theta_i)$ for $0 \leq i,j\leq d$. Then $X$ is west Vandermonde and compatible with $\lbrace \theta_i \rbrace_{i=0}^d$. Moreover $\lbrace f_i \rbrace_{i=0}^{d}$ are the corresponding polynomials from Lemma \ref{lem:polyvand}.
\end{lemma}

\noindent {\it Proof:} Routine. \hfill $\Box$

\begin{lemma}
\label{lem:westvandunique}
Let $\lbrace \theta_i \rbrace_{i=0}^d$ denote a sequence of mutually distinct scalars taken from $\K$.  Let $X \in \Mdf$ denote a west Vandermonde matrix that is compatible with $\lbrace \theta_i \rbrace_{i=0}^d$ and let $\lbrace f_i \rbrace_{i=0}^{d}$ denote the corresponding polynomials from Lemma \ref{lem:polyvand}. Let $X' \in \Mdf$. Then the following {\rm (i)}, {\rm (ii)} are equivalent.
\begin{enumerate}
 \item[\rm (i)]
$X'$ is a west Vandermonde matrix that is compatible with $\lbrace \theta_i \rbrace_{i=0}^d$ and $\lbrace f_i \rbrace_{i=0}^{d}$ are the corresponding polynomials from Lemma \ref{lem:polyvand}.
 \item[\rm (ii)]
There exists an invertible diagonal matrix $D \in \Mdf$ such that $X' = D X$. 
\end{enumerate}
\end{lemma}

\noindent {\it Proof:} (i) $\Rightarrow$ (ii) Using (\ref{eq:vand}), we find that for $0 \leq i,j \leq d$,
\begin{eqnarray}
\label{eq:one}
\frac{X_{ij}'}{X_{i0}'} = \frac{X_{ij}}{X_{i0}} = f_j(\theta_i).
\end{eqnarray}
Define a diagonal matrix $D \in \Mdf$ with $(i,i)$-entry $X_{i0}'/X_{i0}$ for $0 \leq i \leq d$. Observe that $D$ is invertible. Moreover $X' = D X$ by (\ref{eq:one}). \\
(ii) $\Rightarrow$ (i) Since $X' = D X$ we find that for $0 \leq i,j\leq d$,
\begin{eqnarray*}
X_{ij}' = X_{i0}' \frac{X_{ij}}{X_{i0}}. 
\end{eqnarray*}
By this and (\ref{eq:vand}), we find that $X_{ij}'=X_{i0}'f_j(\theta_i)$ and (i) follows. \hfill $\Box$

\begin{lemma}
\label{lem:westvandright}
Let $\lbrace \theta_i \rbrace_{i=0}^d$ denote a sequence of mutually distinct scalars taken from $\K$.  Let $X \in \Mdf$ denote a west Vandermonde matrix that is compatible with $\lbrace \theta_i \rbrace_{i=0}^d$ and let $\lbrace f_i \rbrace_{i=0}^{d}$ denote the corresponding polynomials from Lemma \ref{lem:polyvand}. Let $D \in \Mdf$ denote an invertible diagonal matrix. Then $X D$ is a west Vandermonde matrix that is compatible with $\lbrace \theta_i \rbrace_{i=0}^d$ and $\lbrace D_{ii}f_i/D_{00}  \rbrace_{i=0}^{d}$ are the corresponding polynomials from Lemma \ref{lem:polyvand}.
\end{lemma}

\noindent {\it Proof:} Routine using (\ref{eq:vand}). \hfill $\Box$

\begin{definition}
\rm
By a {\it west} {\it Vandermonde system} in $\Mdf$ we mean a sequence $(X, \lbrace \theta_i \rbrace_{i=0}^d)$ such that 
\begin{enumerate}
 \item[\rm (i)] 
$X$ is a west Vandermonde matrix in $\Mdf$;
 \item[\rm (ii)]
$\lbrace \theta_i \rbrace_{i=0}^d$ is a sequence of mutually distinct scalars taken from $\fld$ that is compatible with $X$. 
\end{enumerate}
\end{definition}

\begin{definition}
\label{def:vandpoly10}
\rm
Let $(X, \lbrace \theta_i \rbrace_{i=0}^d)$ denote a west Vandermonde system in $\Mdf$. In Lemma \ref{lem:polyvand} we associated $(X, \lbrace \theta_i \rbrace_{i=0}^d)$ with some polynomials $\lbrace f_i \rbrace_{i=0}^{d}$. For convenience, let $f_{d+1} = \prod_{i=0}^d (\lambda - \theta_i)$. We call $\lbrace f_i \rbrace_{i=0}^{d+1}$ the {\it polynomials of $(X, \lbrace \theta_i \rbrace_{i=0}^d)$}. 
\end{definition}

\begin{definition}
\rm
Let $X \in \Mdf$. Let $X' \in \Mdf$ denote a matrix that is obtained by rotating $X$ clockwise $90$ degrees. We call $X$ {\it south Vandermonde} whenever $X'$ is west Vandermonde. 
\end{definition}

\noindent The above notions regarding west Vandermonde matrices carry over to south Vandermonde matrices. 

\medskip

\noindent We end this section with a comment. 

\begin{lemma}
\label{lem:vandinvertible}
Let $X \in \Mdf$ denote a west or south Vandermonde matrix. Then $X$ is invertible. 
\end{lemma}

\noindent {\it Proof:} First assume that $X$ is west Vandermonde. Perform invertible row and column operations on $X$ so that the resulting matrix $X'$ has $(i,j)$-entry $\theta_i^j$ for $0 \leq i, j \leq d$. The determinant of $X'$ is equal to $\prod_{0 \leq i < j \leq d} (\theta_j - \theta_i)$.  The $\lbrace \theta_i \rbrace_{i=0}^d$ are mutually distinct so this determinant is nonzero. Therefore $X'$ is invertible so $X$ is invertible. The case of south Vandermonde is similar. \hfill $\Box$

\section{Hessenberg matrices and graded sequences of polynomials}
Recall the notion of a Hessenberg matrix from Section \ref{sec:intro}. In the next section we discuss the role Vandermonde matrices play in the diagonalization of Hessenberg matrices. To prepare for that, in this section we discuss the relationship between Hessenberg matrices and graded sequences of polynomials. 

\begin{lemma}
\label{lem:minpolyofhess}
Let $H \in \Mdf$ denote a Hessenberg matrix. Then the minimal polynomial of $H$ is equal to the characteristic polynomial of $H$. 
\end{lemma}

\noindent {\it Proof:} Using the Hessenberg shape of $H$, we find  $I, H, H^2, \ldots, H^d$ are linearly independent. Therefore the minimal polynomial of $H$ has degree $d+1$. The result follows. \hfill $\Box$

\medskip

Given a Hessenberg matrix $H$, we are interested in finding the polynomial in Lemma \ref{lem:minpolyofhess}.

\begin{notation}
\label{not:cofh}
\rm
Let $H \in \Mdf$ denote a Hessenberg matrix. We denote by $c_H$ the product $\prod_{i=1}^{d} H_{i,i-1}$. Observe that $c_H$ is nonzero.
\end{notation}

\begin{definition}
\label{def:polyofhess}
\rm
Let $H \in \Mdf$ denote a Hessenberg matrix. Define a sequence of polynomials $\lbrace f_i \rbrace_{i=0}^{d+1}$ in $\K[\lambda]$ such that 
\begin{enumerate}
 \item[\rm (i)] 
$f_0 = 1$;
 \item[\rm (ii)]
$\lambda f_j = \sum_{i=0}^{j+1} H_{ij}f_i$ for $0 \leq j \leq d-1$;
 \item[\rm (iii)]
$\lambda f_d = c_H^{-1} f_{d+1} + \sum_{i=0}^d H_{id}f_i$, where $c_H$ is from Notation \ref{not:cofh}.  
\end{enumerate}
We call $\lbrace f_i \rbrace_{i=0}^{d+1}$ the {\it polynomials of $H$}.  
\end{definition}

\begin{definition}
\label{def:standardpoly}
\rm
A graded sequence of polynomials $\lbrace f_i \rbrace_{i=0}^{d+1}$ in $\K[\lambda]$ is called {\it standard} whenever $f_{d+1}$ is monic. 
\end{definition}

\begin{lemma}
\label{lem:polyofH}
Let $H \in \Mdf$ denote a Hessenberg matrix with polynomials $\lbrace f_i \rbrace_{i=0}^{d+1}$. Then the following {\rm (i)--(iii)} hold. 
\begin{enumerate}
 \item[\rm (i)]
For $0 \leq i \leq d$, $f_i$ has degree $i$ with $\lambda^{i}$ coefficient $(\prod_{j=1}^{i} H_{j,j-1})^{-1}$.
 \item[\rm (ii)]
$f_{d+1}$ is monic with degree $d+1$. 
 \item[\rm (iii)]
The sequence $\lbrace f_i \rbrace_{i=0}^{d+1}$ is graded and standard. 
\end{enumerate}
\end{lemma}

\noindent {\it Proof:} Routine. \hfill $\Box$

\medskip

Let $I$ denote the identity matrix in $\Mdf$. For $0 \leq i \leq d$, let $\epsilon_i$ denote the $i^{th}$ column of $I$. Observe that $\lbrace \epsilon_i \rbrace_{i=0}^d$ is a basis for the vector space $\K^{d+1}$.

\begin{lemma}
\label{lem:lab452}
Let $H \in \Mdf$ denote a Hessenberg matrix. Then there exists a unique standard graded sequence of polynomials $\lbrace f_i \rbrace_{i=0}^{d+1}$ in $\K[\lambda]$ such that $f_i(H) \epsilon_0 = \epsilon_i$ for $0 \leq i \leq d$ and $f_{d+1}(H) \epsilon_0 = 0$. The $\lbrace f_i \rbrace_{i=0}^{d+1}$ are the polynomials of $H$ from Definition \ref{def:polyofhess}.
\end{lemma}

\noindent {\it Proof:} Concerning existence, let $\lbrace f_i \rbrace_{i=0}^{d+1}$ denote the polynomials of $H$. By Lemma \ref{lem:polyofH}(iii) the sequence $\lbrace f_i \rbrace_{i=0}^{d+1}$ is graded and standard. We show that $f_i(H) \epsilon_0 = \epsilon_i$ for $0 \leq i \leq d$ and $f_{d+1}(H) \epsilon_0 = 0$.  Abbreviate $v_i=f_i(H) \epsilon_0$ for $0 \leq i \leq d+1$ and note that $v_0 = \epsilon_{0}$. From Definition \ref{def:polyofhess}, we have
\begin{eqnarray}
Hv_j =  \sum_{i=0}^{j+1} H_{ij}v_i \qquad (0 \leq j \leq d), 
\label{eq:rel}
\end{eqnarray}
where $H_{d+1, d} = c_H^{-1}$. 
By the definition of $\lbrace \epsilon_i \rbrace_{i=0}^{d}$, we have
\begin{eqnarray}
H \epsilon_j = \sum_{i=0}^{j+1} H_{ij} \epsilon_i \qquad (0 \leq j \leq d),
\label{eq:rel2}
\end{eqnarray}
where $\epsilon_{d+1} = 0$. 
Comparing (\ref{eq:rel}), (\ref{eq:rel2}) and using $v_0=\epsilon_0$, we find $v_i=\epsilon_i$ for $0 \leq i \leq d+1$. Therefore $f_i(H) \epsilon_0 = \epsilon_i$ for $0 \leq i \leq d$ and $f_{d+1}(H) \epsilon_0 = 0$. Concerning uniqueness, let $\lbrace f_i' \rbrace_{i=0}^{d+1}$ denote a standard graded sequence of polynomials in $\K[\lambda]$ such that $f_i'(H) \epsilon_0 = \epsilon_i$ for $0 \leq i \leq d$ and $f_{d+1}'(H) \epsilon_0 = 0$. We show that $f_i' = f_i$ for $0 \leq i \leq d+1$. Let $i$ be given and define $g_i = f_i'-f_i$. Observe that $g_i(H) \epsilon_0 = 0$. Thus $g_i(H) \epsilon_j = g_i(H)f_j(H) \epsilon_0 = f_j(H)g_i(H) \epsilon_0 = 0$ for $0 \leq j \leq d$, so $g_i(H) = 0$. Therefore the minimal polynomial of $H$ divides $g_i$. The polynomial $g_i$ has degree at most $d$, and the minimal polynomial of $H$ has degree $d+1$ by Lemma \ref{lem:minpolyofhess}. Therefore $g_i = 0$ so $f_i' = f_i$. \hfill $\Box$

\begin{corollary}
\label{cor:fdplusoneminpoly}
Let $H \in \Mdf$ denote a Hessenberg matrix with polynomials $\lbrace f_i \rbrace_{i=0}^{d+1}$. Then $f_{d+1}$ is both the minimal polynomial and the characteristic polynomial of $H$. 
\end{corollary}

\noindent {\it Proof:} Using Lemma \ref{lem:lab452} we find that $f_{d+1}(H) \epsilon_i = f_{d+1}(H) f_i(H) \epsilon_0 = f_i(H)f_{d+1}(H) \epsilon_0 = 0$ for $0 \leq i \leq d$. Therefore $f_{d+1}(H) = 0$. The result follows by Lemma \ref{lem:minpolyofhess} and Lemma \ref{lem:polyofH}(ii). \hfill $\Box$

\medskip

So far, given a Hessenberg matrix we obtain a graded sequence of polynomials. Now turning things around, given a graded sequence of polynomials we obtain a Hessenberg matrix.

\begin{definition}
\label{def:connectioncoeff}
\rm
Let $\lbrace f_i \rbrace_{i=0}^{d+1}$ denote a graded sequence of polynomials in $\K[\lambda]$. Observe that for $0 \leq j \leq d$, $\lambda f_j$ is in the span of $\lbrace f_i \rbrace_{i=0}^{j+1}$. So for $0 \leq j \leq d$, there exists a unique sequence of scalars $\lbrace c_{ij} \rbrace_{i=0}^{j+1}$ taken from $\K$ such that $\lambda f_j = \sum_{i=0}^{j+1} c_{ij}f_i$. We call the scalar $c_{ij}$ the {\it $(i,j)$-connection coefficient} for the given graded sequence of polynomials. 
\end{definition}

\begin{definition}
\label{def:polytohess}
\rm
Let $\lbrace f_i \rbrace_{i=0}^{d+1}$ denote a graded sequence of polynomials in $\K[\lambda]$. By the {\it connection coefficient matrix} of $\lbrace f_i \rbrace_{i=0}^{d+1}$, we mean the Hessenberg matrix $H \in \Mdf$ such that $H_{ij} = c_{ij}$ for $0 \leq i,j \leq d$, $i-j \leq 1$. The scalars $c_{ij}$ are from Definition \ref{def:connectioncoeff}. Observe that the scalar $c_{d+1,d}$ plays no role in the definition of $H$.
\end{definition}

\begin{lemma}
\label{lem:polytohess}
Let $\lbrace f_i \rbrace_{i=0}^{d+1}$ (resp. $\lbrace f_i' \rbrace_{i=0}^{d+1}$) denote a graded sequence of polynomials in $\K[\lambda]$ with connection coefficient matrix $H$ (resp. $H'$). Then the following {\rm (i)}, {\rm (ii)} are equivalent. 
\begin{enumerate}
 \item[\rm (i)]
$H = H'$.
 \item[\rm (ii)]
$f_i = f_i'$ for $0 \leq i \leq d$ and there exists $0 \neq c \in \K$ such that $f_{d+1} = c f_{d+1}'$. 
\end{enumerate}
\end{lemma}

\noindent {\it Proof:} Routine by Definition \ref{def:polytohess}. \hfill $\Box$

\begin{lemma}
\label{lem:polyscatohess}
Let $\lbrace f_i \rbrace_{i=0}^{d+1}$ denote a standard graded sequence of polynomials in $\K[\lambda]$ and let $H \in \Mdf$ denote a Hessenberg matrix. Then the following {\rm (i), (ii)} are equivalent. 
\begin{enumerate}
 \item[\rm (i)] 
$\lbrace f_i \rbrace_{i=0}^{d+1}$ are the polynomials of $H$. 
 \item[\rm (ii)]
$H$ is the connection coefficient matrix of $\lbrace f_i \rbrace_{i=0}^{d+1}$.
\end{enumerate}
\end{lemma}

\noindent {\it Proof:} Routine. \hfill $\Box$

\section{A Vandermonde matrix as a transition matrix}
In this section we discuss the role that Vandermonde matrices play in the diagonalization of a Hessenberg matrix. 

\begin{lemma}
\label{lem:hessmultfree}
Let $H \in \Mdf$ denote a Hessenberg matrix with polynomials $\lbrace f_i \rbrace_{i=0}^{d+1}$. Then the following {\rm (i)--(iii)} are equivalent. 
\begin{enumerate}
 \item[\rm (i)]
$H$ is diagonalizable.
 \item[\rm (ii)]
$H$ is multiplicity-free.
 \item[\rm (iii)]
$f_{d+1}$ has $d+1$ distinct roots in $\K$. 
\end{enumerate}
\end{lemma}

\noindent{\it Proof:} Recall that the polynomial $f_{d+1}$ has degree $d+1$ and it is the minimal polynomial of $H$ by Corollary \ref{cor:fdplusoneminpoly}. By elementary linear algebra a matrix in $\Mdf$ is diagonalizable if and only if its minimal polynomial has distinct roots in $\K$. The result follows.  \hfill $\Box$

\smallskip

Let $H \in \Mdf$ denote a multiplicity-free Hessenberg matrix and let $\lbrace \theta_i \rbrace_{i=0}^d$ denote an ordering of the eigenvalues of $H$. Let $D \in \Mdf$ denote the diagonal matrix with $(i,i)$-entry $\theta_i$ for $0 \leq i \leq d$. By elementary linear algebra, there exists an invertible $X \in \Mdf$ such that $H = X^{-1} D X$. We comment on the uniqueness of $X$. Suppose that $Y \in \Mdf$ is invertible and $H = Y^{-1} D Y$. Then $X^{-1} D X = Y^{-1} D Y$ so $D Y X^{-1} = Y X^{-1} D$. Therefore $Y X^{-1}$ is diagonal by Lemma \ref{lem:multfreecomm}. By construction $Y X^{-1}$ is invertible. By these comments there exists an invertible diagonal matrix $\Delta \in \Mdf$ such that $Y = \Delta X$.

\begin{lemma}
\label{lem:appbw}
Let $H \in \Mdf$ denote a multiplicity-free Hessenberg matrix with polynomials $\lbrace f_i \rbrace_{i=0}^{d+1}$. Let $\lbrace \theta_i \rbrace_{i=0}^d$ denote an ordering of the eigenvalues of $H$ and let $D \in \Mdf$ denote the diagonal matrix with $(i,i)$-entry $\theta_i$ for $0 \leq i \leq d$.  For $X \in \Mdf$, the following {\rm (i)}, {\rm (ii)} are equivalent. 
\begin{enumerate}
 \item[\rm (i)]
$X$ is invertible and $H = X^{-1} D X$.
 \item[\rm (ii)]
$(X, \lbrace \theta_i \rbrace_{i=0}^d)$ is a west Vandermonde system with polynomials $\lbrace f_i \rbrace_{i=0}^{d+1}$.  
\end{enumerate}
\end{lemma}

\noindent{\it Proof:} (i) $\Rightarrow$ (ii) Observe that $\lbrace \theta_i \rbrace_{i=0}^d$ are mutually distinct since $H$ is multiplicity-free. We show that 
\begin{eqnarray} 
\label{eq:vandtemp}
X_{ij}=X_{i0}f_j(\theta_i) \qquad  (0 \leq i,j\leq d).
\end{eqnarray}
Since $H = X^{-1} D X$ we have $f_j(H) = X^{-1} f_j(D) X$ so $X f_j(H) =  f_j(D) X$. Hence $\epsilon_i^t X f_j(H) \epsilon_0 =  \epsilon_i^t f_j(D) X \epsilon_0$. Simplify this equation using $f_j(H) \epsilon_0 = \epsilon_j$ from Lemma \ref{lem:lab452} together with matrix multiplication to obtain (\ref{eq:vandtemp}). By (\ref{eq:vandtemp}) and since $X$ is invertible we find $X_{i0} \neq 0$ for $0 \leq i \leq d$. By Corollary \ref{cor:fdplusoneminpoly} we have $f_{d+1} = \prod_{i=0}^d (\lambda - \theta_i)$. By these comments $(X, \lbrace \theta_i \rbrace_{i=0}^d)$ is a west Vandermonde system with polynomials $\lbrace f_i \rbrace_{i=0}^{d+1}$. \\
(ii) $\Rightarrow$ (i) By Lemma \ref{lem:vandinvertible} $X$ is invertible. We now show that $H = X^{-1} D X$. For $0 \leq i \leq d$, evaluate the equations in Definition \ref{def:polyofhess}(ii),(iii) at $\lambda = \theta_i$. In the resulting equations multiply each side by $X_{i0}$ and simplify using (\ref{eq:vand}) and Corollary \ref{cor:fdplusoneminpoly} to obtain $\theta_i X_{ij} = \sum_{n=0}^{d} H_{nj} X_{in}$ for $0 \leq j \leq d$. Therefore $D X = X H$ so $H = X^{-1} D X$.  \hfill $\Box$

\begin{corollary}
\label{cor:appbwtwist}
Let $(X, \lbrace \theta_i \rbrace_{i=0}^d)$ denote a west Vandermonde system with polynomials $\lbrace f_i \rbrace_{i=0}^{d+1}$. Let $D \in \Mdf$ denote the diagonal matrix with $(i,i)$-entry $\theta_i$ for $0 \leq i \leq d$. Then $X^{-1} D X$ is the connection coefficient matrix of $\lbrace f_i \rbrace_{i=0}^{d+1}$. Moreover $X^{-1} D X$ is multiplicity-free and Hessenberg. 
\end{corollary}

\noindent {\it Proof:} Let $H \in \Mdf$ denote the connection coefficient matrix of $\lbrace f_i \rbrace_{i=0}^{d+1}$. Then $H$ is Hessenberg by Definition \ref{def:polytohess}. We show that $H = X^{-1} D X$ and that $H$ is multiplicity-free. By Lemma \ref{lem:polyscatohess} the $\lbrace f_i \rbrace_{i=0}^{d+1}$ are the polynomials of $H$. Since $f_{d+1} = \prod_{i=0}^d (\lambda - \theta_i)$ and $\lbrace \theta_i \rbrace_{i=0}^d$ are mutually distinct, we find using Corollary \ref{cor:fdplusoneminpoly} that $H$ is multiplicity-free. Now by Lemma \ref{lem:appbw} we find $H = X^{-1} D X$. The result follows. 
\hfill $\Box$

\medskip

\noindent We have been discussing west Vandermonde systems. We now obtain analogous results for south Vandermonde systems.

\begin{definition}
\label{def:dualpolynomial}
\rm
Let $\lbrace f_i \rbrace_{i=0}^{d+1}$ denote a standard graded sequence of polynomials in $\K[\lambda]$. Let $H \in \Mdf$ denote the connection coefficient matrix of $\lbrace f_i \rbrace_{i=0}^{d+1}$. Recall from Definition \ref{def:polytohess} that $H$ is Hessenberg, so $H^{\varsigma}$ is Hessenberg. The polynomials of $H^{\varsigma}$ will be denoted by $\lbrace f^{\varsigma}_i \rbrace_{i=0}^{d+1}$. The two polynomial sequences $\lbrace f_i \rbrace_{i=0}^{d+1}$ and $\lbrace f^{\varsigma}_i \rbrace_{i=0}^{d+1}$  are said to be {\it associated}. Note that $f^{\varsigma}_{d+1} = f_{d+1}$ by Corollary \ref{cor:fdplusoneminpoly}.
\end{definition}

\begin{lemma}
\label{lem:appbs}
Let $H \in \Mdf$ denote a multiplicity-free Hessenberg matrix with polynomials $\lbrace f_i \rbrace_{i=0}^{d+1}$. Let $\lbrace \theta_i \rbrace_{i=0}^d$ denote an ordering of the eigenvalues of $H$ and let $D \in \Mdf$ denote the diagonal matrix with $(i,i)$-entry $\theta_i$ for $0 \leq i \leq d$. For $X \in \Mdf$, the following {\rm (i)}, {\rm (ii)} are equivalent. 
\begin{enumerate}
 \item[\rm (i)] 
$X$ is invertible and $H = X D X^{-1}$.
 \item[\rm (ii)]
$(X, \lbrace \theta_i \rbrace_{i=0}^d)$ is a south Vandermonde system with polynomials $\lbrace f^{\varsigma}_i \rbrace_{i=0}^{d+1}$. 
\end{enumerate}
\end{lemma}

\noindent {\it Proof:} (i) $\Rightarrow$ (ii) In the equation $H = X D X^{-1}$, apply $\varsigma$ to each side to obtain $H^\varsigma = (X^\varsigma)^{-1} D^\varsigma X^\varsigma$. By this and Lemma \ref{lem:appbw} the sequence $(X^\varsigma, \lbrace \theta_{d-i} \rbrace_{i=0}^d)$ is a west Vandermonde system with polynomials $\lbrace f^{\varsigma}_i \rbrace_{i=0}^{d+1}$. Therefore $(X, \lbrace \theta_i \rbrace_{i=0}^d)$ is a south Vandermonde system with polynomials $\lbrace f^{\varsigma}_i \rbrace_{i=0}^{d+1}$. \\
(ii) $\Rightarrow$ (i) By assumption $(X, \lbrace \theta_i \rbrace_{i=0}^d)$ is a south Vandermonde system with polynomials $\lbrace f^{\varsigma}_i \rbrace_{i=0}^{d+1}$. 
Therefore $(X^\varsigma, \lbrace \theta_{d-i} \rbrace_{i=0}^d)$ is a west Vandermonde system with polynomials $\lbrace f^{\varsigma}_i \rbrace_{i=0}^{d+1}$. Applying Lemma \ref{lem:appbw} we find that $X^\varsigma$ is invertible and $H^\varsigma = (X^\varsigma)^{-1} D^\varsigma X^\varsigma$. Applying $\varsigma$ we find that $X$ is invertible and $H = X D X^{-1}$. 
\hfill $\Box$

\begin{corollary}
\label{cor:appbstwist}
Let $(X, \lbrace \theta_i \rbrace_{i=0}^d)$ denote a south Vandermonde system with polynomials $\lbrace f_i \rbrace_{i=0}^{d+1}$. Let $D \in \Mdf$ denote the diagonal matrix with $(i,i)$-entry $\theta_i$ for $0 \leq i \leq d$.  Then $X D X^{-1}$ is the connection coefficient matrix of $\lbrace f^{\varsigma}_i \rbrace_{i=0}^{d+1}$. Moreover $X D X^{-1}$ is multiplicity-free and Hessenberg. 
\end{corollary}

\noindent {\it Proof:} Similar to the proof of Corollary \ref{cor:appbwtwist}. \hfill $\Box$

\section{The inverse of a Vandermonde matrix}
\label{sec:invvand}
Let $X \in \Mdf$ denote a west or south Vandermonde matrix. In Lemma \ref{lem:vandinvertible} we showed that $X$ is invertible. In this section we discuss the matrix $X^{-1}$. 

\begin{proposition}
\label{prop:inversevandisvand}
Let $(X, \lbrace \theta_i \rbrace_{i=0}^d)$ denote a west Vandermonde system with polynomials $\lbrace f_i \rbrace_{i=0}^{d+1}$. Then the following {\rm (i)}, {\rm (ii)} hold. 
\begin{enumerate}
 \item[\rm (i)]
$(X^{-1}, \lbrace \theta_i \rbrace_{i=0}^d)$ is a south Vandermonde system with polynomials $\lbrace f^{\varsigma}_i \rbrace_{i=0}^{d+1}$, where $\lbrace f^{\varsigma}_i \rbrace_{i=0}^{d+1}$ are the associated polynomials of $\lbrace f_i \rbrace_{i=0}^{d+1}$.
 \item[\rm (ii)]
$(X^{-1})_{dj} = \frac{c_H}{\tau_{j}(\theta_j) \eta_{d-j}(\theta_j) X_{j0}}$ for $0 \leq j \leq d$, where $H$ is the connection coefficient matrix of $\lbrace f_i \rbrace_{i=0}^{d+1}$ and $c_H$ is from Notation \ref{not:cofh}. 
\end{enumerate}
\end{proposition}

\noindent {\it Proof:} (i) Let $D \in \Mdf$ denote the diagonal matrix with $(i,i)$-entry $\theta_i$ for $0 \leq i \leq d$. Observe that $H$ is Hessenberg by Definition \ref{def:polytohess} and that $\lbrace f_i \rbrace_{i=0}^{d+1}$ are the polynomials of $H$ by Lemma \ref{lem:polyscatohess}. Since $f_{d+1} = \prod_{i=0}^d (\lambda - \theta_i)$ and $\lbrace \theta_i \rbrace_{i=0}^d$ are mutually distinct, we find using Corollary \ref{cor:fdplusoneminpoly} that $H$ is multiplicity-free and $\lbrace \theta_i \rbrace_{i=0}^d$ is an ordering of the eigenvalues of $H$. Therefore $H = X^{-1} D X$ by Lemma \ref{lem:appbw}. Applying Lemma \ref{lem:appbs} to $X^{-1}$, we find $(X^{-1}, \lbrace \theta_i \rbrace_{i=0}^d)$ is a south Vandermonde system with polynomials $\lbrace f^{\varsigma}_i \rbrace_{i=0}^{d+1}$. \\
(ii) First assume that $X_{i0} = 1$ for $0 \leq i \leq d$. Let $h \in \K[\lambda]$ denote the polynomial $\sum_{i=0}^{d} (X^{-1})_{ij} f_i$. In the equation $XX^{-1} = I$, evaluate the $j^{th}$-column using matrix multiplication to find that $h(\theta_i) = \delta_{ij}$ for $0 \leq i \leq d$. Let $e_j \in \K[\lambda]$ denote the polynomial $\frac{\tau_j \eta_{d-j}}{\tau_j (\theta_j) \eta_{d-j}(\theta_j)}$. Observe that $e_j (\theta_i) = \delta_{ij}$ for $0 \leq i \leq d$. Thus $h(\theta_i) = e_j (\theta_i)$ for $0 \leq i \leq d$. It follows that $h = e_j$ since both $h$ and $e_j$ have degree $d$. In particular, the leading coefficient of $h$ is equal to the leading coefficient of $e_j$. By Lemma \ref{lem:polyofH}(i) the leading coefficient of $f_d$ is $c_H^{-1}$, so the leading coefficient of $h$ is $(X^{-1})_{dj} c_H^{-1}$. The leading coefficient of $e_j$ is $(\tau_j(\theta_j) \eta_{d-j}(\theta_j))^{-1}$. By these comments $(X^{-1})_{dj} c_H^{-1} = (\tau_j(\theta_j) \eta_{d-
 j}(\theta_j))^{-1}$ so $(X^{-1})_{dj} = \frac{c_H}{\tau_{j}(\theta_j) \eta_{d-j}(\theta_j)}$. The result is now proven for the special case in which $X_{i0} = 1$ for $0 \leq i \leq d$. For the general case,
apply the special case to the west Vandermonde system $(\Delta^{-1} X, \lbrace \theta_i \rbrace_{i=0}^d)$,
where $\Delta \in \Mdf$ is the diagonal matrix with $(i,i)$-entry $X_{i0}$ for $0 \leq i \leq d$. \hfill $\Box$

\begin{proposition}
\label{prop:inversevandisvand2}
Let $(X, \lbrace \theta_i \rbrace_{i=0}^d)$ denote a south Vandermonde system with polynomials $\lbrace f_i \rbrace_{i=0}^{d+1}$. Then the following {\rm (i)}, {\rm (ii)} hold. 
\begin{enumerate}
 \item[\rm (i)]
$(X^{-1}, \lbrace \theta_i \rbrace_{i=0}^d)$ is a west Vandermonde system with polynomials $\lbrace f^{\varsigma}_i \rbrace_{i=0}^{d+1}$, where $\lbrace f^{\varsigma}_i \rbrace_{i=0}^{d+1}$
are the associated polynomials of $\lbrace f_i \rbrace_{i=0}^{d+1}$.
 \item[\rm (ii)]
$(X^{-1})_{i0} = \frac{c_H}{\tau_{i}(\theta_i) \eta_{d-i}(\theta_i) X_{di}}$ for $0 \leq i \leq d$, where
$H$ is the connection coefficient matrix of $\lbrace f_i \rbrace_{i=0}^{d+1}$ and $c_H$ is from Notation \ref{not:cofh}. 
\end{enumerate}
\end{proposition}

\noindent {\it Proof:} Similar to the proof of Proposition \ref{prop:inversevandisvand}. \hfill $\Box$

\medskip

\noindent In the next section we return to our discussion of TH systems. 

\section{The transition matrices $P, \mathcal P$ and their Vandermonde structures}
\label{sec:transvand}
We return our attention to TH systems. Let $\Phi$ denote a TH system. Recall the transition matrices $P$ and $\mathcal P$ of $\Phi$ from Definition \ref{def:transmat} and Definition \ref{def:transmatnorm}. We will show that each of $P, \mathcal P$ has a west Vandermonde structure and a south Vandermonde structure. We start by associating with $\Phi$ a graded sequence of polynomials.

\begin{definition}
\label{def:polyofphi}
\rm
Let $\Phi=(A;\lbrace E_i\rbrace_{i=0}^d;A^*;\lbrace E^*_i\rbrace_{i=0}^d)$ denote a TH system on $V$. Let $\lbrace v_i \rbrace_{i=0}^d$ denote a $\Phi$-standard basis for $V$ and let $H \in \Mdf$ denote the matrix representing $A$ with respect to $\lbrace v_i \rbrace_{i=0}^d$. Observe that $H$ is Hessenberg. Let $\lbrace s_i \rbrace_{i=0}^{d+1}$ denote the polynomials of $H$ from Definition \ref{def:polyofhess}, so that $s_i(A)v_0 = v_i$ for $0 \leq i \leq d$ by Lemma \ref{lem:lab452} and $s_{d+1}$ is the minimal polynomial of $A$ by Corollary \ref{cor:fdplusoneminpoly}.
\end{definition}

\noindent The following normalization of the $\lbrace s_i\rbrace_{i=0}^{d+1}$ will be useful. 

\begin{definition}
\label{def:normpolyofphi}
\rm
With reference to Definition \ref{def:polyofphi}, let $\lbrace t_i \rbrace_{i=0}^{d+1}$ denote the sequence of polynomials in $\K[\lambda]$ that satisfies {\rm (i)}, {\rm (ii)} below. 
\begin{enumerate}
 \item[\rm (i)] 
For $0 \leq i \leq d$, $t_i = s_i/\ell_i$ where $\ell_i$ is from Definition \ref{def:kappas}. 
 \item[\rm (ii)]
$t_{d+1} = s_{d+1}$.
\end{enumerate}
We will show in Corollary \ref{lab:234} that $t_i(\theta_d) = 1$ for $0 \leq i \leq d$.
\end{definition}

In Definition \ref{def:polyofphi} we saw how the polynomials $\lbrace s_i\rbrace_{i=0}^{d+1}$ arise naturally from the action of $A$ on a $\Phi$-standard basis for $V$. We now discuss the meaning of the polynomials $\lbrace t_i \rbrace_{i=0}^{d+1}$ from this point of view. Let $\lbrace u_i \rbrace_{i=0}^d$ denote the inverted dual of a ${\tilde \Phi}$-standard basis for ${\tilde V}$. By Corollary \ref{cor:dualvstilde}, $\lbrace \ell_i u_i \rbrace_{i=0}^d$ is a $\Phi$-standard basis for $V$, where $\ell_i$ is from Definition \ref{def:kappas}. Therefore by Definition \ref{def:normpolyofphi}, $t_i(A)u_0 = u_i$ for $0 \leq i \leq d$ and $t_{d+1}$ is the minimal polynomial of $A$.  

Our next goal is to show that the polynomials $\lbrace s_i\rbrace_{i=0}^{d+1}$ and $\lbrace {\tilde t_i} \rbrace_{i=0}^{d+1}$ are associated in the sense of Definition \ref{def:dualpolynomial}. We will use the following fact. 
Let $\{ v_i \}_{i=0}^d$ denote a basis for $V$ and let $R \in {\rm End}(V)$. Let $S \in \Mdf$ denote the matrix representing $R$ with respect to $\{ v_i \}_{i=0}^d$. By elementary linear algebra, $S^t$ is the matrix representing $R^\sigma$ with respect to the dual of $\{ v_i \}_{i=0}^d$, where $\sigma: {\rm End}(V) \to {\rm End}(\tilde{V})$ is the canonical anti-isomorphism from above Definition \ref{def:dualphisetup}. Therefore the matrix $S^\varsigma$ represents $R^\sigma$ with respect to the inverted dual of $\{ v_i \}_{i=0}^d$.

\begin{lemma}
\label{lem:polyassociated}
With reference to Definition \ref{def:polyofphi} and Definition \ref{def:normpolyofphi}, for each column of the table below, the two graded sequences of polynomials are associated in the sense of Definition \ref{def:dualpolynomial}. 

\begin{center}
\begin{tabular}{c|c|c|c}
$\lbrace s_i\rbrace_{i=0}^{d+1}$ & $\lbrace s^{*}_i \rbrace_{i=0}^{d+1}$ & $\lbrace {\tilde s_i} \rbrace_{i=0}^{d+1}$ & $\lbrace {\tilde s^{*}_i} \rbrace_{i=0}^{d+1}$ \\
\hline
$\lbrace {\tilde t_i} \rbrace_{i=0}^{d+1}$ & $\lbrace {\tilde t_i^*} \rbrace_{i=0}^{d+1}$ & $\lbrace t_i \rbrace_{i=0}^{d+1}$ & $\lbrace t_i^* \rbrace_{i=0}^{d+1}$ 
\end{tabular}
\end{center}
\end{lemma}

\noindent {\it Proof:} Let $\lbrace v_i \rbrace_{i=0}^d$ denote a $\Phi$-standard basis for $V$ and let $H \in \Mdf$ denote the matrix representing $A$ with respect to $\lbrace v_i \rbrace_{i=0}^d$. By Definition \ref{def:polyofphi}, $\lbrace s_i \rbrace_{i=0}^{d+1}$ are the polynomials of $H$. 
Let $\lbrace w_i \rbrace_{i=0}^d$ denote the inverted dual  of $\{ v_i \}_{i=0}^d$. Applying the comments below Definition \ref{def:normpolyofphi} to ${\tilde \Phi}$, we find ${\tilde t_i}(A^\sigma)w_0 = w_i$ for $0 \leq i \leq d$ and ${\tilde t_{d+1}}$ is the minimal polynomial of $A^\sigma$. Moreover by the comments above the present lemma, the matrix $H^\varsigma$ represents $A^\sigma$ with respect to $\lbrace w_i \rbrace_{i=0}^d$. Now by Lemma \ref{lem:lab452} the $\lbrace {\tilde t_i} \rbrace_{i=0}^{d+1}$ are the polynomials of $H^\varsigma$. Therefore $\lbrace s_i \rbrace_{i=0}^{d+1}$ and $\lbrace {\tilde t_i} \rbrace_{i=0}^{d+1}$ are associated by Definition \ref{def:dualpolynomial}. We have verified our assertions about the first column of the above table. Our assertions about the remaining
columns follow from Definition \ref{def:conv2}. \hfill $\Box$

\medskip

We recall some elementary linear algebra. Let $\{ u_i \}_{i=0}^d$ and $\{ v_i \}_{i=0}^d$ denote bases for $V$. Let $T \in \Mdf$ denote the transition matrix from $\{ u_i \}_{i=0}^d$ to $\{ v_i \}_{i=0}^d$.  
Pick $A \in {\rm End}(V)$ and let $S \in \Mdf$ denote the matrix that represents $A$ with respect to $\lbrace u_i \rbrace_{i=0}^d$. Then the matrix $T^{-1} S T$ represents $A$ with respect to $\lbrace v_i \rbrace_{i=0}^d$.

\medskip

\noindent We now display a west Vandermonde structure for $P$. 

\begin{proposition}
\label{prop:pdoublevandwest}
Let $\Phi$ denote a TH system with eigenvalue sequence $\lbrace \theta_i \rbrace_{i=0}^d$ and dual eigenvalue sequence $\lbrace \theta^*_i \rbrace_{i=0}^d$.  Let $P$ denote the transition matrix of $\Phi$ from Definition \ref{def:transmat}. Then $(P, \lbrace \theta_i \rbrace_{i=0}^d)$ is a west Vandermonde system, and the corresponding polynomials are the $\lbrace s_i\rbrace_{i=0}^{d+1}$ from Definition \ref{def:polyofphi}. For each relative of $P$  we display a west Vandermonde system along with the corresponding polynomials. 

\begin{center}
\begin{tabular}{c|c}
west Vandermonde system & corresponding polynomials \\
\hline
$(P, \lbrace \theta_i \rbrace_{i=0}^d)$ & $\lbrace s_i\rbrace_{i=0}^{d+1}$ \\ 

$(P^{*}, \lbrace \theta^*_i \rbrace_{i=0}^d)$ & $\lbrace s^{*}_i \rbrace_{i=0}^{d+1}$ \\
 
$({\tilde P}, \lbrace \theta_{d-i} \rbrace_{i=0}^d)$ & $\lbrace {\tilde s_i} \rbrace_{i=0}^{d+1}$ \\ 

$({\tilde P^{*}}, \lbrace \theta^*_{d-i} \rbrace_{i=0}^d)$  & $\lbrace {\tilde s^{*}_i} \rbrace_{i=0}^{d+1}$
\end{tabular}
\end{center}
\end{proposition}

\noindent {\it Proof:} Write $\Phi=(A;\lbrace E_i\rbrace_{i=0}^d;A^*;\lbrace E^*_i\rbrace_{i=0}^d)$ and assume $V$ is the vector space underlying $\Phi$. Let $H \in \Mdf$ (resp. $D \in \Mdf$) denote the matrix representing $A$ with respect to a $\Phi$-standard (resp. $\Phi^*$-standard) basis for $V$. 
By construction $H$ is Hessenberg and multiplicity-free with an ordering of the eigenvalues $\lbrace \theta_i \rbrace_{i=0}^d$. By construction $D$ is diagonal with $(i,i)$-entry $\theta_i$ for $0 \leq i \leq d$. By Definition \ref{def:polyofphi} $\lbrace s_i\rbrace_{i=0}^{d+1}$ are the polynomials of $H$. By Definition \ref{def:transmat} and the comments above this proposition, we have $H = P^{-1} D P$. Therefore by Lemma \ref{lem:appbw} $(P, \lbrace \theta_i \rbrace_{i=0}^d)$ is a west Vandermonde system with polynomials $\lbrace s_i\rbrace_{i=0}^{d+1}$. We have verified our assertions about the first row of the above table. Our assertions about the remaining rows follow from Corollary \ref{cor:parelatives}. \hfill $\Box$

\medskip

\noindent We now display a south Vandermonde structure for $P$. 

\begin{proposition}
\label{prop:pdoublevandsouth}
Let $\Phi$ denote a TH system with eigenvalue sequence $\lbrace \theta_i \rbrace_{i=0}^d$ and dual eigenvalue sequence $\lbrace \theta^*_i \rbrace_{i=0}^d$.  Let $P$ denote the transition matrix of $\Phi$ from Definition \ref{def:transmat}.   
Then $(P, \lbrace \theta^*_i \rbrace_{i=0}^d)$ is a south Vandermonde system, and the corresponding polynomials are the $\lbrace {\tilde t^*_i} \rbrace_{i=0}^{d+1}$ from Definition \ref{def:normpolyofphi}. For each relative of $P$ we display a south Vandermonde system along with the corresponding polynomials. 

\begin{center}
\begin{tabular}{c|c}
south Vandermonde system & corresponding polynomials \\
\hline
$(P, \lbrace \theta^*_i \rbrace_{i=0}^d)$ & $\lbrace {\tilde t^*_i} \rbrace_{i=0}^{d+1}$ \\ 

$(P^{*}, \lbrace \theta_i \rbrace_{i=0}^d)$ & $\lbrace {\tilde t_i} \rbrace_{i=0}^{d+1}$ \\
 
$({\tilde P}, \lbrace \theta^*_{d-i} \rbrace_{i=0}^d)$ & $\lbrace t^*_i \rbrace_{i=0}^{d+1}$ \\ 

$({\tilde P^{*}}, \lbrace \theta_{d-i} \rbrace_{i=0}^d)$  & $\lbrace t_i \rbrace_{i=0}^{d+1}$
\end{tabular}
\end{center}
\end{proposition}

\noindent {\it Proof:} By Proposition \ref{prop:pdoublevandwest} $(P^*, \lbrace \theta_i^* \rbrace_{i=0}^d)$ is a west Vandermonde system with polynomials $\lbrace s_i^* \rbrace_{i=0}^{d+1}$. Thus by Proposition \ref{prop:inversevandisvand} and Lemma \ref{lem:polyassociated}, $((P^*)^{-1}, \lbrace \theta_i^* \rbrace_{i=0}^d)$ is a south Vandermonde system with polynomials $\lbrace {\tilde t_i^*} \rbrace_{i=0}^{d+1}$. By this and since $P P^*= \nu I$, $(\nu^{-1} P, \lbrace \theta^*_i \rbrace_{i=0}^d)$ is a south Vandermonde system with  polynomials $\lbrace {\tilde t^*_i} \rbrace_{i=0}^{d+1}$. Therefore by Lemma \ref{lem:westvandunique} $(P, \lbrace \theta^*_i \rbrace_{i=0}^d)$ is a south Vandermonde system with  polynomials $\lbrace {\tilde t^*_i} \rbrace_{i=0}^{d+1}$. We have verified our assertions about the first row of the above table. Our assertions about the remaining rows follow from Corollary \ref{cor:parelatives}. \hfill $\Box$

\medskip

We now turn to the matrix $\mathcal P$. Below we display a west Vandermonde structure and a south Vandermonde structure for $\mathcal P$. We begin with the west Vandermonde structure.

\begin{corollary}
\label{cor:pnormdoublevandwest}
Let $\Phi$ denote a TH system with eigenvalue sequence $\lbrace \theta_i \rbrace_{i=0}^d$ and dual eigenvalue sequence $\lbrace \theta^*_i \rbrace_{i=0}^d$.  Let $\mathcal P$ denote the transition matrix of $\Phi$ from Definition \ref{def:transmatnorm}. Then $(\mathcal P, \lbrace \theta_i \rbrace_{i=0}^d)$ is a west Vandermonde system, and the corresponding polynomials are the $\lbrace t_i \rbrace_{i=0}^{d+1}$ from Definition \ref{def:normpolyofphi}. For each relative of $\mathcal P$ we display a west Vandermonde system along with the corresponding polynomials. 

\begin{center}
\begin{tabular}{c|c}
west Vandermonde system & corresponding polynomials \\
\hline
$(\mathcal P, \lbrace \theta_i \rbrace_{i=0}^d)$ & $\lbrace t_i \rbrace_{i=0}^{d+1}$ \\ 

$(\mathcal P^{*}, \lbrace \theta^*_i \rbrace_{i=0}^d)$ & $\lbrace t^{*}_i \rbrace_{i=0}^{d+1}$ \\
 
$({\tilde \mathcal P}, \lbrace \theta_{d-i} \rbrace_{i=0}^d)$ & $\lbrace {\tilde t_i} \rbrace_{i=0}^{d+1}$ \\ 

$({\tilde \mathcal P^{*}}, \lbrace \theta^*_{d-i} \rbrace_{i=0}^d)$  & $\lbrace {\tilde t^{*}_i} \rbrace_{i=0}^{d+1}$
\end{tabular}
\end{center}
\end{corollary}

\noindent {\it Proof:} Routine by Lemma \ref{lem:westvandright} and Proposition \ref{prop:pdoublevandwest}. \hfill $\Box$

\medskip

\noindent We now display a south Vandermonde structure for $\mathcal P$.

\begin{corollary}
\label{cor:pnormdoublevandsouth}
Let $\Phi$ denote a TH system with eigenvalue sequence $\lbrace \theta_i \rbrace_{i=0}^d$ and dual eigenvalue sequence $\lbrace \theta^*_i \rbrace_{i=0}^d$. Let $\mathcal P$ denote the transition matrix of $\Phi$ from Definition \ref{def:transmatnorm}.
Then $(\mathcal P, \lbrace \theta^*_i \rbrace_{i=0}^d)$ is a south Vandermonde system, and the corresponding polynomials are the $\lbrace {\tilde t^*_i} \rbrace_{i=0}^{d+1}$ from Definition \ref{def:normpolyofphi}. For each relative of $\mathcal P$ we display a south Vandermonde system along with the corresponding polynomials. 

\begin{center}
\begin{tabular}{c|c}
south Vandermonde system & corresponding polynomials \\
\hline
$(\mathcal P, \lbrace \theta^*_i \rbrace_{i=0}^d)$ & $\lbrace {\tilde t^*_i} \rbrace_{i=0}^{d+1}$ \\ 

$(\mathcal P^{*}, \lbrace \theta_i \rbrace_{i=0}^d)$ & $\lbrace {\tilde t_i} \rbrace_{i=0}^{d+1}$ \\
 
$({\tilde \mathcal P}, \lbrace \theta^*_{d-i} \rbrace_{i=0}^d)$ & $\lbrace t^*_i \rbrace_{i=0}^{d+1}$ \\ 

$({\tilde \mathcal P^{*}}, \lbrace \theta_{d-i} \rbrace_{i=0}^d)$  & $\lbrace t_i \rbrace_{i=0}^{d+1}$
\end{tabular}
\end{center}
\end{corollary}

\noindent {\it Proof:} Similar to the proof of Corollary \ref{cor:pnormdoublevandwest} using Proposition \ref{prop:pdoublevandsouth}. \hfill $\Box$

\medskip

\noindent We have now displayed the west Vandermonde and south Vandermonde structures for $P$ and $\mathcal P$. As corollaries to these results, we now obtain some facts involving the polynomials $\lbrace s_i \rbrace_{i=0}^{d+1}$ and $\lbrace t_i \rbrace_{i=0}^{d+1}$. 

\begin{corollary}
\label{cor:aw}
Let $\Phi$ denote a TH system with eigenvalue sequence $\lbrace \theta_i \rbrace_{i=0}^d$ and dual eigenvalue sequence $\lbrace \theta^*_i \rbrace_{i=0}^d$.  Let $P$  denote the transition matrix of $\Phi$ from Definition \ref{def:transmat}  and let $\mathcal P$ denote the corresponding matrix from Definition \ref{def:transmatnorm}. Let $\lbrace s_i \rbrace_{i=0}^{d+1}$  (resp.  $\lbrace t_i \rbrace_{i=0}^{d+1}$) denote the polynomials of $\Phi$ from Definition \ref{def:polyofphi} (resp. Definition \ref{def:normpolyofphi}). Then the following {\rm (i)}, {\rm (ii)} hold  for $0 \leq i,j \leq d$. 
\begin{enumerate}
 \item[\rm (i)]
$P_{ij} = \ell_j t_j(\theta_i) = \ell_j {\tilde t^*_{d-i}}(\theta_j^*) = s_j(\theta_i) = \ell_j {\tilde s^*_{d-i}}(\theta_j^*)/{\tilde \ell^*_{d-i}}$. 
 \item[\rm (ii)]
$\mathcal P_{ij} = t_j(\theta_i) = {\tilde t^*_{d-i}}(\theta_j^*) = s_j(\theta_i)/\ell_j = {\tilde s^*_{d-i}}(\theta_j^*)/{\tilde \ell^*_{d-i}}$.
\end{enumerate}
Here $\ell_j$, ${\tilde \ell_{d-i}^*}$ are from Definition \ref{def:kappas} and Lemma \ref{lem:ellistar} respectively.
\end{corollary}

\noindent {\it Proof:} (i) Using Corollary \ref{cor:pentries}, (\ref{eq:vand}), and Proposition \ref{prop:pdoublevandwest}, we find $P_{ij} = s_j(\theta_i)$. 
Similarly using Proposition \ref{prop:pdoublevandsouth} in place of Proposition \ref{prop:pdoublevandwest} we find $P_{ij} = \ell_j {\tilde t^*_{d-i}}(\theta_j^*)$. The remaining assertions follow using Definition \ref{def:normpolyofphi}. \\
(ii) Use (i) and the fact that $P_{ij} = \mathcal P_{ij} \ell_j$ for $0 \leq i,j \leq d$. \hfill $\Box$

\medskip

We emphasize one aspect of Corollary \ref{cor:aw} which is telling us that the $\lbrace s_i \rbrace_{i=0}^{d+1}$ and the $\lbrace t_i \rbrace_{i=0}^{d+1}$ each satisfy a variation on the Askey-Wilson duality \cite[Theorems 14.7--14.9]{qrac}. 

\begin{corollary}
Let $\Phi$ denote a TH system with eigenvalue sequence $\lbrace \theta_i \rbrace_{i=0}^d$ and dual eigenvalue sequence $\lbrace \theta^*_i \rbrace_{i=0}^d$. Let $\lbrace s_i \rbrace_{i=0}^{d+1}$ (resp. $\lbrace t_i \rbrace_{i=0}^{d+1}$) denote the corresponding polynomials from Definition \ref{def:polyofphi} (resp. Definition \ref{def:normpolyofphi}). Then the following {\rm (i)}, {\rm (ii)} hold  for $0 \leq i,j \leq d$. 
\begin{enumerate}
\item[\rm (i)] 
$t_j(\theta_i) = {\tilde t^*_{d-i}}(\theta_j^*)$. 
\item[\rm (ii)]
$s_j(\theta_i)/\ell_j = {\tilde s^*_{d-i}}(\theta_j^*)/{\tilde \ell^*_{d-i}}$, where $\ell_j$, ${\tilde \ell_{d-i}^*}$ are from Definition \ref{def:kappas} and Lemma \ref{lem:ellistar} respectively.
\end{enumerate}
\end{corollary}

\noindent We have a comment on how the polynomials $\lbrace s_i \rbrace_{i=0}^{d+1}$ and  $\lbrace t_i \rbrace_{i=0}^{d+1}$ are normalized. 

\begin{corollary}
\label{lab:234}
Let $\Phi$ denote a TH system with eigenvalue sequence $\lbrace \theta_i \rbrace_{i=0}^d$ and dual eigenvalue sequence $\lbrace \theta^*_i \rbrace_{i=0}^d$. Let $\lbrace s_i \rbrace_{i=0}^{d+1}$ (resp. $\lbrace t_i \rbrace_{i=0}^{d+1}$) denote the corresponding polynomials from Definition \ref{def:polyofphi} (resp. Definition \ref{def:normpolyofphi}). Then the following {\rm (i)}, {\rm (ii)} hold. 
\begin{enumerate}
 \item[\rm (i)]
$s_i(\theta_d) = \ell_i$ for $0 \leq i \leq d$,  where $\ell_i$ is from Definition \ref{def:kappas}. 
 \item[\rm (ii)]
$t_i(\theta_d) = 1$ for $0 \leq i \leq d$.
\end{enumerate}
\end{corollary}

\noindent {\it Proof:} Use Corollary \ref{cor:pentries} and Corollary \ref{cor:aw}(i). \hfill $\Box$

\medskip

\noindent The polynomials $\lbrace t_i \rbrace_{i=0}^{d+1}$ are not orthogonal in general; however we do have the following. 

\begin{corollary}
\label{cor:viorth}
Let $\Phi$ denote a TH system with eigenvalue sequence $\lbrace \theta_i \rbrace_{i=0}^d$ and dual eigenvalue sequence $\lbrace \theta^*_i \rbrace_{i=0}^d$. Let $\lbrace t_i \rbrace_{i=0}^{d+1}$ denote the corresponding polynomials from Definition \ref{def:normpolyofphi}. Then the following {\rm (i), (ii)} hold. 
\begin{enumerate}
\item[\rm (i)]
${\displaystyle \sum_{n=0}^d t_i(\theta_n) {\tilde t_{j}}(\theta_n) \ell^*_n = \delta_{i+j,d} \, \nu \ell_i^{-1} \qquad \qquad \  \ \ \ \ (0 \leq i,j \leq d)}$.                     
\item[\rm (ii)]
${\displaystyle \sum_{i=0}^d  t_{i}(\theta_m) {\tilde t_{d-i}}(\theta_n) \ell_{i}  = \delta_{mn} \, \nu (\ell_m^{*})^{-1}   \qquad \qquad (0 \leq m,n \leq d)}$.                    
\end{enumerate}                            
Here $\{ \ell_i \}_{i=0}^d$, $\{ \ell_i^* \}_{i=0}^d$ are from Definition \ref{def:kappas} and Lemma \ref{lem:ellistar} respectively, and $\nu$ is from Definition \ref{def:nu}. 
\end{corollary} 

\noindent {\it Proof:} (i) Let $P$ denote the transition matrix of $\Phi$ from Definition \ref{def:transmat}. In the equation $P^*P = \nu I$, compare the $(d-j,i)$-entry of each side to obtain
$\sum_{n=0}^d P^*_{d-j,n} P_{ni} = \delta_{i+j,d} \nu$. In this equation, evaluate $P_{ni}$ and $P^*_{d-j,n}$ using Corollary \ref{cor:aw} to obtain $\sum_{n=0}^d \ell^*_n {\tilde t_{j}}(\theta_n) \ell_i t_{i}(\theta_n) = \delta_{i+j,d} \nu$.  The result follows. \\
\noindent (ii) Let $P$ denote the transition matrix of $\Phi$ from Definition \ref{def:transmat}. In the equation $PP^* = \nu I$,  compare the $(m,n)$-entry of each side to obtain
$\sum_{i=0}^d P_{mi} P^*_{in} = \delta_{mn} \nu$. In this equation, evaluate $P_{mi}$ and $P^*_{in}$ using Corollary \ref{cor:aw} to obtain  $\sum_{i=0}^d \ell_i t_i(\theta_m) \ell^*_n {\tilde t_{d-i}}(\theta_n) = \delta_{mn} \nu$. The result follows. \hfill $\Box$

\medskip

\noindent We now give an analogue of Corollary \ref{cor:viorth} that applies to the polynomials $\lbrace s_i\rbrace_{i=0}^{d+1}$. 

\begin{corollary}
\label{cor:viorth2}
Let $\Phi$ denote a TH system with eigenvalue sequence $\lbrace \theta_i \rbrace_{i=0}^d$ and dual eigenvalue sequence $\lbrace \theta^*_i \rbrace_{i=0}^d$. Let $\lbrace s_i \rbrace_{i=0}^{d+1}$ denote the corresponding polynomials from Definition \ref{def:polyofphi}. Then the following {\rm (i), (ii)} hold.
\begin{enumerate}
\item[\rm (i)]
${\displaystyle \sum_{n=0}^d s_i(\theta_n) {\tilde s_{j}}(\theta_n)\ell^*_n = \delta_{i+j,d} \, \nu   {\tilde \ell_{j}} \qquad \qquad \qquad \qquad (0 \leq i,j \leq d)}$.                                      \item[\rm (ii)]
${\displaystyle \sum_{i=0}^d s_{i}(\theta_m) {\tilde s_{d-i}}(\theta_n) ({\tilde \ell_{d-i}})^{-1} = \delta_{mn} \, \nu (\ell_m^{*})^{-1}  \qquad \qquad (0 \leq m,n \leq d)}$.                             
\end{enumerate}
Here $\{{\tilde \ell_i} \}_{i=0}^d$, $\{ \ell_i^* \}_{i=0}^d$ are from Lemma \ref{lem:ellistar} and $\nu$ is from Definition \ref{def:nu}.  
\end{corollary}

\noindent {\it Proof:} Use Definition \ref{def:normpolyofphi}(i) and Corollary \ref{cor:viorth}. \hfill $\Box$

\medskip

\noindent We now express the polynomials $\lbrace t_i \rbrace_{i=0}^{d+1}$ and $\lbrace s_i\rbrace_{i=0}^{d+1}$ in terms of the parameter array of $\Phi$. To do this we will use the two-variable polynomial $p$ of $\Phi$ from Definition \ref{def:scriptp}. 

\begin{corollary}
\label{cor:pivsscriptp}
Let $\Phi$ denote a TH system with eigenvalue sequence $\lbrace \theta_i \rbrace_{i=0}^d$ and dual eigenvalue sequence $\lbrace \theta^*_i \rbrace_{i=0}^d$. Let $\lbrace t_i \rbrace_{i=0}^{d+1}$ denote the corresponding polynomials from Definition \ref{def:normpolyofphi}. Then the following {\rm (i)}, {\rm (ii)} hold. 
\begin{enumerate}
 \item[\rm (i)]
For $0 \leq i \leq d$, $t_i= p(\lambda, \theta^*_i)$ where $p$ is from Definition \ref{def:scriptp}.
 \item[\rm (ii)]
$t_{d+1} = \prod_{i=0}^d (\lambda - \theta_i)$. 
\end{enumerate}
\end{corollary}

\noindent {\it Proof:} Let $\mathcal P$ denote the transition matrix of $\Phi$ from Definition \ref{def:transmatnorm}. Since $\mathcal P_{ij} = p(\theta_i, \theta_j^*)$ for $0 \leq i, j \leq d$, we find that $(\mathcal P, \lbrace \theta_i \rbrace_{i=0}^d)$ is a west Vandermonde system with polynomials $\{ f_i \}_{i=0}^{d+1}$ where $f_i = p(\lambda, \theta^*_i)$ for $0 \leq i \leq d$ and $f_{d+1} = \prod_{i=0}^d (\lambda - \theta_i)$. The result follows by Corollary \ref{cor:pnormdoublevandwest}. \hfill $\Box$

\medskip

\begin{example}
\rm
With reference to Definition \ref{def:normpolyofphi}, assume $d = 2$. Then 
\begin{eqnarray*}
t_0 &=& 1, \\
t_1 &=& 1 +  \frac{(\lambda - \theta_2)(\theta_1^* - \theta_0^*)}{\phi_1}, \\
t_2 &=& 1 +  \frac{(\lambda - \theta_2)(\theta_2^* - \theta_0^*)}{\phi_1} +  \frac{(\lambda - \theta_2)(\lambda - \theta_1)(\theta_2^* - \theta_0^*)(\theta_2^* - \theta_1^*)}{\phi_1 \phi_2}.
\end{eqnarray*}
\end{example}

\begin{corollary}
\label{cor:qivsscriptp}
Let $\Phi$ denote a TH system with eigenvalue sequence $\lbrace \theta_i \rbrace_{i=0}^d$ and dual eigenvalue sequence $\lbrace \theta^*_i \rbrace_{i=0}^d$. Let $\lbrace s_i \rbrace_{i=0}^{d+1}$ denote the corresponding polynomials from Definition \ref{def:polyofphi}. Then the following {\rm (i)}, {\rm (ii)} hold. 
\begin{enumerate}
 \item[\rm (i)] 
For $0 \leq i \leq d$, $s_i= \ell_i p(\lambda, \theta^*_i)$ where $\ell_i$ is from Definition \ref{def:kappas} and $p$ is from Definition \ref{def:scriptp}.
 \item[\rm (ii)]
$s_{d+1} = \prod_{i=0}^d (\lambda - \theta_i)$. 
\end{enumerate}
\end{corollary}
 
\noindent {\it Proof:} Use Definition \ref{def:normpolyofphi} and Corollary \ref{cor:pivsscriptp}. \hfill $\Box$

\medskip

\begin{remark}
\rm
In view of Corollary \ref{cor:pivsscriptp}, one may wonder about the polynomial $p(\theta_i, \lambda)$. By Lemma \ref{lem:twovarpolydualrel} and Corollary \ref{cor:pivsscriptp}, $p(\theta_i, \lambda) = {\tilde p^*}(\lambda, \theta_i) = {\tilde t^*_{d-i}}$ for $0 \leq i \leq d$.
\end{remark}

\noindent We now give the results promised at the end of Section \ref{sec:bilform}. Let $\Phi=(A;\lbrace E_i\rbrace_{i=0}^d;A^*;\lbrace E^*_i\rbrace_{i=0}^d)$ denote a TH system on $V$. 
Let $0 \neq \xi_0 \in E_0V$ and recall the $\Phi$-standard basis $\lbrace E^*_i \xi_0 \rbrace_{i=0}^d$ for $V$ from above (\ref{eq:defui}).  
Let $0 \neq {\tilde \xi^*_d} \in E^{*\sigma}_d{\tilde V}$ and recall the ${\tilde \Phi^*}$-standard basis $\lbrace E^{\sigma}_{d-i} {\tilde \xi^*_d} \rbrace_{i=0}^d$ for ${\tilde V}$ from above Proposition \ref{prop:orthogstar}. These two bases are related as follows. 

\begin{proposition}
\label{prop:basiscross}
With reference to the TH system $\Phi$ in Definition \ref{def:dualphisetup}, let $0 \neq \xi_0 \in E_0V$ and $0 \neq {\tilde \xi_d^*} \in E^{*\sigma}_d{\tilde V}$. Then for $0 \leq i,j \leq d$,
\begin{eqnarray*}
\langle E_i^*\xi_0,E^{\sigma}_j {\tilde \xi_d^*} \rangle = \nu^{-1} \ell_i \ell_j^* t_{i}(\theta_j) \b{\xi_0, {\tilde \xi^*_d}}.
\end{eqnarray*}
Here $\ell_i$, $\ell^*_j$ are from Definition \ref{def:kappas} and Lemma \ref{lem:ellistar} respectively, and $\nu$, $t_i$ are from Definition \ref{def:nu} and Definition \ref{def:normpolyofphi} respectively. 
\end{proposition}

\noindent {\it Proof:} Let $P$ denote the transition matrix of $\Phi$ from Definition \ref{def:transmat}. Let $\xi_0^* = E_0^* \xi_0$ and observe by Lemma \ref{lem:nusig} that
\begin{eqnarray}
\label{lab:tempo}
E_0 \xi^*_0 = E_0 E_0^* \xi_0 = E_0 E_0^* E_0 \xi_0 = \nu^{-1} E_0 \xi_0 = \nu^{-1} \xi_0.
\end{eqnarray}
We may now argue 
\begin{eqnarray*}
\langle E_i^*\xi_0,E^{\sigma}_j {\tilde \xi_d^*}   \rangle
&=&
\sum_{n=0}^d P_{ni} \langle E_n\xi_0^*,E^{\sigma}_j {\tilde \xi_d^*}  \rangle   \qquad \qquad \ \ \  (\hbox{by Definition \ref{def:transmat}}) \\
&=&
\ell_j^* P_{ji} \b{E_0 \xi^*_0, {\tilde \xi^*_d}}  \qquad \qquad \qquad \ \ \  (\hbox{by Proposition \ref{prop:orthogstar}})  \\
&=&
\ell_j^* s_{i}(\theta_j) \b{E_0 \xi^*_0, {\tilde \xi^*_d}} \qquad \qquad \ \ \ \ \ \  (\hbox{by  Corollary \ref{cor:aw}}) \\
&=&
\nu^{-1} \ell_j^* s_{i}(\theta_j) \b{\xi_0, {\tilde \xi^*_d}} \qquad \qquad \ \ \ \ \   (\hbox{by  (\ref{lab:tempo})}) \\
&=&
\nu^{-1} \ell_i \ell_j^* t_{i}(\theta_j) \b{\xi_0, {\tilde \xi^*_d}} \qquad \qquad \ \ \  (\hbox{by  Definition \ref{def:normpolyofphi}}).
\end{eqnarray*} \hfill $\Box$ 

Let $\Phi=(A;\lbrace E_i\rbrace_{i=0}^d;A^*;\lbrace E^*_i\rbrace_{i=0}^d)$ denote a TH system on $V$. Let $0 \neq \xi^*_0 \in E^*_0V$ and recall the $\Phi^*$-standard basis $\lbrace E_i \xi^*_0 \rbrace_{i=0}^d$ for $V$ from above Proposition \ref{prop:orthogstar}.  
Let $0 \neq {\tilde \xi_d} \in E^{\sigma}_d{\tilde V}$ and recall the ${\tilde \Phi}$-standard basis $\lbrace E^{*\sigma}_{d-i} {\tilde \xi_d} \rbrace_{i=0}^d$ for ${\tilde V}$ from above Proposition \ref{prop:orthog}. These two bases are related as follows. 

\begin{proposition}
With reference to the TH system $\Phi$ in Definition \ref{def:dualphisetup}, let $0 \neq \xi_0^* \in E_0^*V$ and $0 \neq {\tilde \xi_d} \in E^{\sigma}_d{\tilde V}$. Then for $0 \leq i,j \leq d$, 
\begin{eqnarray*}
\langle E_i\xi_0^*,E^{*\sigma}_j {\tilde \xi_d} \rangle = \nu^{-1} \ell_i^* \ell_j t_{i}^*(\theta_j^*) \b{\xi_0^*, {\tilde \xi_d}}.
\end{eqnarray*}
Here $\ell_i^*, \ell_j$ are from Lemma \ref{lem:ellistar} and Definition \ref{def:kappas} respectively, and $\nu$,  $t_i^*$ are from Definition \ref{def:nu} and Definition \ref{def:normpolyofphi} respectively. 
\end{proposition}

\noindent {\it Proof:} Apply Proposition \ref{prop:basiscross} to $\Phi^*$. \hfill $\Box$

\section{TH systems and Vandermonde systems}
In the previous sections we discussed TH systems and Vandermonde systems. In this section we give a natural correspondence between these two objects.

\begin{definition}
\label{def:normalizedWSVand}
\rm
A matrix $X \in \Mdf$ is called {\it west-south {\rm (or {\it double})} Vandermonde} whenever $X$ is both west Vandermonde and south Vandermonde. Assume $X$ is west-south Vandermonde. We say that $X$ is {\it west} (resp. {\it south}) {\it normalized} whenever $X_{i0} = 1$ (resp. $X_{di} = 1$) for $0 \leq i \leq d$. We say that $X$ is {\it normalized} whenever it is both west normalized and south normalized.
\end{definition}

\begin{definition}
\label{def:wsvandsystem}
\rm
By a {\it west-south {\rm (or {\it double})} Vandermonde system} in $\Mdf$, we mean a sequence $(X, \lbrace \theta_i \rbrace_{i=0}^d, \lbrace \theta^*_i \rbrace_{i=0}^d)$ such that $(X, \lbrace \theta_i \rbrace_{i=0}^d)$ is a west Vandermonde system in $\Mdf$ and $(X, \lbrace \theta^*_i \rbrace_{i=0}^d)$ is a south Vandermonde system in $\Mdf$. Let $(X, \lbrace \theta_i \rbrace_{i=0}^d, \lbrace \theta^*_i \rbrace_{i=0}^d)$ denote a west-south Vandermonde system. Observe that $X$ is west-south Vandermonde. We say that $(X, \lbrace \theta_i \rbrace_{i=0}^d, \lbrace \theta^*_i \rbrace_{i=0}^d)$ is {\it west normalized} (resp. {\it south normalized}) (resp. {\it normalized}) whenever $X$ is west normalized (resp. south normalized) (resp.  normalized) in the sense of Definition \ref{def:normalizedWSVand}. 
\end{definition}

\noindent Our main goal in this section is to establish a bijection between the following two sets:
\begin{eqnarray}
\label{eq:thsystem}
\hbox{the set of isomorphism classes of TH systems over} \ \K \hbox{\ of diameter} \ d, \\
\label{eq:wsvand}
\hbox{the set of normalized west-south Vandermonde systems in} \ \Mdf.
\end{eqnarray}

\noindent To do this we define a map $\rho$ from (\ref{eq:thsystem}) to (\ref{eq:wsvand}) and a map $\chi$ from (\ref{eq:wsvand}) to (\ref{eq:thsystem}), and show that they are inverses of each other. We start with an observation. Let $\Phi$ denote a TH system over $\K$ of diameter $d$, with eigenvalue sequence $\lbrace \theta_i \rbrace_{i=0}^d$ and dual eigenvalue sequence $\lbrace \theta^*_i \rbrace_{i=0}^d$. Let $\mathcal P$ denote the transition matrix of $\Phi$ from Definition \ref{def:transmatnorm}. By Corollary \ref{cor:pnormdoublevandwest} the sequence $(\mathcal P, \lbrace \theta_i \rbrace_{i=0}^d)$ is a west Vandermonde system in $\Mdf$, and by Corollary \ref{cor:pnormdoublevandsouth} the sequence $(\mathcal P, \lbrace \theta^*_i \rbrace_{i=0}^d)$ is a south Vandermonde system in $\Mdf$. 
By Corollary \ref{cor:scriptpendpoint}, $\mathcal P$ is both west normalized and south normalized. Therefore $(\mathcal P, \lbrace \theta_i \rbrace_{i=0}^d, \lbrace \theta^*_i \rbrace_{i=0}^d)$ is a normalized west-south Vandermonde system in $\Mdf$. 

\begin{definition}
\label{def:rho}
\rm
We define a map $\rho$ from (\ref{eq:thsystem}) to (\ref{eq:wsvand}). We do this as follows. 
Let $\Phi$ denote a TH system over $\K$ of diameter $d$, with eigenvalue sequence $\lbrace \theta_i \rbrace_{i=0}^d$ and dual eigenvalue sequence $\lbrace \theta^*_i \rbrace_{i=0}^d$. Let $\mathcal P$ denote the transition matrix of $\Phi$ from Definition \ref{def:transmatnorm}. By the above comment $(\mathcal P, \lbrace \theta_i \rbrace_{i=0}^d, \lbrace \theta^*_i \rbrace_{i=0}^d)$ is a normalized west-south Vandermonde system in $\Mdf$. The map $\rho$ sends the isomorphism class of $\Phi$ to $(\mathcal P, \lbrace \theta_i \rbrace_{i=0}^d, \lbrace \theta^*_i \rbrace_{i=0}^d)$. 
\end{definition}

We now define the map $\chi$. We start with the following construction.  Let $(X, \lbrace \theta_i \rbrace_{i=0}^d, \lbrace \theta^*_i \rbrace_{i=0}^d)$ denote a normalized west-south Vandermonde system in $\Mdf$. We construct a TH system $\Phi$ on $V$ as follows. Recall that $X$ is invertible by Lemma \ref{lem:vandinvertible}. Therefore there exist bases $\lbrace u_i \rbrace_{i=0}^d$, $\lbrace v_i \rbrace_{i=0}^d$ for $V$ such that $X$ is the transition matrix from $\lbrace u_i \rbrace_{i=0}^d$ to $\lbrace v_i \rbrace_{i=0}^d$.  Define $A \in {\rm End}(V)$ such that $A u_i = \theta_i u_i$ for $0 \leq i \leq d$. Define $A^* \in {\rm End}(V)$ such that $A^* v_i = \theta^*_i v_i$ for $0 \leq i \leq d$. For $0 \leq i \leq d$ let $E_i$ (resp. $E_i^*$) denote the primitive idempotent of $A$ (resp. $A^*$) corresponding to $\theta_i$ (resp. $\theta_i^*$). Now define $\Phi = (A;\lbrace E_i\rbrace_{i=0}^d;A^*;\lbrace E^*_i\rbrace_{i=0}^d)$. We claim that $\Phi$ is a TH system on $V$. To prove the claim we show that $\Phi$ satisfies conditions (i)--(v) in Definition \ref{def:HS}. Conditions (i)--(iii) hold by construction. Let $D^* \in \Mdf$ denote the diagonal matrix with $(i,i)$-entry $\theta^*_i$ for $0 \leq i \leq d$. Observe that $D^*$ represents $A^*$ with respect to $\lbrace v_i \rbrace_{i=0}^d$. Hence by the comment above Proposition \ref{prop:pdoublevandwest}, the matrix $X D^* X^{-1}$ represents $A^*$ with respect to $\lbrace u_i \rbrace_{i=0}^d$.
Moreover since $(X, \lbrace \theta_i^* \rbrace_{i=0}^d)$ is a south Vandermonde system, $X D^* X^{-1}$ is Hessenberg by Corollary \ref{cor:appbstwist}. By these comments, condition (iv) holds. Let $D \in \Mdf$ denote the diagonal matrix with $(i,i)$-entry $\theta_i$ for $0 \leq i \leq d$. Observe that $D$ represents $A$ with respect to $\lbrace u_i \rbrace_{i=0}^d$. Hence by the comment above Proposition \ref{prop:pdoublevandwest}, the matrix $X^{-1} D X$ represents $A$ with respect to $\lbrace v_i \rbrace_{i=0}^d$. 
Moreover since $(X, \lbrace \theta_i \rbrace_{i=0}^d)$ is a west Vandermonde system, $X^{-1} D X$ is Hessenberg by Corollary \ref{cor:appbwtwist}. By these comments, condition (v) holds. Therefore $\Phi$ is a TH system on $V$. By  construction $\Phi$ has eigenvalue sequence $\lbrace \theta_i \rbrace_{i=0}^d$ and dual eigenvalue sequence $\lbrace \theta_i^* \rbrace_{i=0}^d$. 

\begin{definition}
\label{def:pi}
\rm
We define a map $\chi$ from (\ref{eq:wsvand}) to (\ref{eq:thsystem}). We do this as follows. 
Let $(X, \lbrace \theta_i \rbrace_{i=0}^d, \lbrace \theta^*_i \rbrace_{i=0}^d)$ denote a normalized west-south Vandermonde system in $\Mdf$. Let $\Phi$ denote the corresponding TH system constructed above. The map $\chi$ sends $(X, \lbrace \theta_i \rbrace_{i=0}^d, \lbrace \theta^*_i \rbrace_{i=0}^d)$ to the isomorphism class of $\Phi$. 
\end{definition}

\noindent Our next goal is to show that the maps $\rho$ and $\chi$ are inverses of each other. We first recall some elementary linear algebra. Let $\lbrace u_i \rbrace_{i=0}^d$, $\lbrace v_i \rbrace_{i=0}^d$, $\lbrace w_i \rbrace_{i=0}^d$ denote bases for $V$. Let $T \in \Mdf$ denote the transition matrix from $\lbrace u_i \rbrace_{i=0}^d$ to $\lbrace v_i \rbrace_{i=0}^d$ and let $S \in \Mdf$ denote the transition matrix from $\lbrace v_i \rbrace_{i=0}^d$ to $\lbrace w_i \rbrace_{i=0}^d$. Then $TS$ is the transition matrix from $\lbrace u_i \rbrace_{i=0}^d$ to $\lbrace w_i \rbrace_{i=0}^d$. 

\begin{lemma}
\label{lem:chiconstruction}
Let $(X, \lbrace \theta_i \rbrace_{i=0}^d, \lbrace \theta^*_i \rbrace_{i=0}^d)$ denote a normalized west-south Vandermonde system in $\Mdf$ and let $\Phi$ denote the TH system constructed above Definition \ref{def:pi}. Then $X$ is the transition matrix of $\Phi$ from Definition \ref{def:transmatnorm}. 
\end{lemma}

\noindent {\it Proof:} Let $\mathcal P$ denote the transition matrix of $\Phi$ from Definition \ref{def:transmatnorm}. We show that $\mathcal P = X$. In what follows we refer to the construction of $\Phi$ above Definition \ref{def:pi}. By the construction of $A$ (resp. $A^*$) we find that $u_i \in E_iV$ (resp. $v_i \in E^*_iV$) for $0 \leq i \leq d$. Recall that $X$ is the transition matrix from $\lbrace u_i \rbrace_{i=0}^d$ to $\lbrace v_i \rbrace_{i=0}^d$. By these comments, Definition \ref{def:transmatnorm}, and the comment above this lemma, there exist invertible diagonal matrices $D_1, D_2 \in \Mdf$ such that $\mathcal P = D_1 X D_2$. The matrices $\mathcal P$ and $X$ are west normalized, meaning that $P_{i0} = X_{i0} = 1$ for $0 \leq i \leq d$.  The matrices $\mathcal P$ and $X$ are south normalized, meaning that $P_{di} = X_{di} = 1$ for $0 \leq i \leq d$. Evaluating the equation  $\mathcal P = D_1 X D_2$ using these comments, we find that $D_1$ is a nonzero scalar multiple of the identity and $D_2$ is the inverse of $D_1$. Therefore $\mathcal P = X$. \hfill $\Box$

\begin{theorem}
\label{thm:thsysvandbij}
The map $\rho$ from Definition \ref{def:rho} and the map $\chi$ from Definition \ref{def:pi} are inverses of each other. Moreover each of $\rho$, $\chi$ is bijective. 
\end{theorem}

\noindent {\it Proof:} We show that $\rho \circ \chi$ is the identity map on (\ref{eq:wsvand}) and 
$\chi \circ \rho$ is the identity map on (\ref{eq:thsystem}).
We first show that $\rho \circ \chi$ is the identity map on (\ref{eq:wsvand}). 
Let $(X, \lbrace \theta_i \rbrace_{i=0}^d, \lbrace \theta^*_i \rbrace_{i=0}^d)$ denote a normalized west-south Vandermonde system in $\Mdf$. Let $\Phi$ denote the corresponding TH system constructed above Definition \ref{def:pi}. The map $\chi$ sends $(X, \lbrace \theta_i \rbrace_{i=0}^d, \lbrace \theta^*_i \rbrace_{i=0}^d)$ to the isomorphism class of $\Phi$. Recall from above Definition \ref{def:pi} that $\Phi$ has eigenvalue sequence $\lbrace \theta_i \rbrace_{i=0}^d$ and dual eigenvalue sequence $\lbrace \theta_i^* \rbrace_{i=0}^d$. By Lemma \ref{lem:chiconstruction} $X$ is the transition matrix of $\Phi$ from Definition \ref{def:transmatnorm}. By these comments and Definition \ref{def:rho} the map $\rho$ sends the isomorphism class of $\Phi$ to $(X, \lbrace \theta_i \rbrace_{i=0}^d, \lbrace \theta^*_i \rbrace_{i=0}^d)$. Therefore $\rho \circ \chi$ is the identity map on (\ref{eq:wsvand}). 

Next we show that $\chi \circ \rho$ is the identity map on (\ref{eq:thsystem}). 
Let $\Phi$ denote a TH system over $\K$ of diameter $d$, with eigenvalue sequence $\lbrace \theta_i \rbrace_{i=0}^d$ and dual eigenvalue sequence $\lbrace \theta^*_i \rbrace_{i=0}^d$. Let $\mathcal P$ denote the transition matrix of $\Phi$ from Definition \ref{def:transmatnorm}. The map $\rho$ sends the isomorphism class of $\Phi$ to $(\mathcal P, \lbrace \theta_i \rbrace_{i=0}^d, \lbrace \theta^*_i \rbrace_{i=0}^d)$. The map $\chi$ sends  
$(\mathcal P, \lbrace \theta_i \rbrace_{i=0}^d, \lbrace \theta^*_i \rbrace_{i=0}^d)$ to the isomorphism class of $\Phi'$, where $\Phi'$ is the corresponding TH system constructed above Definition \ref{def:pi}. Recall from above Definition \ref{def:pi} that $\Phi'$ has eigenvalue sequence $\lbrace \theta_i \rbrace_{i=0}^d$ and dual eigenvalue sequence $\lbrace \theta_i^* \rbrace_{i=0}^d$. We show that $\Phi$ and $\Phi'$ are isomorphic. To do this we will invoke Lemma
\ref{lem:isomaltsystem}. 
Write $\Phi= (A;\lbrace E_i\rbrace_{i=0}^d;A^*;\lbrace E^*_i\rbrace_{i=0}^d)$ and $\Phi'=(A';\lbrace E_i' \rbrace_{i=0}^d;A^{*\prime};\lbrace E^{*\prime}_i\rbrace_{i=0}^d)$. Let $V$ (resp. $V'$) denote the vector space underlying $\Phi$ (resp. $\Phi'$). 
Let $0 \neq \xi^*_0 \in E^*_0V$ (resp. $0 \neq \xi^{*\prime}_0 \in E^{*\prime}_0V'$) and recall the $\Phi^*$-standard basis (resp. $\Phi^{*\prime}$-standard basis) $\lbrace E_i \xi^*_0 \rbrace_{i=0}^d$ (resp. $\lbrace E'_{i} \xi^{*\prime}_0 \rbrace_{i=0}^d$) for $V$ (resp. $V'$) from above Proposition \ref{prop:orthogstar}. 
Let $\Gamma: V \rightarrow V'$ denote the $\K$-vector space isomorphism which sends $E_i \xi^*_0$ to $E_i' \xi^{*\prime}_0$ for $0 \leq i \leq d$.
We show that
\begin{eqnarray}
\label{lab:isosys}
A' \Gamma = \Gamma A, \qquad A^{*\prime} \Gamma = \Gamma A^*, \qquad E_i' \Gamma = \Gamma E_i, \qquad  E^{*\prime}_i \Gamma = \Gamma E^*_i  \ \ \ (0 \leq i \leq d).
\end{eqnarray}
We first show that $E_i' \Gamma = \Gamma E_i$ for $0 \leq i \leq d$. Let $i$ be given. In order to show that $E_i' \Gamma = \Gamma E_i$, we show that $E_i' \Gamma$ and $\Gamma E_i$ agree at each vector in the $\Phi^*$-standard basis $\lbrace E_j \xi^*_0 \rbrace_{j=0}^d$. 
Observe that for $0 \leq j \leq d$, $E_i' \Gamma E_j \xi^*_0 =
E_i' E_j' \xi^{*\prime}_0 = \delta_{ij} E_i' \xi^{*\prime}_0$ and $\Gamma E_i E_j \xi^*_0 = \delta_{ij} \Gamma E_i \xi^*_0 = \delta_{ij} E_i' \xi^{*\prime}_0$. Thus $E_i' \Gamma = \Gamma E_i$.
Next we show that $A' \Gamma = \Gamma A$. 
Recall $A=\sum_{i=0}^d \theta_i E_i$. Observe that $A'=\sum_{i=0}^d \theta_i E_i'$ since $\Phi'$ has eigenvalue sequence $\lbrace \theta_i \rbrace_{i=0}^d$. 
By these comments $A' \Gamma = \Gamma A$. Next we show that $E^{*\prime}_i \Gamma = \Gamma E^*_i$ for $0 \leq i \leq d$. 
Let $P$ (resp. $P'$) denote the transition matrix of $\Phi$ (resp. $\Phi'$) from Definition \ref{def:transmat},
and let $L$ (resp. $L'$) denote the matrix associated with $\Phi$ (resp. $\Phi'$) from Definition \ref{def:diagli}. Observe that $L = L'$ by Definition \ref{def:kappas}, since $\Phi$ and $\Phi'$ have the same dual eigenvalue sequence $\lbrace \theta_i^* \rbrace_{i=0}^d$. By Lemma \ref{lem:chiconstruction} $\mathcal P$ is the transition matrix of $\Phi'$ from Definition \ref{def:transmatnorm}. By these comments and Definition \ref{def:transmatnorm}, we have $P = P'$.
Let $0 \neq \xi_0 \in E_0V$ (resp. $0 \neq \xi_0' \in E_0'V'$) such that $\xi^*_0=E^*_0\xi_0$ (resp. $\xi^{*\prime}_0=E^{*\prime}_0\xi_0'$). Recall the $\Phi$-standard basis (resp. $\Phi'$-standard basis) $\lbrace E^*_i \xi_0 \rbrace_{i=0}^d$ (resp. $\lbrace E^{*\prime}_i \xi_0' \rbrace_{i=0}^d$) for $V$ (resp. $V'$) from above (\ref{eq:defui}). By Definition \ref{def:transmat} and since $P = P'$, $P$ is the transition matrix from $\lbrace E_i \xi^*_0 \rbrace_{i=0}^d$ (resp. $\lbrace E_i' \xi^{*\prime}_0 \rbrace_{i=0}^d$) to $\lbrace E^*_i \xi_0 \rbrace_{i=0}^d$ (resp. $\lbrace E^{*\prime}_i \xi'_0 \rbrace_{i=0}^d$). We can now easily show that $E^{*\prime}_i \Gamma = \Gamma E^*_i$ for $0 \leq i \leq d$. Let $i$ be given. In order to show that $E_i^{*\prime} \Gamma = \Gamma E^*_i$ we show that $E_i^{*\prime} \Gamma$ and $\Gamma E^*_i$ agree at each vector in the $\Phi$-standard basis $\lbrace E^*_j \xi_0 \rbrace_{j=0}^d$. For $0 \leq j \leq d$,
\begin{eqnarray*}
E^{*\prime}_i \Gamma E^*_j \xi_0 = E^{*\prime}_i \Gamma \sum_{h=0}^d P_{hj} E_h \xi^*_0 = E^{*\prime}_i \sum_{h=0}^d P_{hj} E'_h \xi^{*\prime}_0 = E^{*\prime}_i E^{*\prime}_j \xi'_0  = \delta_{ij} E^{*\prime}_j \xi'_0, \\
\Gamma E^*_i E^*_j \xi_0 = \delta_{ij} \Gamma E^*_j \xi_0 = \delta_{ij} \Gamma \sum_{h=0}^d P_{hj} E_h \xi^*_0 = \delta_{ij} \sum_{h=0}^d P_{hj} E_h' \xi^{*\prime}_0 = \delta_{ij} E^{*\prime}_j \xi'_0.
\end{eqnarray*}
We have now shown that $E_i^{*\prime} \Gamma$ and $\Gamma E^*_i$ agree at each vector in the $\Phi$-standard basis $\lbrace E^*_j \xi_0 \rbrace_{j=0}^d$.
Therefore $E_i^{*\prime} \Gamma = \Gamma E^*_i$. Next we show that $A^{*\prime} \Gamma = \Gamma A^*$.  Recall $A^* =\sum_{i=0}^d \theta_i^* E_i^*$. Observe that $A^{*\prime}=\sum_{i=0}^d \theta_i^* E_i^{*\prime}$ since $\Phi'$ has dual eigenvalue sequence $\lbrace \theta_i^* \rbrace_{i=0}^d$. By these comments $A^{*\prime} \Gamma = \Gamma A^*$. 
We have now shown (\ref{lab:isosys}). Now $\Phi$ and $\Phi'$ are isomorphic in view of Lemma \ref{lem:isomaltsystem}. Therefore $\chi \circ \rho$ is the identity map on (\ref{eq:thsystem}). 
The result follows. \hfill $\Box$

\medskip

\noindent Combining Corollary \ref{cor:thsyspabij} and Theorem \ref{thm:thsysvandbij}, we get a bijection between any two of the following three sets: 

\begin{itemize}
 \item 
The set of isomorphism classes of TH systems over $\K$ of diameter $d$.
 \item
The set of normalized west-south Vandermonde systems in $\Mdf$. 
 \item
The set of parameter arrays over $\K$ of diameter $d$. 
\end{itemize}

\section{Reduced TH systems and Vandermonde matrices}
In the previous section we explained how double Vandermonde systems correspond with TH systems. In this section we turn our attention to double Vandermonde matrices and explain how these correspond with objects called reduced TH systems.

\begin{definition} 
\label{def:RHS}
\rm
A sequence $(\lbrace E_i\rbrace_{i=0}^d; \lbrace E^*_i\rbrace_{i=0}^d)$ is called a {\it reduced TH system} (or {\it RTH system}) on $V$ whenever there exist $A, A^* \in {\rm End}(V)$ such that $(A;\lbrace E_i\rbrace_{i=0}^d;A^*;\lbrace E^*_i\rbrace_{i=0}^d)$ is a TH system on $V$. Let $\Phi = (A;\lbrace E_i\rbrace_{i=0}^d;A^*;\lbrace E^*_i\rbrace_{i=0}^d)$ denote a TH system on $V$. Then $(\lbrace E_i\rbrace_{i=0}^d;\lbrace E^*_i\rbrace_{i=0}^d)$ is an RTH system on $V$, called the {\it reduction} of $\Phi$.
\end{definition} 

\begin{definition}
\label{def:isomredsystem}
\rm
Let $\Lambda=(\lbrace E_i\rbrace_{i=0}^d; \lbrace E^*_i\rbrace_{i=0}^d)$ denote an RTH system on $V$. Let $W$ denote a vector space over $\fld$ with dimension $d+1$, and let $\Omega=(\lbrace F_i \rbrace_{i=0}^d; \lbrace F^*_i\rbrace_{i=0}^d)$ denote an RTH system on $W$. By an {\it isomorphism of RTH systems} from $\Lambda$ to $\Omega$ we mean a $\K$-algebra isomorphism $\gamma:{\rm End}(V) \rightarrow {\rm End}(W)$ such that $F_i = E_i^\gamma$ and $F^*_i = E^{*\gamma}_i$ for $0 \leq i \leq d$. We say that the RTH systems $\Lambda$ and $\Omega$ are {\it isomorphic} whenever there exists an isomorphism of RTH systems from $\Lambda$ to $\Omega$. 
\end{definition}

\begin{proposition}
\label{prop:reduction1}
Let $\Phi$ and $\Phi'$ denote TH systems over $\K$. Then the following {\rm (i)}, {\rm (ii)} are equivalent. 
\begin{enumerate}
 \item[\rm (i)] 
The reduction of $\Phi$ is isomorphic to the reduction of $\Phi'$. 
 \item[\rm (ii)]
$\Phi$ is affine isomorphic to $\Phi'$. 
\end{enumerate}
\end{proposition}

\noindent{\it Proof:} (i) $\Rightarrow$ (ii) Let $P$ (resp. $P'$) denote the transition matrix of $\Phi$ (resp. $\Phi'$) from Definition \ref{def:transmat}. From its definition, we see that $P$ (resp. $P'$) is determined by the primitive idempotents of $\Phi$ (resp. $\Phi'$). Hence by our assumption $P = P'$. Let $\lbrace \theta_i \rbrace_{i=0}^d$ (resp. $\lbrace \theta_i' \rbrace_{i=0}^d$) denote the eigenvalue sequence of $\Phi$ (resp. $\Phi'$). By Proposition \ref{prop:pdoublevandwest} $P$ and $\lbrace \theta_i \rbrace_{i=0}^d$ are compatible. Similarly $P'$ and $\lbrace \theta_i' \rbrace_{i=0}^d$ are compatible. Since $P = P'$, we conclude that $P$ is compatible with each of $\lbrace \theta_i \rbrace_{i=0}^d$ and $\lbrace \theta_i' \rbrace_{i=0}^d$.
Hence by Lemma \ref{lem:thetarelated} there exist $\alpha, \beta \in \K$ with $\alpha \neq 0$ such that $\theta_i' = \alpha \theta_i + \beta$ for $0 \leq i \leq d$.
Let $\lbrace \theta_i^* \rbrace_{i=0}^d$ (resp. $\lbrace \theta_i^{*\prime} \rbrace_{i=0}^d$) denote the dual eigenvalue sequence of $\Phi$ (resp. $\Phi'$). By a similar argument,
there exist $\alpha^*, \beta^* \in \K$ with $\alpha^* \neq 0$ such that $\theta_i^{*\prime} = \alpha^* \theta_i^* + \beta^*$ for $0 \leq i \leq d$. It follows that $\Phi$ is affine isomorphic to $\Phi'$. \\
(ii) $\Rightarrow$ (i) Clear.  \hfill $\Box$

\begin{corollary}
\label{cor:isouniquesystem2}
Let $\Lambda$ and $\Omega$ denote isomorphic RTH systems over $\fld$. Then the isomorphism of RTH systems from $\Lambda$ to $\Omega$ is unique. 
\end{corollary}

\noindent {\it Proof:} Let $\gamma$ and $\gamma'$ denote isomorphisms of RTH systems from $\Lambda$ to $\Omega$. We show that $\gamma = \gamma'$. 
Let $\Phi$ (resp. $\Psi$) denote a TH system over $\K$ whose reduction is $\Lambda$ (resp. $\Omega$). By Proposition \ref{prop:reduction1} $\Phi$ is affine isomorphic to $\Psi$. In other words, $\Phi$ is isomorphic to an affine transformation $\Psi'$ of $\Psi$. By construction $\Phi$ and $\Psi'$ have the same eigenvalue sequence and dual eigenvalue sequence. By this and the comment (iv) above (\ref{eq:defEi}), we find that each of $\gamma$ and $\gamma'$ is an isomorphism of TH systems from $\Phi$ to $\Psi'$. Now $\gamma = \gamma'$ in view of Lemma \ref{lem:isouniquesystem}. The result follows. \hfill $\Box$

\medskip

\noindent We now give a correspondence between TH systems and reduced TH systems. 

\begin{corollary}
\label{cor:one}
The map which sends a TH system to its reduction induces a bijection from the set of affine isomorphism classes of TH systems over $\K$ to the set of isomorphism classes of RTH systems over $\K$. 
\end{corollary}

\noindent{\it Proof:} Immediate from Proposition \ref{prop:reduction1}. \hfill $\Box$

\medskip

\noindent We now turn our attention to double Vandermonde systems and double Vandermonde matrices. 

\begin{lemma}         
\label{lem:affinevand}  
Let $\Omega = (X, \lbrace \theta_i \rbrace_{i=0}^d, \lbrace \theta^*_i \rbrace_{i=0}^d)$ denote a normalized west-south Vandermonde system in $\Mdf$. Let $\alpha, \beta, \alpha^*, \beta^*$ denote scalars in $\K$ with $\alpha, \alpha^*$ nonzero. Then the sequence
\begin{equation}         
\label{eq:affOmega}
(X, \lbrace \alpha \theta_i + \beta \rbrace_{i=0}^d; \lbrace \alpha^* \theta_i^* + \beta^* \rbrace_{i=0}^d)
\end{equation}
is a normalized west-south Vandermonde system in $\Mdf$.
\end{lemma}

\noindent{\it Proof:} Routine by Lemma \ref{lem:thetarelated}. \hfill $\Box$

\begin{definition}            
\label{def:affvand}        
\rm
Referring to Lemma \ref{lem:affinevand}, we call (\ref{eq:affOmega}) the {\em affine transformation of $\Omega$ associated with $\alpha, \beta, \alpha^*, \beta^*$}.
\end{definition}

\begin{definition}          
\label{def:affisovand}        
\rm 
Let $\Omega$ and $\Omega'$ denote normalized west-south Vandermonde systems in $\Mdf$. We say that $\Omega$ and $\Omega'$ are {\em affine related} whenever $\Omega$ is an affine transformation of $\Omega'$. Observe that the affine relation is an equivalence relation.
\end{definition}

\begin{lemma}
\label{lem:reduction2}
Let $\Omega = (X, \lbrace \theta_i \rbrace_{i=0}^d, \lbrace \theta^*_i \rbrace_{i=0}^d)$ and $\Omega' = (X', \lbrace \theta_i' \rbrace_{i=0}^d, \lbrace \theta^{*\prime}_i \rbrace_{i=0}^d)$ denote normalized west-south Vandermonde systems in $\Mdf$. Then the following {\rm (i)}, {\rm (ii)} are equivalent.
\begin{enumerate}
 \item[\rm (i)] 
$X = X'$.
 \item[\rm (ii)]
$\Omega$ is affine related to $\Omega'$.
\end{enumerate}
\end{lemma}

\noindent {\it Proof:} (i) $\Rightarrow$ (ii) Immediate from Lemma \ref{lem:thetarelated}. \\
(ii) $\Rightarrow$ (i) Clear.  \hfill $\Box$

\medskip

\noindent We now give a correspondence between normalized double Vandermonde systems and normalized double Vandermonde matrices. 

\begin{corollary}
\label{cor:two}
The map which sends a normalized west-south Vandermonde system \\
$(X, \lbrace \theta_i \rbrace_{i=0}^d, \lbrace \theta^*_i \rbrace_{i=0}^d)$ to the matrix $X$ induces a bijection from the set of affine classes of normalized west-south Vandermonde systems in $\Mdf$ to the set of normalized west-south Vandermonde matrices in $\Mdf$. 
\end{corollary}

\noindent{\it Proof:} Immediate from Lemma \ref{lem:reduction2}. \hfill $\Box$

\medskip

\noindent Next we give a correspondence between affine isomorphism classes of TH systems and affine classes of normalized double Vandermonde systems. Recall the map $\rho$ from Definition \ref{def:rho}. 

\begin{lemma}
\label{lem:rhoaffine}
Let $\Phi=(A;\lbrace E_i\rbrace_{i=0}^d;A^*;\lbrace E^*_i\rbrace_{i=0}^d)$ denote a TH system over $\K$. Let \\ $(X, \lbrace \theta_i \rbrace_{i=0}^d, \lbrace \theta^*_i \rbrace_{i=0}^d)$ denote the image under $\rho$ of the isomorphism class of $\Phi$. Let $\alpha, \beta, \alpha^*, \beta^*$ denote scalars in $\K$ with $\alpha, \alpha^*$ nonzero and consider the TH system (\ref{eq:affPhi}). Then $\rho$ sends the isomorphism class of (\ref{eq:affPhi}) to $(X, \lbrace \alpha \theta_i + \beta \rbrace_{i=0}^d, \lbrace \alpha^* \theta_i^* + \beta^* \rbrace_{i=0}^d)$. 
\end{lemma}

\noindent{\it Proof:} Immediate from Lemma \ref{lem:affinepa1}. \hfill $\Box$

\begin{corollary}
\label{cor:three}
The bijection $\rho$ from  Definition \ref{def:rho} induces a bijection from the set of affine isomorphism classes of TH systems over $\K$ of diameter $d$, to the set of affine classes of normalized west-south Vandermonde systems in $\Mdf$.  
\end{corollary}

\noindent{\it Proof:} Immediate from Lemma \ref{lem:rhoaffine}. \hfill $\Box$

\medskip

\noindent We now bring the parameter arrays into the discussion. 

\begin{definition}
\label{def:eqclasspa}
\rm
We define a binary relation on the set of parameter arrays over $\fld$ of diameter $d$. We do this as follows. Let $p=(\lbrace \theta_i \rbrace_{i=0}^d, \lbrace \theta^*_i \rbrace_{i=0}^d, \lbrace \phi_i \rbrace_{i=1}^d)$ and $p'=(\lbrace \theta_i' \rbrace_{i=0}^d, \lbrace \theta^{*\prime}_i \rbrace_{i=0}^d, \lbrace \phi_i' \rbrace_{i=1}^d)$ denote parameter arrays over $\fld$ of diameter $d$. We say that $p$ and $p'$ are {\it affine related} whenever there exist scalars $\alpha, \beta, \alpha^*, \beta^*$ in $\fld$ with $\alpha, \alpha^*$ nonzero such that the following {\rm (i)--(iii)} hold.
\begin{enumerate}
\item[\rm (i)] 
$ \theta_i' = \alpha \theta_i + \beta \qquad \ \ (0 \leq i \leq d).$
\item[\rm (ii)]
$\theta^{*\prime}_i = \alpha^* \theta^*_i + \beta^* \qquad (0 \leq i \leq d).$
\item[\rm (iii)] 
$\phi_i' = \alpha \alpha^* \phi_i \qquad \ \  (1 \leq i \leq d).$
\end{enumerate}
Observe that the affine relation is an equivalence relation. By a {\it reduced parameter array} we mean an equivalence class of this relation. 
\end{definition}

\begin{corollary}
\label{cor:four}
The bijection from Corollary \ref{cor:thsyspabij} induces a bijection from the set of affine isomorphism classes of TH systems over $\K$ of diameter $d$, to the set of reduced parameter arrays over $\K$ of diameter $d$.
\end{corollary}

\noindent{\it Proof:} Immediate from Lemma \ref{lem:affinepa1}. \hfill $\Box$

\medskip

Combining Corollaries \ref{cor:one}, \ref{cor:two}, \ref{cor:three}, \ref{cor:four}, we get a bijection between any two of the following five sets:

\begin{itemize}
\item
The set of affine isomorphism classes of TH systems over $\K$ of diameter $d$. 
\item 
The set of isomorphism classes of RTH systems over $\K$ of diameter $d$.
\item
The set of affine classes of normalized west-south Vandermonde systems in $\Mdf$.
 \item
The set of normalized west-south Vandermonde matrices in $\Mdf$. 
 \item
The set of reduced parameter arrays over $\K$ of diameter $d$.
\end{itemize}

\section{Acknowledgement}
This paper was written while the author was a graduate student at the University of Wisconsin-Madison. The author would like to thank his advisor Paul Terwilliger for his many valuable ideas and suggestions.

\bigskip

\noindent Ali Godjali \hfil\break
\noindent Department of Mathematics \hfil\break
\noindent University of Wisconsin \hfil\break
\noindent Van Vleck Hall \hfil\break
\noindent 480 Lincoln Drive \hfil\break
\noindent Madison, WI 53706-1388 USA \hfil\break
\noindent email: {\tt godjali@math.wisc.edu }\hfil\break

\end{document}